\documentclass[12pt]{article}
\usepackage{amsmath, amsthm, amsfonts, amssymb, color}
\usepackage{mathrsfs}
\usepackage{amsfonts, amsmath}
\usepackage{amsmath,amstext,amsthm,amssymb,amsxtra}
\usepackage{txfonts}
\usepackage{tikz}

\usepackage{latexsym,amssymb,amsmath,amsfonts,graphicx,subfigure,epsfig,psfrag,CJK}
\usepackage{float}
\usepackage{psfrag}

\allowdisplaybreaks[4]
\usepackage[colorlinks, citecolor=blue,hypertexnames=false]{hyperref}
\usepackage{pgf,tikz}

\newtheorem{theo}{Theorem}[section]
\newtheorem{lemm}[theo]{Lemma}

\newtheorem{rema}[theo]{Remark}

\numberwithin{equation}{section}

\allowdisplaybreaks

\def\th2{\frac{\theta}{2}}

\title{\bf Four limit cycles in a predator-prey system of Leslie type with generalized Holling type III functional response}
\author{Yanfei Dai $^{1}$,  \quad Yulin Zhao $^{2}$\thanks{Corresponding author.  E-mail:
mcszyl@mail.sysu.edu.cn.}, \quad Bo Sang $^{3}$}
\date{\it\footnotesize $^{1}$ School of Mathematics, Sun Yat-sen University, Guangzhou, 510275, P. R. China\\
$^{2}$School of Mathematics (Zhuhai), Sun Yat-sen University, Zhuhai, 519082, P. R. China\\
$^{3}$ School of Mathematical Sciences, Liaocheng University, Liaocheng 252059, P. R. China}
\begin{document}
\maketitle 

\begin{abstract}
This paper, as a complement to the works by Hsu et al [SIAM. J. Appl. Math. 55 (1995)] and Huang et al [J. Differential Equations 257 (2014)], aims to examine the Hopf bifurcation and global
dynamics of a predator-prey model of Leslie type with generalized Holling type III functional response for the two cases: (A) when it has a unique non-degenerate positive equilibrium; (B)
when it has three distinct positive equilibria.
For each case, the type and stability of each equilibrium, Hopf bifurcation at each weak focus, and the number and distribution of limit cycles in the first quadrant are studied.
It is shown that every equilibrium is not a center.
For the case (A), $i$ limit cycle(s) can appear near the unique positive equilibrium, $i=1, \cdots, 4$.
For $i=3$ or $4$, the model has two stable limit cycles, which gives a positive answer to the open problem proposed by Coleman [Differential equations model,1983]: finding at least two ecologically stable cycles in any natural (or laboratory) predator-prey system.
For the case (B), one positive equilibrium is a saddle and the others are both anti-saddle.
If one of the two anti-saddles is weak focus and the other is not, then the order of the weak focus is at most $3$.
If both anti-saddles are weak foci, then they are unstable weak foci of order one.
Moreover, one limit cycle can bifurcate from each of them simultaneously.
Numerical simulations demonstrate  that there is also a big stable limit cycle enclosing these two limit cycles.
Our results indicate that the maximum number of limit cycles in the model of this kind is at least $4$, which improves the preceding results that this number is at least $2$.
\end{abstract}
\maketitle
\rule[1.5mm]{\textwidth}{0.1pt}

\section{Introduction and statement of the main results}

In this work we consider the predator-prey system of Leslie type with generalized Holling type III or sigmoidal functional response \cite{Huang2014},
\begin{equation}\label{model}
\left\{
\begin{aligned}
&\dot{x}=r x \left(1-\dfrac{x}{K}\right)-\dfrac{m x^2 y}{a x^2+b x+1}, \\
&\dot{y}=s y \left(1-\dfrac{y}{h x}\right)
\end{aligned}
\right.
\end{equation}
with $(x, y)\in \mathcal{A}=\{(x, y)\left|x>0, y\geq 0 \right.\}$ and $r,\ K,\ m,\ a,\ s,\ h>0$, $b>-2\sqrt{a}$.

In system \eqref{model} $x(t)$ and $y(t)$ represent the population densities of the prey and predator at time $t>0$, respectively; the effect of the predation is given by the function
\begin{equation*}
p(x)=\frac{m x^2}{a x^2+b x+1},
\end{equation*}
and is called  the generalized Holling type III or sigmoidal functional response \cite{Bazykin1998}. The parameters  $r$ and $s$ are the intrinsic growth rates or biotic potential of the prey and
predator, respectively,  $K$ is the prey environment carrying capacity, and   $h$ is a measure of the food quality of the prey for conversion into
predator births.

Here the predator growth equation is of Leslie form originated by Leslie \cite{Leslie1948}, but the conventional environmental carrying capacity $h x$ is proportional to prey abundance $x$ \cite{May1973}.

The system \eqref{model} has been investigated by Hsu et al \cite{Hsu1995} and Huang et al \cite{Huang2014}. It has been conjectured that for predator-prey systems with a unique positive
equilibrium, local and global stability are equivalent \cite{Saez1999}.
When $b>0$, Hsu et al \cite{Hsu1995} prove that it is true under certain conditions.
When $b>-2\sqrt{a}$ (so that $a x^2+b x+1>0$ and hence $p(x)>0$ for all $x>0$), Huang et al \cite{Huang2014} show that it is not true by analysing the nonlinear dynamics
of system \eqref{model} when it has at least one degenerate positive equilibrium.

The predator-prey systems have been widely studied, such as investigating the
stability and bifurcations, proving the global stability of the unique positive equilibrium, and studying the uniqueness or nonexistence of limit cycles, see in \cite{May1973,Bazykin1998,Hsu1995} and the references therein.
However, it's not an easy task to study the number and distribution of limit cycles in a given predation model. This problem is related to Hilbert's 16th Problem \cite{Gaiko}, and it is a question that has remained unsolved for the predation model \cite{Eduardo2011}.

For some two-dimensional predator-prey systems, many authors proved that there can exist at least two limit cycles, see  \cite{Kuang1988,Saez1999,Ruan2001,Wrzosek1990,Zhu2002}, etc.
In \cite{Aguirre2009}, Aguirre et al showed that there can exist three limit cycles including only one stable cycle.
In \cite{Lu2001,Lloyd1996} the existence of some Kolmogorov type systems with at least two stable limit cycles surrounding the singularity in the positive quadrant is given. We note that the systems are not predator-prey systems. So far, whether there are more than $3$ limit cycles in predator-prey systems is still unknown.
The study of the number and distribution of limit cycles in a given predation model may be the most difficult part. This is principally because these systems are generally not polynomial differential equations, and the expressions of some positive equilibria of such systems are so complicated and even not formulated explicitly that one can not study them.
In this paper, we will propose some available methods for studying the existence, stability, number and distribution of limit cycles in a given predator-prey systems.

This paper aims to study the system \eqref{model} when it has no degenerate positive equilibrium. If $E(x_{*}, y_{*})$ is a positive equilibrium of system \eqref{model}, then $y_{*}=hx_{*}$ and $x_{*}$ is a root of the equation
\begin{equation}\label{dai1.2}
a x^3+\left(\frac{K m h}{r}+b-K a\right) x^2+(1-K b) x-K=0
\end{equation}
in the interval $(0, K)$.
We denote the determinant and trace of the Jacobian matrix of system \eqref{model} at $E$ by $\textrm{det}\left(J(E)\right)$ and $\textrm{tr}\left(J(E)\right)$, respectively.
If $\textrm{det}\left(J(E)\right)\neq0$, then  $E$ is called an elementary equilibrium,  otherwise  it is a degenerate equilibrium.
Specially, $E$ is called a hyperbolic saddle  if $\textrm{det}\left(J(E)\right)<0$ and called \emph{center or focus type} if $\textrm{det}\left(J(E)\right)>0$ and $\textrm{tr}\left(J(E)\right)=0$, respectively.

The number of positive  equilibria of
system \eqref{model} is determined by the number of roots of Eq. \eqref{dai1.2} in the interval $(0, K)$.
Note that  Eq. \eqref{dai1.2} can have one, two or
three positive roots in the interval $(0, K)$. Correspondingly, system \eqref{model} can have one, two, or three positive equilibria.
For system \eqref{model}, it can be inferred from Hsu and Huang in \cite{Hsu1995} that if $x_{*}$ is a multiple root of Eq. \eqref{dai1.2} in the interval $(0, K)$, then $E(x_{*}, y_{*})$ must be a degenerate equilibrium of system \eqref{model}.
Let
\begin{equation*}
\Delta=-4A^3+\left[27 a^2 K+9 a (1-K b) \left(\frac{K m h}{r}+b-K a\right)-
2 \left(\frac{K m h}{r}+b-K a \right)^3\right]^2,
\end{equation*}
where $A=\left(K m h/r+b-K a\right)^2 +3a(Kb-1)$.
Then the case $\Delta=0$ corresponds to the case that system has a degenerate positive equilibrium, which has been studied by Huang et al \cite{Huang2014}.
Hence we only need to consider the case $\Delta\neq 0$, i.e., system \eqref{model} has either a unique non-degenerate positive equilibrium or three distinct positive equilibria.

First, we give the main results about the system with a unique non-degenerate positive equilibrium, which are stated by the following two theorems.
\begin{theo}\label{T3.2}
Suppose $\Delta>0$, then system \eqref{model} has a unique positive equilibrium, which is an elementary and anti-saddle equilibrium.
If the unique positive equilibrium is unstable, then there exists at least one stable limit cycle in the first quadrant.
\end{theo}
\begin{theo}\label{T4.1}
If  $\Delta>0$ and the unique positive equilibrium is center or focus type, then
\begin{itemize}
\item[(1)] it is not a center;
\item[(2)] it is a weak focus of order at most $4$ and there exists a unique class of  parameter values  such that its order is $4$;
\item[(3)] there exist some parameter values  such that system \eqref{model} has $i$ small limit cycles around it, for each $i=1, \cdots, 4$.
\end{itemize}
\end{theo}

It's worth to point out that the unique class of parameter values in the statement (2) of Theorem \ref{T4.1} satisfies $b>0$, see the proof of this theorem in Section $3$. Therefore, for the case $b\leq 0$, if the unique non-degenerate positive equilibrium of system \eqref{model} is center or focus type, then it is a weak focus of order at most $3$, implying at most $3$ limit cycles can bifurcate from it.

The phenomenon that one limit cycle can appear around the unique positive equilibrium was also proved by Hsu et al \cite{Hsu1995}  and  Huang et al \cite{Huang2014}.
The phenomenon that two limit cycles can appear near the unique positive equilibrium
was also observed by Huang et al \cite{Huang2014} through subcritical Hopf bifurcation and numerical simulations.
However, the phenomena that three or four limit cycles can appear near the unique positive equilibrium, observed in the present paper, have not been found by other authors.
Furthermore, the coexistence of four limit cycles including two stable cycles or three limit cycles including two stable cycles, which is not yet reported in the predation model, gives a positive answer to the open problem proposed by Coleman \cite{Coleman1983}, suggesting: ``Find a predator-prey or other interacting system in nature, or construct one in the laboratory, with at least two ecologically stable cycles''.

Second, we will consider the system with three different positive equilibria. The main results are given by the following two theorems.
\begin{theo}\label{T4.3}
Suppose  $\Delta<0$ and $-2\sqrt{a}<b<K a-K m h/r$, then system \eqref{model} has  three different positive equilibria, of which one is a saddle and the others are anti-saddles. If one of the two positive anti-saddles is center or focus type and the other is not, then
\begin{itemize}
\item[(1)] the center or focus type equilibrium is not a center. It is a weak focus of order at most $3$;
\item[(2)] two limit cycles can bifurcate from the weak focus.
\end{itemize}
\end{theo}

\begin{rema}\label{rem3}
When system has three distinct positive equilibria, two of them are anti-saddles.
The phenomenon that one limit cycle can bifurcate from one of the two positive anti-saddles was also observed by Huang et al \cite{Huang2014}.
However, the phenomenon that two limit cycles can bifurcate from one of the two positive anti-saddles, observed in the present paper, has
not been found by other authors. Unfortunately, we can't prove whether there exist three limit cycles surrounding one of the two positive anti-saddles.
\end{rema}
\begin{figure}\label{Fig1.1}
\begin{center}
{\includegraphics[height=7cm,width=9cm]{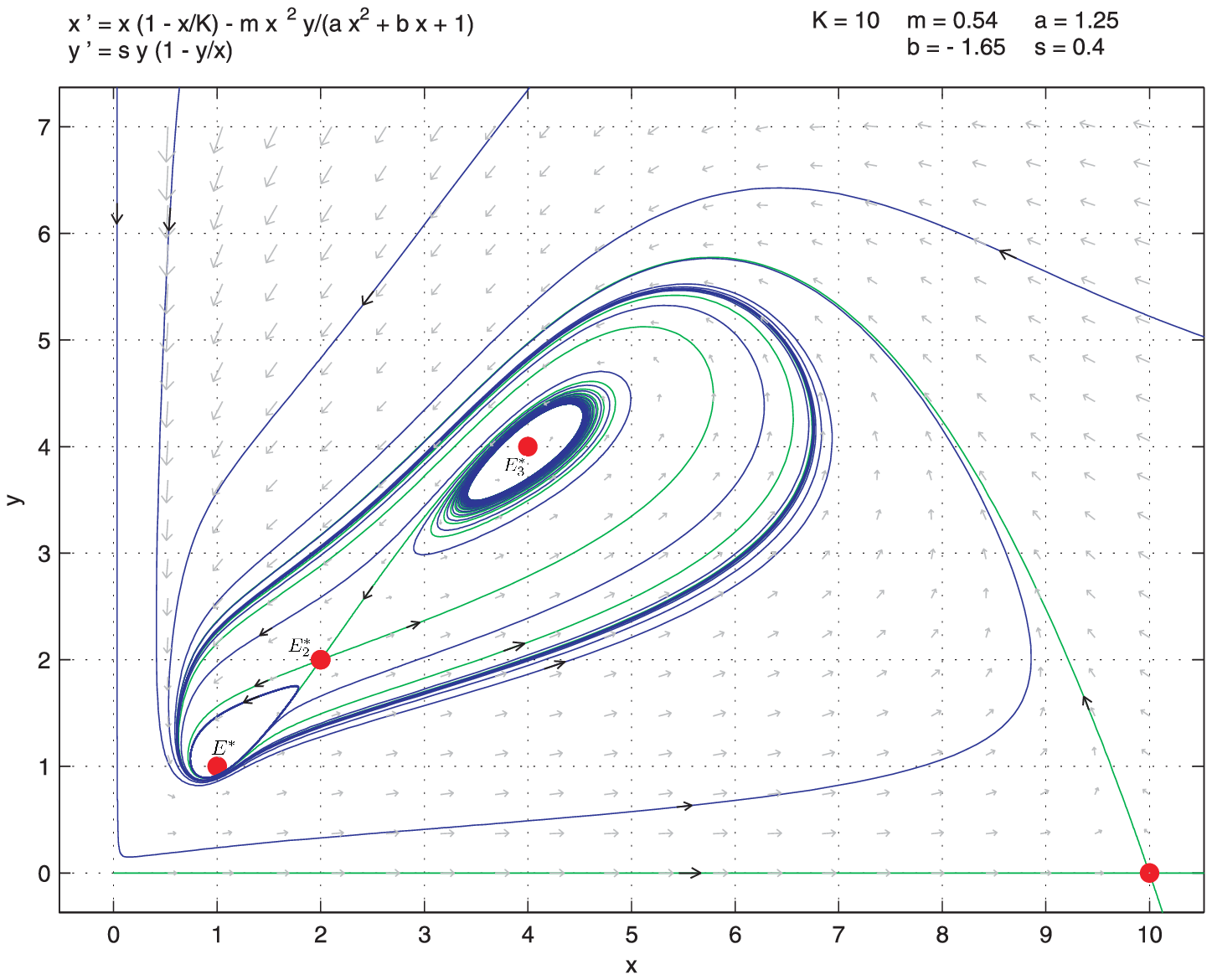}}
\end{center}
{\small {\bf Fig. 1.1.} Coexistence of one big stable limit cycle enclosing two unstable limit cycles in system \eqref{model}. Here, we fix $r=1$, $K=10$, $m=0.54$, $a=1.25$, $b=-1.65$, $h=1$ and $s=0.4$.}
\end{figure}

\begin{theo}\label{T4.4}
Suppose $\Delta<0$ and $-2\sqrt{a}<b<K a-K m h/r$, i.e., system \eqref{model} has one  positive saddle and two positive anti-saddles. If the two positive anti-saddles are center or focus type, then
\begin{itemize}
\item[(1)] they are unstable weak foci of order one;
\item[(2)] one limit cycle can bifurcate from each of them simultaneously.
\end{itemize}
\end{theo}

When system has three distinct positive equilibria, of which two are anti-saddles, one limit cycle can bifurcate from each of the two positive anti-saddles simultaneously (see Fig. 1.1). This has not been observed  by Hsu et al \cite{Hsu1995} and
Huang et al. \cite{Huang2014}. Numerical simulations show that there is also a big stable limit cycle enclosing these two limit cycles (see Fig. 1.1), which is a new phenomenon observed in the present paper. Unfortunately, we can't prove it.

The rest of this paper is organized as follows. Section $2$ is devoted to some  preliminary results including model reduction, analysis of equilibria of an equivalent polynomial differential system of \eqref{model} and  computation of Lyapunov constants at positive equilibrium. In Section $3$,  we give the proof of our main results.
The paper ends with a discussion.

\section{Some preliminary results}

The preliminary results provided in this section are helpful for the proof of our results.

\subsection{Model reduction}

In order to calculate the Liapunov constants at the nonhyperbolic focus and
study the Hopf bifurcation in the parameter space in a simpler way, it is necessary to reduce system \eqref{model} to a polynomial differential system.

Without loss of generality, assume that $E(x_{*}, y_{*})$ is an arbitrary positive equilibrium of system \eqref{model}.
Let us consider the change of variables
\begin{equation}\label{change1}
u=\dfrac{x}{x_{*}},\,\,\,\, v=\dfrac{y}{y_{*}},\,\,\,\, \tau=r t,
\end{equation}
and the parameter transformation given by $\zeta: \mathbb{R}_+^6\times \mathbb{R}\rightarrow \mathbb{R}_+^6\times \mathbb{R}$,
\begin{equation*}
\zeta (r, K, m, a, s, h, b)=\left(r, \frac{K}{x_{*}}, \frac{m x_{*}y_{*}}{r}, a x_{*}^2, \frac{s}{r}, \frac{h x_{*}}{y_{*}}, b x_{*}\right):=(r, K, m, a, s, h, b)
\end{equation*}
with Jacobian $\textrm{det}D \zeta (r, K, m, a, s, h, b)=\frac{x_{*}^4}{r^2}>0$, which implies that $\zeta$ is invertible.
For simplicity, let us rename the new variables $u\rightarrow x$, $v\rightarrow y$,  $\tau\rightarrow t$. Then in the new variables, system \eqref{model} is reduced to the following system
\begin{equation}\label{norm1}
\left\{
\begin{aligned}
&\dot{x}=x\left(1-\dfrac{x}{K}\right)-\dfrac{m x^2 y}{a x^2+b x+1}, \\
&\dot{y}=s y\left(1-\dfrac{y}{h x}\right).
\end{aligned}
\right.
\end{equation}
Noting that system \eqref{norm1} has an equilibrium at $E^{*}(1, 1)$, we have $m=(a+b+1) (1-1/K)$ and  $h=1$. Furthermore, $K>x_{*}$ becomes $K>1$ and $a>0$, $b>-2\sqrt{a}$,
$s>0$ remain the same, respectively, under the projection $\zeta$.
That is to say, system \eqref{norm1} can be written as
\begin{equation}\label{norm2}
\left\{
\begin{aligned}
&\dot{x}=x\left(1-\dfrac{x}{K}\right)-
\dfrac{(a+b+1)(K-1)x^2y}
{K(ax^2+bx+1)}, \\
&\dot{y}=sy\left(1-\dfrac{y}{x}\right)
\end{aligned}
\right.
\end{equation}
with $(x, y)\in \mathcal{A}$ and the new vector of parameters $\eta=(s, K, a, b)\in \mathcal{D}_0$ is given by the natural projection, where
\begin{equation*}
\mathcal{D}_0=\left\{\eta\in \mathbb{R}_+^3\times \mathbb{R}\left|s>0,\, K>1, \, a>0,\,  b>-2\sqrt{a}\right.\right\}.
\end{equation*}
Since the transformation \eqref{change1} is a linear sign-reserving transformation,  system \eqref{norm2} and system \eqref{model} have the same qualitative property.
By scaling
\begin{equation*}
d\tau=\dfrac{dt}{K x (a x^2+b x+1)},
\end{equation*}
(we will still use $t$ to denote $\tau$ for ease of notation) system \eqref{norm2} can be transformed into  the following quintic polynomial differential system
\begin{equation}\label{norm}
\left\{
\begin{aligned}
&\dot{x}=(ax^2+bx+1)\left(K-x\right)x^2-(K-1)(a+b+1)x^3y, \\
&\dot{y}=K s y\left(x-y\right)(ax^2+bx+1),
\end{aligned}
\right.
\end{equation}
with $(x, y)\in \bar{\mathcal{A}}=\{(x, y)\left|x, y\geq 0 \right.\}$ and $\eta\in \mathcal{D}_0$.
Notice that system \eqref{norm} has the same topological structure as system \eqref{norm2} in $\mathcal{A}$ because $K x (ax^2+bx+1)>0$ for all $x>0$.
In other words, system \eqref{norm} is equivalent to system \eqref{model} in $\mathcal{A}$. Hence we only need to study the polynomial system \eqref{norm} in the region $\bar{\mathcal{A}}$ with $\eta\in \mathcal{D}_0$.

\subsection{Analysis of equilibria  of  system \eqref{norm}}

Obviously, the $x$-axis, $y$-axis and the interior of $\bar{\mathcal{A}}$ are all invariant under system \eqref{norm}.
It's  standard  to  show  that  all solutions of  \eqref{norm} with positive initial values are positive and bounded, and  will eventually tend into the region $\Omega=\{(x(t),  y(t)): 0 \leq x(t)\leq K,\, 0\leq y(t)\leq K\}$. Therefore, $\Omega$ is a positive invariant set of system \eqref{norm} and the limit cycle of  \eqref{norm}, if it exists, must be inside $\Omega$.

Notice that system \eqref{norm} always has a boundary equilibrium $E_0=(K, 0)$ for all parameters, which is always a hyperbolic saddle and divides the positive $x$-axis into two parts which are two stable manifolds of $E_0$. Furthermore, there exists a unique unstable manifold of $E_0$ in the interior of  $\bar{\mathcal{A}}$.

The origin $O(0, 0)$ is an isolated critical point of higher order of system \eqref{norm}.
By introducing the polar coordinates $x=r\cos\theta$, $y=r\sin\theta$, we get the characteristic equation of system \eqref{norm} as follows.
\begin{equation}\label{dai3.1}
\begin{aligned}
G(\theta)=K \sin\theta\cos\theta[(s-1)\cos\theta-s\sin\theta]=0.
\end{aligned}
\end{equation}
Noting that $G(\theta^{*})=0$ is a necessary condition for $\theta=\theta^{*}$ to be a characteristic direction \cite{ZhangZF},  there are at most $3$ possible directions $\theta=0,  \pi/2, \theta_3$, where
$\theta_3=\arctan(1-1/s)$, along which an orbit of system \eqref{norm} may approach the origin in the first quadrant.
Let  $S^+_{\delta} (O)=\{(r, \theta):0<r<\delta,\
0\leq\theta\leq\pi/2\}$, where $0<\delta\ll 1$.
By Theorems 3.4, 3.7 and 3.8 in \cite{ZhangZF}, we get the following results.
\begin{lemm}\label{T3.1}  For system \eqref{norm}, the origin $O$ is an isolated critical point.
\begin{itemize}
\item[(1)] If  $0<s<1$, then
\begin{itemize}
\item[(a)] the positive $x$-axis to the left of $E_0$ is a unique orbit of system \eqref{norm} tending to $O$ along $\theta=0$ as  $t\rightarrow-\infty$;
\item[(b)] the positive $y$-axis is a unique orbit of system \eqref{norm} tending to $O$ along $\theta=\pi/2$ as $t\rightarrow+\infty$.
\end{itemize}
\item[(2)] If $s=1$, then
\begin{itemize}
\item[(a)] there is an infinite number of orbits of system \eqref{norm} in $S^+_{\delta} (O)$ tending to $O$ along $\theta=0$ as $t\rightarrow-\infty$;
\item[(b)]  the positive $y$-axis is a unique orbit of system \eqref{norm} tending to $O$ along $\theta=\pi/2$ as $t\rightarrow+\infty$.
\end{itemize}
\item[(3)] If $s>1$, then
\begin{itemize}
\item[(a)]  there is an infinite number of orbits of system \eqref{norm} in $S^+_{\delta} (O)$ tending to $O$ along  $\theta=0$ as $t\rightarrow-\infty$;
\item[(b)]  the positive $y$-axis is a unique orbit of system \eqref{norm} tending to $O$ along $\theta=\pi/2$ as
  $t\rightarrow+\infty$;
\item[(c)] there is a unique orbit of system \eqref{norm} tending to $O$ along $\theta=\theta_3$ as $t\rightarrow-\infty$. And this orbit is a separatrix that divides $S^+_{\delta} (O)$
    into two parts: one is a hyperbolic sector, the other is a repulsing parabolic sector.
\end{itemize}
\end{itemize}

The phase portrait of system \eqref{norm} near $O$ is shown in Fig. 2.1 for all cases.
\end{lemm}
\begin{figure}\begin{center}
{\includegraphics[height=5cm,width=13cm]{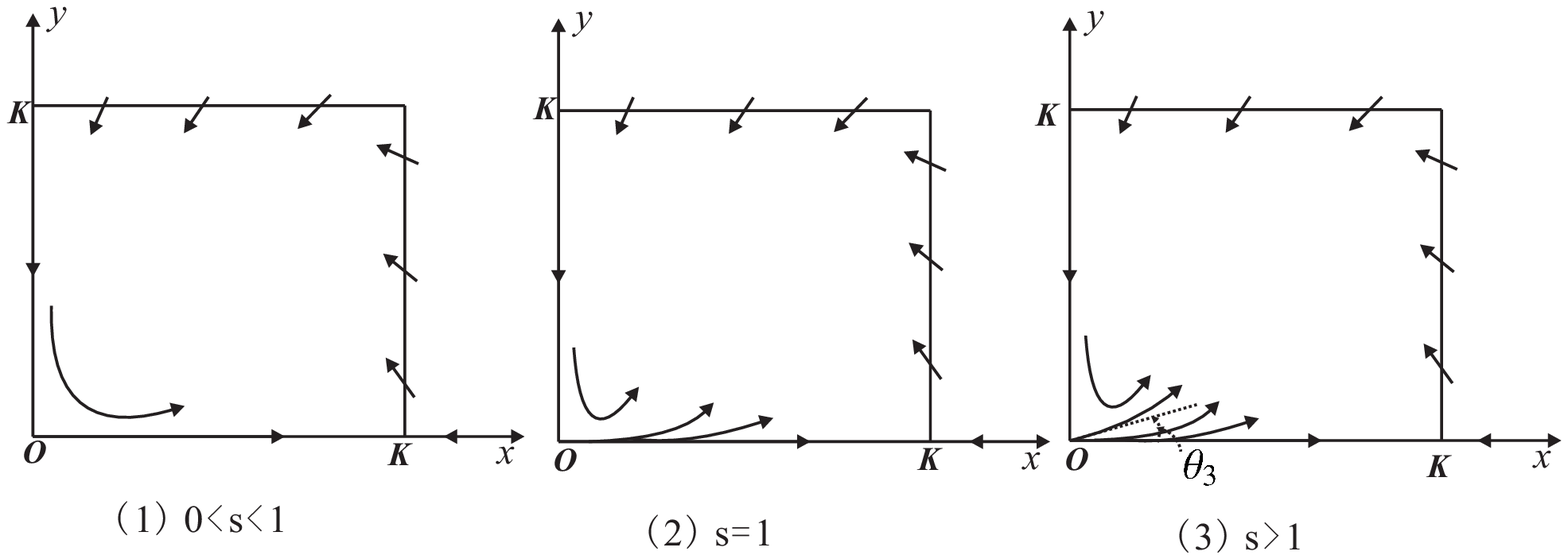}}
\end{center}
{\small {\bf Fig. 2.1.} The sketch graph of the vector field on the boundary of $\Omega$ and topological structure of the orbits of system \eqref{norm} near the origin. }
\end{figure}

Now we are going to study the positive (i.e., interior) equilibria of system \eqref{norm}.
Assume that $\tilde{E}(\tilde{x}, \tilde{y})$ is a positive equilibrium of system \eqref{norm}, then $\tilde{y}=\tilde{x}$ and $\tilde{x}$ is a positive root of the equation
\begin{equation}\label{dai2.13}
\begin{aligned}(x-1)(a x^2+(K b+K-1)x+K)=0
\end{aligned}
\end{equation}
on the interval $(0, K)$.
The Jacobian matrix  of system \eqref{norm} at $E^*(1, 1)$ is
\begin{equation*}
\textrm{J}(E^*)=\left(\begin{array}{cc}
K a-K-2 a-b & -(K-1) (a+b+1)\\
K s (a+b+1)  & -K s (a+b+1)
\end{array}\right).
\end{equation*}
Then the characteristic polynomial of $\textrm{J}(E^*)$ is
\begin{equation*}
\Theta(\lambda)=\lambda^2-\textrm{tr}\left(J(E^{*})\right)\lambda+\textrm{det}\left(J(E^{*})\right),
\end{equation*}
where
\begin{equation}\label{TR}
\begin{aligned}
&\textrm{tr}\left(J(E^{*})\right)=
K a-K-2 a-b-K s (a+b+1),\\
&\textrm{det}\left(J(E^{*})\right)=
K s (a+b+1) (K b+2 K+a-1).
\end{aligned}
\end{equation}
Hence, $E^*$ is a saddle if $\textrm{det}\left(J(E^{*})\right)<0$, is a degenerate equilibrium if $\textrm{det}(J(E^{*})) =0$, and is an elementary equilibrium if
$\textrm{det}\left(J(E^{*})\right)>0$, respectively. Specially, $E^*$ is \emph{center or focus type} if $\textrm{det}\left(J(E^{*})\right)>0$ and $\textrm{tr}\left(J(E^{*})\right)=0$.

Now we assume that $\textrm{det}\left(J(E^{*})\right)>0$, i.e., $K b+2 K+a-1>0$. Then the stability of $E^*$ is determined by the sign of $\textrm{tr}\left(J(E^{*})\right)$.
This gives that $E^*$ is locally asymptotically stable (or unstable) if $\textrm{tr}\left(J(E^{*})\right)<$ (or $>$) $0$.
Specially, $E^{*}$ is locally asymptotically stable if $K a-K-2 a-b\leq 0$.

Furthermore, a straightforward calculation shows that
\begin{equation}\label{dai2.10}
\left(\textrm{tr}\left(J(E^{*})\right)\right)^2-4 \textrm{det}\left(J(E^{*})\right)=\mu_2 {s}^{2}+ \mu_1 s+ \mu_0 \triangleq \psi(s),
\end{equation}
where $\mu_2={K}^{2} (a+b+1)^{2}>0$, $\mu_1=-2 K ( a+b+1)
 ( K a+2 K b+3 K-b-2 )$, $\mu_0=( K a-K-2 a-b )^{2}\geq 0$.
Note that $\tilde{\Delta}=\mu_1^2-4 \mu_2 \mu_0=16\,{K}^{2}(a+b+1)^{3} (K-1) (K b+2 K+a-1)>0$. From \eqref{dai2.10}, if $K a+2 K b+3 K-b-2\leq 0$, then $\psi(s)>0$ and hence $E^{*}$ is a node.

Let  $K a+2 K b+3 K-b-2>0$.
If $K a-K-2 a-b=0$, then $\psi(s)$ has a unique positive zero at $s_1$, where $s_1=-\mu_1/\mu_2$. Furthermore, $\psi(s)\geq 0$  for $s\in[s_1, +\infty)$ and
$\psi(s)<0$ for $s\in(0, s_1)$, respectively.
Notice that $E^{*}$ is locally asymptotically stable for this case. Hence, if $s>s_1$, then $E^*$ is a locally asymptotically stable node. Otherwise, it is a locally asymptotically stable focus.
Let $K a-K-2 a-b\neq 0$, then $\psi(s)$ has two positive zeros at $s_2$ and $s_3$ with $0<s_2<s_3$, where
\begin{align*}
&s_2=\frac{K a+2 K b+3 K-b-2-2 \sqrt{( K-1)( a+b+1)(K b+2 K+a-1)}}{ K (a+b+1)},\\
&s_3=\frac{K a+2 K b+3 K-b-2+2 \sqrt{( K-1)( a+b+1)(K b+2 K+a-1)}}{ K (a+b+1)}.
\end{align*}
Furthermore, $\psi(s)\geq 0$ for $s\in(0, s_2]\cup[s_3, +\infty)$ and $\psi(s)<0$ for $s\in(s_2, s_3)$, respectively. Thus $E^*$ is a node for $s\in(0, s_2)\cup(s_3, +\infty)$ and is a focus or center for  $s\in(s_2, s_3)$.

We note that, if $K a-K-2 a-b>0$, then $K a+2 K b+3 K-b-2=2 (K b+2 K+a-1)+(K a-K-2 a-b)>0$.
For this case, $\textrm{tr}\left(J(E^{*})\right)$ as a function of $s$ has a positive zero at $s^{*}$,  where $s^{*}=(K a-K-2a-b)/(K(a+b+1))$. Since $\psi(s^{*})=-4 (K a-K-2 a-b) (K b+2
K+a-1)<0$, we have $s^{*}\in(s_2, s_3)$.

Summarizing the above  discussions, we arrive at the following results.
\begin{lemm}\label{L2.1} For the positive equilibrium $E^{*}(1, 1)$, the following statements hold.
\begin{itemize}
\item[(a)] Suppose $K b+2 K+a-1<0$, then $E^{*}(1, 1)$ is a hyperbolic saddle;
\item[(b)] Suppose $K b+2 K+a-1=0$, then $E^{*}(1, 1)$ is a degenerate equilibrium;
\item[(c)] Suppose $K b+2 K+a-1>0$, then $E^{*}(1, 1)$ is a node, if either $K a+2 K b+3 K-b-2\leq 0$ or $K a+2 K b+3 K-b-2>0$ and $s\in(0, s_2)\cup(s_3, +\infty)$.  And it is a focus or center, if $s\in(s_2, s_3)$. More precisely,
\begin{itemize}
\item[(c1)] if $K a-K-2a-b<0$,  then $E^*$ is a locally asymptotically stable node, if either $K a+2 K b+3 K-b-2\leq 0$ or $K a+2 K b+3 K-b-2>0$ and $s\in(0, s_2)\cup(s_3, +\infty)$.  And it is a locally asymptotically stable focus, if  $K a+2 K b+3 K-b-2>0$ and  $s\in(s_2, s_3)$.
\item[(c2)] if $K a-K-2a-b=0$, then $E^{*}$ is a locally asymptotically stable node (or focus), if $s\in(s_1, +\infty)$ (or $s\in(0, s_1)$).
\item[(c3)]  Let  $K a-K-2a-b>0$.
\begin{itemize}
\item[(i)] If $s^{*}<s<s_3$ (or $s>s_3$), then $E^{*}(1, 1)$ is a locally asymptotically stable focus (or node);
\item[(ii)]  If $s_2<s<s^{*}$ (or $s<s_2$), then $E^{*}(1, 1)$ is an unstable focus (or node);
\item[(iii)] If $s=s^{*}$, then $E^{*}(1, 1)$ is a weak focus or center.
\end{itemize}
\end{itemize}
\end{itemize}
\end{lemm}
The number of positive  equilibria of
system \eqref{norm} is determined by the number of roots of Eq. \eqref{dai2.13} on the interval $(0, K)$.
Obviously,  Eq. \eqref{dai2.13} can have one, two or
three positive roots on the interval $(0, K)$. Correspondingly, system \eqref{norm} can have one, two, or three
positive equilibria. For system \eqref{norm}, it can be inferred from \cite{Hsu1995} that if $\tilde{x}$ is a multiple positive root of Eq. \eqref{dai2.13}, then $\tilde{E}(\tilde{x},
\tilde{y})$ must be a degenerate positive equilibrium of system \eqref{norm}.
Applying the results of Lemma \ref{L2.1}, after tedious analysis, we can get the following results which are equivalent to Lemma 2.1 of \cite{Huang2014}.
\begin{lemm}\label{L2.2}
Let $\bar{\Delta}=(K b+K-1)^2-4aK$. Then the following statements hold.
\begin{itemize}
\item[(a)] Suppose  $b>1/K-1-2\sqrt{a/K}$, then system \eqref{norm} has a unique positive equilibrium $E^{*}(1, 1)$, which is an elementary and anti-saddle equilibrium;
\item[(b)] Suppose  $b=1/K-1-2\sqrt{a/K}$.
\begin{itemize}
\item[(b1)] If $a=K$, i.e., $b=1/K-3$, then system \eqref{norm} has a unique positive equilibrium $E^{*}(1, 1)$, which is a degenerate equilibrium;
\item[(b2)] If $a\neq K$, then system \eqref{norm} has two different positive equilibria: a degenerate equilibrium $E_2^{*}(\sqrt{K/a}, \sqrt{K/a})$ and an elementary anti-saddle
    equilibrium $E^{*}(1, 1)$;
\end{itemize}
\item[(c)] Suppose  $b<1/K-1-2\sqrt{a/K}$.
\begin{itemize}
\item[(c1)] If $b=(1-a)/K-2$, then system \eqref{norm} has two different positive equilibria: a degenerate equilibrium $E^{*}(1, 1)$ and an elementary anti-saddle equilibrium $E_2^{*}(K/a,
    K/a)$;
\item[(c2)] If $b\neq (1-a)/K-2$, then system \eqref{norm} has three distinct positive equilibria $E^{*}(1, 1)$, $E_2^{*}((1-Kb-K-\sqrt{\bar{\Delta}})/(2a), (1-Kb-K-\sqrt{\bar{\Delta}})/(2a))$ and
    $E_3^{*}((1-K b-K+\sqrt{\bar{\Delta}})/(2 a), (1-K b-K+\sqrt{\bar{\Delta}})/(2 a))$, which are all elementary equilibria. Without loss of generality, assume that $1$ is the minimum root of  Eq.
    \eqref{dai2.13}. Then both $E^{*}$ and $E_3^{*}$ are anti-saddle, and $E_2^{*}$ is a saddle.
\end{itemize}
\end{itemize}
\end{lemm}
The cases (b) and (c1) of Lemma \ref{L2.2} have been investigated by Huang et al \cite{Huang2014}. In the rest of this paper  we only need to focus on the cases (a) and (c2) of Lemma \ref{L2.2}.
Noting that Huang et al \cite{Huang2014} considered all possible cases when system has at least one degenerate positive equilibrium, we only need to consider the case that $E^{*}(1, 1)$
is non-degenerate.
In order words, the type and stability of $E^{*}(1, 1)$ have been determined for all the cases of Lemma \ref{L2.1} except the subcase (iii). Hence, we only need to  determine if $E^{*}(1, 1)$ is a center or a weak focus under the corresponding conditions.

\begin{rema}\label{rem2.2}
The qualitative behavior near the origin shows that the origin is always a repeller in the interior of the first quadrant, which corrects the results in \cite{Eduardo2014} that the origin can be a global attractor.
It's pointed out that we analyze the isolated critical point $O$ of higher order with a different approach. By Lemma \ref{L2.1}, we know that the origin $O$ cannot be an $\omega$ limit set of the unstable manifold of the hyperbolic saddle $E_0$, i.e., there is no heteroclinic orbit. This also corrects the wrong results in \cite{Eduardo2014} that there exists a heteroclinic curve joining the hyperbolic saddle $E_0$ and the origin $O$.
\end{rema}

\subsection{Computation of Lyapunov constants}

In this subsection, we are going to study Hopf bifurcation of system \eqref{norm}.
Without loss of generality,  assume that $E^{*}(1, 1)$ is center or focus type, which implies that the condition (iii) in Lemma \ref{L2.1} (c3) holds, i.e.,
\begin{equation}\label{cond2.3}
\begin{aligned}
K b+2 K+a-1>0,\ K a-K-2 a-b>0,\  s=s^{*}.
\end{aligned}
\end{equation}

To determine if $E^{*}(1, 1)$ is a center or a weak focus, one needs to calculate the focal values of system \eqref{norm} at this equilibrium.
In this paper, we will use the Lyapunov
constants instead of the focal values to solve these problems.
The equivalence between the Lyapunov constants and the focal values can be seen in \cite{Liu2001} and Chapter $1$ in \cite{Liu2008}.

Denote the $k$-th Lyapunov constants of system \eqref{norm} at  $E^{*}(1, 1)$ by $V_{2k+1}$, $k=0, 1, \cdots$. Then $V_1=\textrm{Tr}\left(J(E^{*}(1, 1))\right)$. Under the hypothesis $V_1=0$,
i.e., $s=s^{*}$, we will compute higher Lyapunov constants of system \eqref{norm} at $E^{*}(1, 1)$. To this end, we consider a
conjugation $\psi: \mathbb{R}^2\rightarrow \mathbb{R}^2$, such that
\begin{equation*}
\psi (u, v)=(u, v)+(1, 1)=(x, y).
\end{equation*}
For brevity, we denote $d\triangleq Kb+2K+a-1>0,\ e\triangleq K a-K-2a-b>0$.
Applying the algorithm in \cite{Sang1} and using the software Maple for symbolic calculas, we have
\begin{equation}\label{dai2.19}
\begin{aligned}
V_3=-\dfrac{(K-1)^2 \varphi_1}{4d},
V_5=\dfrac{(K-1)^3 \varphi_2}{48d^3e},
V_7=-\dfrac{(K-1)^4 \varphi_3}{9216d^5e^2},
V_9=\dfrac{(K-1)^5 \varphi_4}{1105920d^7e^3},
\end{aligned}
\end{equation}
where the quantity $V_k$ is reduced w.r.t. the Gr\"{o}bner basis of $\{V_j: j<k\}$ and the expressions of $\varphi_i$'s, $i=1, 2, 3, 4$, we omit here for their cumbersome expressions are available upon request.
This requires a large-scale computation. For example, for $i=1, 2, 3, 4$, the number of  terms of $\varphi_i$ is, respectively, $29$, $340$, $1216$ and $2947$. From \eqref{cond2.3}, we shall consider the Lyapunov constants in \eqref{dai2.19} with $\mu=(K, a, b)\in \Lambda$, where
\begin{equation}\label{basecond}
\Lambda=\left\{\mu\in\mathbb{R}^3\left|K>1, a>0, \textup{max}\bigg\{\frac{1-a}{K}-2, -2\sqrt{a}\bigg\}<b<K a-K- 2 a\right.\right\}.
\end{equation}

\begin{rema}\label{R3.2}
From Lemmas \ref{L2.1} and \ref{L2.2}, if the conditions in the case (a) of Lemma \ref{L2.2} and the case (ii) of Lemma \ref{L2.1} (c2) hold, i.e.,
\begin{equation*}
\begin{aligned}
\textup{max}\bigg\{\frac{1}{K}-1-2\sqrt{\frac{a}{K}},\ -2\sqrt{a}\bigg\}<b<Ka-K-2a, \ s<s^{*},
\end{aligned}
\end{equation*}
then $E^{*}(1, 1)$ is the unique positive equilibrium of system \eqref{norm} and is unstable. Hence system \eqref{norm} has at least one stable limit cycle in $\Omega$ from Theorem \ref{T3.2}.
\end{rema}

\section{Proofs}

Notice that system \eqref{norm2} is equivalent to system \eqref{model} in $\mathcal{A}$. Next we only need to consider the polynomial system \eqref{norm}, instead of system \eqref{model}.
\subsection{Proof of Theorem \ref{T3.2}}
If $\Delta>0$, then the equation \eqref{dai1.2} has a unique positive root in the interval $(0, K)$.
Correspondingly, system \eqref{model} has a unique positive equilibrium.
Similar to the analysis of Lemma 2.1 in \cite{Huang2014}, the first assertion is derived immediately.

Next we will prove the second assertion. From Lemma \ref{T3.1}, every solution $(x(t),  y(t))$ of system \eqref{norm} with positive initial values will eventually be away from the origin (see Fig. 2.1). Hence, the boundary
of the region $\Omega$ can be used as the outer boundary of a Poincar\'{e}-Bendixson annular region. By applying the Poincar\'{e}-Bendixson Theorem, the conclusion of this theorem is clearly established.

\subsection{Proof of Theorem \ref{T4.1}}

From the analysis in Section $2$, the conditions of this theorem imply that $E^{*}(1, 1)$ is the unique positive equilibrium of system \eqref{norm} and is center or focus type.
Noting that $b>1/K-1-2\sqrt{a/K}$ implies that $b>(1-a)/K-2$. From Lemmas \ref{L2.1} and \ref{L2.2}, if the conditions in Theorem \ref{T4.1} hold, then $s=s^{*}$ and $\mu\in\Sigma_1$, where
\begin{equation*}
\Sigma_1=\left\{\mu\in\mathbb{R}^3\left|K>1, a>0, \textup{max}\left\{\frac{1}{K}-1-2\sqrt{\frac{a}{K}},\
-2\sqrt{a}\right\}<b<K a-K-2 a\right.\right\}.
\end{equation*}

To prove this theorem, we first prove that $\varphi_1, \varphi_2, \varphi_3$ have a unique common real root in $\Sigma_1$.
We consider these polynomials in the ring $\textbf{R}[K, a, b]$.

Firstly, computing the mutual resultants of $\varphi_{1}, \varphi_{2}, \varphi_{3}$ with respect to $b$, we get
\begin{equation}\label{dai4.2}
\left\{
\begin{aligned}
&\varphi_{12}:=\textrm{Res}(\varphi_{1}, \varphi_{2}, b)
=32K^{4}a(K-1)^{13}(a-1)^{18}(a-K)^2g_{1b}g_{2}^2g_{31}, \\
&\varphi_{13}:=\textrm{Res}(\varphi_{1}, \varphi_{3}, b)=512K^{5}a(K-1)^{21}(a-1)^{28}(a-K)^3g_{1b}g_{2}^2g_{32}, \\
&\varphi_{23}:=\textrm{Res}(\varphi_{2}, \varphi_{3}, b)=-480K^{10}a(K-1)^{63}(a-1)^{64}(a-K)^{10}g_{1b}g_{2}^2g_{33},
\end{aligned}
\right.
\end{equation}
where $g_{1b}=K^3 a-4 K^2 a^2-K^3+10 K a^2+6 K^2-11 K a-4 a^2-3 K+6 a$, $g_{2}=2{K}^{2}+K-1-(3 K-1)a$, and
$g_{31}$, $g_{32}$, $g_{33}$ are three complicated polynomials in $K$, $a$ of degree $19$, $39$,  $103$, respectively. And more, $\varphi_{12}, \varphi_{13}, \varphi_{23}$ are well factored
over the rational field.
Using the resultant elimination theory (Theorem \ref{L4.3} in Appendix A),  we have $\textbf{V}(\varphi_1, \varphi_2, \varphi_3)=\textbf{V}(\varphi_1, \varphi_2, \varphi_3, \varphi_{12},
\varphi_{13}, \varphi_{23})$.
To find the common real roots of $\{\varphi_1, \varphi_2, \varphi_3, \varphi_{12}, \varphi_{13}, \varphi_{23}\}$, from \eqref{dai4.2}, we have the following five cases to discuss.

\textbf{Case} (a): $a=1$. Substituting it to $\varphi_1$, we get
$\varphi_{1}=-(b+4) (b+2)^2(K^2-4 K+1-K b)$.
Hence $\varphi_1=0$ implies that $b=-2$, $-4$ or $(K^2-4 K+1)/K$.
The cases $b=-2$ and  $b=-4$ can be excluded since $b>-2\sqrt{a}$.
The case $b=(K^2-4 K+1)/K$ follows that $K a-K-2 a-b=-(K-1)^2/K<0$, which contradicts to
$K a-K-2 a-b>0$. Therefore, for this case, we have $\textbf{V}(\varphi_1)=\emptyset$.

\textbf{Case} (b): $a=K$. With this substitution, we have
\begin{equation}\label{case(b)}
\varphi_{1}=(K b+3K-1)\phi_1, \,\,
\varphi_{2}=(K b+3K-1)^2(K^2-3 K-b)\phi_2,
\end{equation}
where $\phi_{1}$ and $\phi_{2}$ are polynomials in $K$, $b$ of degree $4$ and $9$, respectively.
If $b=1/K-3$, then $K b+2 K+a-1=0$, which contradicts to \eqref{cond2.3}.
If $b=K^2-3K$, then $K a-K-2 a-b=0$, which contradicts to  \eqref{cond2.3}.
Thus, it follows from \eqref{case(b)} that $\textbf{V}(\varphi_1, \varphi_2)=\textbf{V}(\phi_1, \phi_2)$. By Sturm's Theorem, we conclude that the resultant
$\textrm{Res}(\phi_{1}, \phi_{2}, b)$ is a nonzero constant, implying that $\textbf{V}(\phi_1, \phi_2)=\emptyset$. Therefore, for this case, we have $\textbf{V}(\varphi_1,
\varphi_2)=\emptyset$.

\textbf{Case} (c): $g_2=0$, i.e., $a=(2K^2+K-1)/(3K-1)$. Then we have
\begin{equation}\label{case(c)}
\varphi_{1}=
(3 K b+8 K-b-4)\phi_{4}/(3K-1)^3, \,\,
\varphi_{2}=(3 K b+8 K-b-4)\phi_{5}/(3K-1)^8,
\end{equation}
where $\phi_{4}$, $\phi_{5}$ are  polynomials in $K$, $b$ of degree $7$, $21$, respectively.
If $3 K b+8 K-b-4=0$, then $b=-4(2K-1)/(3K-1)$. This follows $K b+2 K+a-1=0$, which contradicts to \eqref{cond2.3}. Thus, it follows from \eqref{case(c)} that $\textbf{V}(\varphi_1,
\varphi_2)=\textbf{V}(\phi_4, \phi_5)$. Similar to the case (b), we conclude that $\textbf{V}(\phi_4, \phi_5)=\emptyset$. Hence, for this case, we have  $\textbf{V}(\varphi_1,
\varphi_2)=\emptyset$.

\textbf{Case} (d): $g_{31}=g_{32}=g_{33}=0$. Calculating the  mutual resultants of $g_{31}, g_{32}, g_{33}$  with respect to $a$, we get
\begin{equation}\label{dai4.3}
\begin{aligned}
&\textrm{Res}(g_{31}, g_{32}, a)=
c_1 K^{12}(2K-1)^2(K-1)^{147}(19 K^2-6 K+2) h_1 h_2 h_3 h_4, \\
&\textrm{Res}(g_{31}, g_{33}, a)=
c_2 K^{12}(2K-1)^2(K-1)^{372}(19 K^2-6 K+2) h_1 h_2 h_3 h_5,
\end{aligned}
\end{equation}
where $c_1$, $c_2$ are two constants, and $h_1, h_2, h_3$, $h_4$ and $h_5$ are polynomials in $K$ of degree $3$, $8$, $42$, $274$ and $830$, respectively.  Obviously, we have $19 K^2-6 K+2>0$
since $K>1$. In addition, by Sturm's Theorem, we conclude that $h_1>0$ and $h_2>0$ for $(1, +\infty)$. Furthermore, $h_4$ and $h_5$ have no common zeros because the resultant $\textrm{Res}(h_4,
h_5, K)$ is a nonzero constant.
Thus, from \eqref{dai4.3} and Theorem \ref{L4.3}  in Appendix A, in this case we have $\textbf{V}(\varphi_1, \varphi_2, \varphi_3)=\textbf{V}(\varphi_1, \varphi_2, \varphi_3, g_{31}, g_{32},
g_{33}, h_3)=\textbf{V}(h_3, g_{31}, \varphi_1)\cap\textbf{V}(g_{32}, g_{33}, \varphi_2, \varphi_3)$.
Next, we first find out the common real roots of $\{h_3, g_{31}, \varphi_1\}$ in $\Sigma_1$, and then verify whether they are really the common real roots of $\{\varphi_2, \varphi_3\}$. Even
though $\{h_3, g_{31}, \varphi_1\}$ is an irreducible ascending set, we can not use the $\textsl{RealRootIsolate}$ command in Maple to obtain the isolated real roots in $\Sigma_1$ directly for
their extremely complicated expressions. Here, we will use the real root isolation algorithm of multivariate polynomial systems, proposed by Lu et al \cite{Lu2007},   to solve this problem.
According to the sign of $b$, we have the following three subcases to discuss.

\textbf{Subcase} (d1): $b>0$. Without specification, taking the accuracy of command \emph{realroot} in Maple always being  $1/10^{20}$.
By command $\textsl{realroot}$, we get that $h_3(K)$ has two real roots in $(1, +\infty)$ as follows.
\begin{align*}
&K_{1}\in \left[{\dfrac {354982973071688257601}{295147905179352825856}}, {\dfrac {
1419931892286753030405}{1180591620717411303424}}\right]\triangleq \left[\underline{K}_{1}, \overline{K}_{1}\right],\\
&K_{2}\in \left[{\dfrac {244229161108537248689527983}{1180591620717411303424}},{\dfrac
{15264322569283578043095499}{73786976294838206464}}\right]\triangleq \left[\underline{K}_{2}, \overline{K}_{2}\right].
\end{align*}
For the real root interval $[\underline{K}_{1}, \overline{K}_{1}]$ of $h_3(K)$,
we get the following two positive real root isolation intervals of $\{h_3(K), g_{31}(K, a)\}$.
\begin{align*}
&a_{1}\in\left[{\dfrac {75253211057447250009}{295147905179352825856}},{\dfrac {
150506422114898607437}{590295810358705651712}}\right]\triangleq \left[\underline{a}_{1}, \overline{a}_{1}\right],\\
&a_{2}\in\left[{\dfrac {494428213075491687541}{295147905179352825856}},{\dfrac {
1977712852306468880883}{1180591620717411303424}}\right]\triangleq \left[\underline{a}_{2}, \overline{a}_{2}\right].
\end{align*}
From Definition 2.4 of \cite{Lu2007}, we denote the maximal and minimal polynomials
of $\varphi_{1}$ on $G_{1}=[\underline{K}_{1}, \overline{K}_{1}]\times[\underline{a}_{1}, \overline{a}_{1}]$
by $\overline{\varphi}_1(b)$ and $\underline{\varphi}_1(b)$, respectively.
By command \emph{realroot}, we know that both $\overline{\varphi}_{1}(b)$ and $\underline{\varphi}_{1}(b)$ have no positive real root. Hence $\varphi_{1}$ has no positive real  root for the
real root isolation interval $G_{1}$. Similarly, $\varphi_{1}$ has no positive real root for the real root isolation interval $G_{2}=[\underline{K}_{1},
\overline{K}_{1}]\times[\underline{a}_{2}, \overline{a}_{2}]$.

Similarly, for the real root interval $[\underline{K}_{2}, \overline{K}_{2}]$ of $h_3(K)$, we can get three positive real root isolation intervals of $\{h_3(K), g_{31}(K, a)\}$ as below.
\begin{align*}
&a_{3}\in\left[{\dfrac {302357397559963930735}{590295810358705651712}},{\dfrac {604714795119927861471}{1180591620717411303424}}
\right]\triangleq \left[\underline{a}_{3}, \overline{a}_{3}\right],\\
&a_{4}\in\left[{\dfrac {111474188614417001592773}{590295810358705651712}}
,{\dfrac {222948377228834003185547}{1180591620717411303424}}
\right]\triangleq \left[\underline{a}_{4}, \overline{a}_{4}\right],\\
&a_{5}\in\left[{\dfrac {4494291717948679361239971}{1180591620717411303424}}
,{\dfrac {1123572929487169840309993}{295147905179352825856}}\right]\triangleq \left[\underline{a}_{5}, \overline{a}_{5}\right].
\end{align*}
For the real root isolation interval $G_{3}=[\underline{K}_{2}, \overline{K}_{2}]\times[\underline{a}_{3}, \overline{a}_{3}]$ of $\{h_3, g_{31}\}$, similar to the approach for $G_{1}$, we get
that $\varphi_{1}$ has a positive real root at $b_{1}$, where
\begin{equation*}
b_{1}\in\left[{\dfrac {125098419344119784764828183}{1180591620717411303424}}
,{\dfrac {62549209672059892382517527}{590295810358705651712}}
\right]\triangleq \left[\underline{b}_{1}, \overline{b}_{1}\right].
\end{equation*}
By Theorem 2.3 of \cite{Lu2007}, we have $Ka-K-2a-b\leq \overline{K}_{2}\overline{a}_{3}-\underline{K}_{2}-2\underline{a}_{3}-\underline{b}_{1}
<0$, which implies that $(K_{2}, a_3, b_1)\not\in\Sigma_1$.

For the real root isolation interval $G_{4}=[\underline{K}_{2}, \overline{K}_{2}]\times[\underline{a}_{4}, \overline{a}_{4}]$ of $\{h_3(K), g_{31}(K, a)\}$, $\varphi_{1}$ has three positive
real roots at $b_{i}$ with $b_{i}\in[\underline{b}_{i}, \overline{b}_{i}]$, $i=2, 3, 4$, where
\begin{align*}
&\left[\underline{b}_{2}, \overline{b}_{2}\right]=\left[{\dfrac {192301334985162977663}{2361183241434822606848}}
,{\dfrac {769205339940651910653}{9444732965739290427392}}
\right],\\
&\left[\underline{b}_{3}, \overline{b}_{3}\right]=\left[{\dfrac {12252658027754888276693}{1180591620717411303424}}
,{\dfrac {6126329013877444138347}{590295810358705651712}}
\right],\\
&\left[\underline{b}_{4}, \overline{b}_{4}\right]=\left[{\dfrac {46120025694613890773675092515}{1180591620717411303424}}
,{\dfrac {46120025694613890773675299959}{1180591620717411303424}}
\right].
\end{align*}
By Theorem 2.3 of \cite{Lu2007}, we have
$\varphi_2(K_{2}, a_4, b_3)\geq {\varphi}^+_{2}(\underline{K}_{2}, \underline{a}_{4}, \underline{b}_{3})+{\varphi}^-_{2}(\overline{K}_{2}, \overline{a}_{4}, \overline{b}_{3})\approx 1.83\times
10^{46}>0$, implying that $(K_{2}, a_4, b_3)$ is not a common real root of $\{\varphi_1, \varphi_2, \varphi_3\}$. In addition,
$(K_{2}, a_4, b_4)\not\in \Sigma_1$ since $K a-K-2 a-b\leq \overline{K}_{2}\overline{a}_{4}-\underline{K}_{2}-2\underline{a}_{4}-\underline{b}_{4}<0$.

It's not difficult to verify that the real root $(K_2, a_4, b_2)\in\Sigma_1$. Similarly, after tedious calculation, we can verify that $(K_2, a_4, b_2)$ is a common real root of $\{h_3, g_{31},
\varphi_2\}$ and also a common real root of $\{h_3, g_{31}, \varphi_3\}$. Thus, we conclude that $(K_2, a_4, b_2)$ is a common real root of $\{\varphi_1, \varphi_2, \varphi_3\}$ in $\Sigma_1$.

For the real root isolation interval $G_{5}=[\underline{K}_{2}, \overline{K}_{2}]\times[\underline{a}_{5}, \overline{a}_{5}]]$ of $\{h_3(K), g_{31}(K, a)\}$, $\varphi_{1}$ has two positive real
roots at $b_{i}$ with $b_{i}\in[\underline{b}_{i}, \overline{b}_{i}]$, $i=5, 6$, where
\begin{align*}
&\left[\underline{b}_{5}, \overline{b}_{5}\right]=\left[{\dfrac {71394370384755146991531}{1180591620717411303424}}
,{\dfrac {17848592596188786747883}{295147905179352825856}}
\right],\\
&\left[\underline{b}_{6}, \overline{b}_{6}\right]=\left[{\dfrac {464853909302125803012485143759}{590295810358705651712}}
,{\dfrac {929707818604251606024970505815}{1180591620717411303424}}
\right].
\end{align*}
The real root $(K_{2}, a_{5}, b_{5})$ is not a common real root of $\{\varphi_1, \varphi_2, \varphi_3\}$ since
$\varphi_3(K_{2}, a_{5}, b_{5})\geq {\varphi}^+_{3}(\underline{K}_{2}, \underline{a}_{5}, \underline{b}_{5})+{\varphi}^-_{3}(\overline{K}_{2}, \overline{a}_{5}, \overline{b}_{5})\approx
2.085\times 10^{91}>0.$
On the other hand, $(K_{2}, a_{5}, b_{6})\not\in \Sigma_1$ because
$Ka-K-2a-b\leq \overline{K}_{2}\overline{a}_{5}-\underline{K}_{2}-2\underline{a}_{5}-\underline{b}_{6}\approx -1.92\times 10^{5}
<0$.

Therefore, for the subcase (d1), we know that $p^{*}=(K_2, a_4, b_2)$ is the unique common real root of $\{\varphi_1, \varphi_2, \varphi_3\}$ in $\Sigma_1$.

\textbf{Subcase} (d2): $b<0$. For this subcase, we will find the common real roots of $\{h_3, g_{31}, \varphi_1\}$ in the cone $(1, +\infty)\times(0, +\infty)\times(-\infty, 0)$  by taking the
transformation $b\rightarrow -b$.
We denote
$\hat{\varphi}_1(K, a, b)=\varphi_1(K, a, -b)$.
Since the variable $b$ exactly appears in $\varphi_1(K, a, b)$ for the triangular set $\{h_3, g_{31}, \hat{\varphi}_1\}$, we only need to find the common real roots of  $\{h_3, g_{31},
\hat{\varphi}_1\}$ in the cone $(1, +\infty)\times(0, +\infty)\times(0, +\infty)$.
Repeating the same arguments as the subcase (d1), we conclude that $\{h_3, g_{31}, \hat{\varphi}_1\}$ have fourteen classes of real roots in the cone $(1, +\infty)\times(0, +\infty)\times(0,
+\infty)$: $(K_1, a_1, b_i)$, $(K_1, a_2, b_j)$, $(K_2, a_3, b_k)$, $(K_2, a_4, b_{18})$, $(K_2, a_5, b_{19})$, $(K_2, a_5, b_{20})$, $i=7, \cdots, 10, j=11, \cdots, 14, k=15, \cdots, 17$,
where $K_1$, $K_2$, $a_3$,  $a_4$ and $a_5$ are the same as in the subcase (d1) and
\begin{align*}
&b_{7}\in\left[{\dfrac {342692738901203887719}{590295810358705651712}}
,{\dfrac {685385477803542424951}{1180591620717411303424}}
\right]\triangleq \left[\underline{b}_{7}, \overline{b}_{7}\right],\\
&b_{8}\in\left[{\dfrac {1338789934782848805353}{1180591620717411303424}}
,{\dfrac {334697483698195198185}{295147905179352825856}}
\right]\triangleq \left[\underline{b}_{8}, \overline{b}_{8}\right],\\
&b_{9}\in\left[{\dfrac {194243729812573167855}{147573952589676412928}}
,{\dfrac {776974919255073511755}{590295810358705651712}}
\right]\triangleq \left[\underline{b}_{9}, \overline{b}_{9}\right],\\
&b_{10}\in\left[{\dfrac {3967870668570439988007}{1180591620717411303424}}
,{\dfrac {3967870668571263500067}{1180591620717411303424}}
\right]\triangleq \left[\underline{b}_{10}, \overline{b}_{10}\right],\\
&b_{11}\in\left[{\dfrac {2889636985159625502057}{1180591620717411303424}}
,{\dfrac {2889637013691542396369}{1180591620717411303424}}
\right]\triangleq \left[\underline{b}_{11}, \overline{b}_{11}\right],\\
&b_{12}\in\left[{\dfrac {386107933062469966493}{147573952589676412928}}
,{\dfrac {1544431759404950266357}{590295810358705651712}}
\right]\triangleq \left[\underline{b}_{12}, \overline{b}_{12}\right],\\
&b_{13}\in\left[{\dfrac {3400585485349590843693}{1180591620717411303424}},{\dfrac {1700292756101773392567}{590295810358705651712}}
\right]\triangleq \left[\underline{b}_{13}, \overline{b}_{13}\right],\\
&b_{14}\in\left[{\dfrac {3105247046679173264103}{590295810358705651712}},{\dfrac {6210494094466506191137}{1180591620717411303424}}\right]\triangleq \left[\underline{b}_{14},
\overline{b}_{14}\right],\\
&b_{15}\in\left[{\dfrac {821227566792260176395}{590295810358705651712}},{\dfrac {1642455133584520352861}{1180591620717411303424}}\right]\triangleq \left[\underline{b}_{15},
\overline{b}_{15}\right],\\
&b_{16}\in\left[{\dfrac {497030357386605526709}{295147905179352825856}},{\dfrac {1988121429546422106927}{1180591620717411303424}}\right]\triangleq \left[\underline{b}_{16},
\overline{b}_{16}\right],\\
&b_{17}\in\left[{\dfrac {4031280671528406533751}{590295810358705651712}}
,{\dfrac {8062561343056813067555}{1180591620717411303424}}\right]\triangleq \left[\underline{b}_{17}, \overline{b}_{17}\right],\\
&b_{18}\in\left[{\dfrac {19438401444953006468837}{1180591620717411303424}}
,{\dfrac {9719200722476503234419}{590295810358705651712}}\right]\triangleq \left[\underline{b}_{18}, \overline{b}_{18}\right],\\
&b_{19}\in\left[{\dfrac {973745532080307675583}{4722366482869645213696}},{\dfrac {15214773938754807431}{73786976294838206464}}\right]\triangleq \left[\underline{b}_{19},
\overline{b}_{19}\right],\\
&b_{20}\in\left[\dfrac {78104754415268652314327}{1180591620717411303424}
,\dfrac {9763094301908581539291}{147573952589676412928}\right]\triangleq \left[\underline{b}_{20}, \overline{b}_{20}\right].
\end{align*}
It's pointed out that $(K_i, a_j, b_k)$ is a common real root of  $\{h_3, g_{31}, \hat{\varphi}_1\}$ if and only if $(K_i, a_j, -b_k)$ is a common real root of  $\{h_3, g_{31}, \varphi_1\}$.
It's easy to verify that the following eight classes of
real roots of $\{h_3, g_{31}, \varphi_1\}$ contradict with $b>-2\sqrt{a}$.
\begin{equation*}
(K_1, a_1, -b_{i}), (K_1, a_2, -b_{j}), (K_2, a_3, -b_{k}),\ i=8, 9, 10,\ j=12, 13, 14,\ k=16, 17.
\end{equation*}
By Theorem 2.3 of \cite{Lu2007}, we can verify that the real roots $(K_1, a_1, -b_{7})$, $(K_1, a_2, -b_{11})$ and $(K_2, a_3, -b_{15})$  contradict to $K a-K-2 a-b>0$. In addition, the real
roots $(K_2, a_4, -b_{18})$ and $(K_2, a_5, -b_{20})$ contradict to $K b+2 K+a-1>0$.
Furthermore, $(K_{2}, a_{5}, b_{19})$ is not a common real root of $\{\varphi_1, \varphi_2, \varphi_3\}$ since
$\hat{\varphi}_2(K_{2}, a_{5}, b_{19})\leq \hat{\varphi}^+_{2}(\overline{K}_{2}, \overline{a}_{5}, \overline{b}_{19})+\hat{\varphi}^-_{2}(\underline{K}_{2}, \underline{a}_{5},
\underline{b}_{19})\approx -2.07\times 10^{50}<0,$
where $\hat{\varphi}_2(K, a, b)=\varphi_2(K, a, -b)$.

Therefore, for the subcase (d2), we have  $\textbf{V}(\varphi_1, \varphi_2, \varphi_3)=\emptyset$.

\textbf{Subcase} (d3): $b=0$.  With this substitution,  we have
$\varphi_{1}=-4a \psi_1$ and $\varphi_{2}=-4a \psi_2$,
where $\psi_1$ and $\psi_2$ are polynomials in $K$, $a$ of degree $3$ and $10$, respectively.
Thus, $\textbf{V}(\varphi_1, \varphi_2)=\textbf{V}(\psi_1, \psi_2)$.
The resultant of $\psi_1$  and $\psi_2$ with respect to $a$ is
\begin{equation}\label{H2}
\textrm{Res}(\psi_1, \psi_2, a)=4096(3 K-1)^{2}(K-1)^{13}H_{1}H_{2},
\end{equation}
where $H_{1}$ is a polynomial in $K$ of degree $5$ and $H_{2}=K^2-6 K+3$. By Sturm's Theorem, we know that $\textrm{Res}(\psi_1, \psi_2, a)=0$ implies $H_{2}=0$. This follows $K=3+\sqrt{6}$.
With this substitution, we get that $\psi_1=(8+3\sqrt{6})(a+3)(3+2\sqrt{6}-5a)/5=0$ implies $a=(3+2\sqrt{6})/5$.
Hence, $K a-K-2 a-b=0$, which contradicts to \eqref{cond2.3}. Therefore, for this case, we have $\textbf{V}(\varphi_1, \varphi_2)=\emptyset$.

To summarize, $\{\varphi_1, \varphi_2, \varphi_3\}$ have a unique common real root
$p^{*}=(K_2, a_4, b_2)\in\Sigma_1$ for all the cases above.

\textbf{Case} (e): $g_{1b}=0$.
From above, if $\{\varphi_1, \varphi_2, \varphi_3\}$ have other common real roots in $\Sigma_1$ except the real root $p^{*}$, then there must be $g_{1b}=0$.

Computing the mutual resultants of $\varphi_{1}, \varphi_{2}, \varphi_{3}$ with respect to $a$, we get
\begin{equation}\label{dai4.4}
\left\{
\begin{aligned}
&\textrm{Res}(\varphi_{1}, \varphi_{2}, a)=-8b(b+2)^{18}(K-1)^{13}(K b+3 K-1)^2 g_{1a}\tilde{g}^2_{2}\tilde{g}_{31}, \\
&\textrm{Res}(\varphi_{1}, \varphi_{3}, a)=-64b(b+2)^{28}(K-1)^{21}(K b+3 K-1)^2 g_{1a}\tilde{g}^2_{2}\tilde{g}_{32}, \\
&\textrm{Res}(\varphi_{2}, \varphi_{3}, a)=-120b(b+2)^{64}(K-1)^{63}
(K b+3K-1)^8 g_{1a}\tilde{g}^2_{2}\tilde{g}_4^2\tilde{g}_{33},
\end{aligned}
\right.
\end{equation}
where $g_{1a}={K}^{3}b+4\,{K}^{3}-8{K}^{2}b-4K{b}^{2}-24{K}^{2}-7Kb+2{b}^{2}+12K+6b$, $\tilde{g}_{2}=3 K b+8 K-b-4$,
$\tilde{g}_{4}=K^2-3 K-b$, and $\tilde{g}_{31}$, $\tilde{g}_{32}$, $\tilde{g}_{33}$ are three polynomials in $K$, $b$ of degree $18$, $37$, $102$, respectively.
Similar to the proof of eliminating  the variable $b$ above, we can prove that there must be $g_{1a}=0$, provided that $\{\varphi_1, \varphi_2, \varphi_3\}$ have other common real roots in
$\Sigma_1$  except $p^{*}$.

In what follows, we will prove, by contradiction, that $\{\varphi_1, \varphi_2, \varphi_3\}$ have no common real roots in $\Sigma_1$ such that  $g_{1b}=g_{1a}=0$.
Otherwise, assume that there exist parameter values in $\Sigma_1$ such that $\varphi_1= \varphi_2=\varphi_3=g_{1b}=g_{1a}=0$ hold.

Firstly, applying  pseudo division and by command \emph{prem} in Maple, we get
\begin{equation*}
(2-4K)^3\varphi_1=q_1g_{1a}+r_1,
\end{equation*}
where $q_1$ and $r_1$ are polynomials in $K$, $a$, $b$ of degree $7$ and  $11$, respectively. It follows from $\varphi_1=g_{1a}=0$ that $r_1=0$, which yields $\textbf{V}(\varphi_1,
g_{1a})=\textbf{V}(g_{1a}, r_1)$. Noting that $r_1=H_3b-H_4=0$,
where both $H_3$ and $H_4$ are polynomials in $K$, $a$ of degree $10$.
If $H_3\neq 0$, then $b=H_4/H_3$ and hence $Ka-K-2a-b=-g_{1b}H_5/H_3=0$, where  $H_5$ is a polynomial in $K$, $a$ of degree $8$. This contradicts to \eqref{cond2.3}. So there must be $H_3=0$,
then $r_1=H_4=0$, which yields $\textbf{V}(g_{1a}, r_1)=\textbf{V}(g_{1a}, H_3, H_4)$.
Applying  pseudo division and by command \emph{prem} in Maple, we obtain
\begin{equation}\label{dai4.5}
(K-2)^2H_4=q_2g_{1b}+4(K-1)^2H_2H_6,
\end{equation}
where $H_6$ and $q_2$ are polynomials in $K$, $a$ of degree $8$ and $6$, respectively, and $H_2$ is given in Eq. \eqref{H2}.
If $K=2$, then $g_{1b}=-2(4a-5)$ and $g_{1a}=-2(3b+10)(b+2)$.
It follows from $g_{1b}=g_{1a}=0$ that $a=5/4$ and $b=-10/3$ or $-2$.
These two cases can be excluded since $b>-2\sqrt{a}$ and $Ka-K-2a-b>0$. Thus, we only need to consider the case $K\neq 2$.  From \eqref{dai4.5} and $H_4=g_{1b}=0$, we have $H_2H_6=0$. This
follows $H_2=0$ or $H_6=0$.
The case $H_2=0$ can be excluded. In fact, if $H_2=0$, i.e., $K=3+\sqrt{6}$, then we have
\begin{equation*}
g_{1b}=2(7\sqrt{6}+17) a (2\sqrt{6}-5a+3)/5, \,\, g_{1a}=-2(2\sqrt{6}+5) b(\sqrt{6}+b+3).
\end{equation*}
From the subcase (3) above, the case $b=0$ can be excluded.
Then $g_{1b}=g_{1a}=0$ implies $a=(3+2\sqrt{6})/5$ and $b=-3-\sqrt{6}$. This follows $K b+2 K+a-1=-(47+18\sqrt{6})/5<0$, which contradicts to \eqref{cond2.3}. Therefore, there must be $H_6=0$
and hence $\textbf{V}(g_{1a}, g_{1b}, H_4)=\textbf{V}(g_{1a}, g_{1b}, H_6)$. Notice that $H_6=-H_7a+H_8$, where $H_7$, $H_8$ are polynomials in $K$ of degree $7$, $8$, respectively.
The case $H_7=0$, implies $H_8=0$, can be excluded since the resultant $\textrm{Res}(H_7, H_8, K)$ is a nonzero constant. Thus, we have $H_7\neq0$ and $a=H_8/H_7$. Substituting it to
$g_{1b}=0$, we get
\begin{equation}\label{dai4.6}
g_{1b}=
2K(2K-1)(K-1)^2(K-2)^2H_9/H_7^2=0,
\end{equation}
where $H_9$ is a polynomial in $K$ of degree $9$.
The analysis above tells us that $K\neq 2$, then from \eqref{dai4.6}, there must be $H_9=0$ and hence $\textbf{V}(g_{1b}, H_6)=\textbf{V}(H_6, H_9)$.

To sum up,  we have $\textbf{V}(\varphi_1, g_{1a}, g_{1b})=\textbf{V}(g_{1a}, g_{1b}, r_1)=\textbf{V}(g_{1a}, g_{1b}, H_3, H_4)=\textbf{V}(g_{1a}, g_{1b}, H_3, H_6)=\textbf{V}(g_{1a}, H_3, H_6,
H_9)=\textbf{V}(g_{1a}, H_6, H_9)\cap\textbf{V}(H_3)$.

Now we will show that $\textbf{V}(g_{1a}, H_6, H_9)=\emptyset$. For convenience, denote  by  $S_a$ a semi-algebraic system whose polynomial equations, non-negative polynomial inequalities,
positive polynomial inequalities and polynomial inequations are given by $F_a:=[H_9, H_6, g_{1a}]$, $N:=[\ ]$, $P:=[K-1, a, K a-K-2a-b, K b+2 K+a-1]$, and $H:=[\ ]$, respectively, where $[\ ]$
represents the null set.
It's obvious that the regular chain $\{H_9, H_6, g_{1a}\}$ is squarefree.
Using the \textit{RealRootIsolate} program in Maple and taking the accuracy being $1/10^{20}$, we find that
the solution set of the semi-algebraic system $S_a$ is null, implying that $\textbf{V}(g_{1a}, H_6, H_9)=\emptyset$.  Hence, $\textbf{V}(\varphi_1, \varphi_2, \varphi_3, g_{1a},
g_{1b})=\textbf{V}(\varphi_2, \varphi_3, H_3)\cap\textbf{V}(g_{1a}, H_6, H_9)=\textbf{V}(\varphi_2, \varphi_3, H_3)\cap\emptyset=\emptyset$.
This leads to a contradiction and proves our assertion.

To sum up, we know that $V_{2i+1}=0$, $i=1, 2, 3$, if and only if
\begin{equation}\label{dai4.7}
(K, a, b)=(K_{2}, a_{4}, b_{2}).
\end{equation}
By Theorem 2.3 of \cite{Lu2007},  we have
\begin{equation*}
\varphi_4(K_{2}, a_{4}, b_{2})\leq {\varphi}^+_{4}(\overline{K}_{2}, \overline{a}_{4}, \overline{b}_{2})+{\varphi}^-_{4}(\underline{K}_{2}, \underline{a}_{4}, \underline{b}_{2})\approx -{
1.725035818\times 10^{103}}<0,
\end{equation*}
which yields $V_9(K_2, a_4, b_2)<0$.
This implies that the unique positive equilibrium $E^{*}(1, 1)$ of system \eqref{norm} is not a center but a weak focus of order at most $4$, provided that it is center or focus type.
Furthermore, it is a weak focus of order $4$ if and only if $s=s^{*}$ and \eqref{dai4.7} hold.
By Theorem 2.3.2 of \cite{Han2013}, at most $4$ limit cycles can bifurcate from $E^{*}(1, 1)$. This proves the assertions (1) and (2).

Next we will prove the assertion (3).
Similar to determining the sign of $\varphi_4(K_{2}, a_{4}, b_{2})$, we can check by Maple that the Jacobian determinant of $\varphi_1$, $\varphi_2$,  $\varphi_3$ with respect to $K$, $a$, $b$
at $(K_2, a_4, b_2)$ is negative, which means that the Jacobian matrix of  $V_1$, $V_3$, $V_5$,  $V_7$ with respect to $s$, $K$, $a$, $b$ has its full rank $4$, i.e.,
\begin{eqnarray}\label{dai4.8}
\textrm{rank}\left[\frac{\partial (V_1, V_3, V_5, V_7)}{\partial (s, K, a, b)}\right]_{(\hat{s}, K_2, a_4, b_2)}=4,
\end{eqnarray}
where $\hat{s}=s^{*}(K_2, a_4, b_2)$.
By Theorem 2.3.2 of \cite{Han2013}, there exist some parameter values such that system \eqref{norm} has $4$ small limit cycles around  $E^{*}(1, 1)$
and two of them are stable (see Fig. 3.1. (a)).
Obviously, we can choose appropriate parameter perturbations such that system \eqref{norm} has $i$ small limit cycles around  $E^{*}(1, 1)$ for each $i=1, 2, 3$. This completes the proof.

\subsection{Proof of Theorem \ref{T4.3}}

If $\Delta<0$ and $-2\sqrt{a}<b<K a-K m h/r$, then the equation \eqref{dai1.2} has three different positive roots in the interval $(0, K)$.
Correspondingly, system \eqref{model} has three different positive equilibria, i.e., system \eqref{norm} has three distinct positive equilibria.
By Lemma \ref{L2.2}, we know that  one of the three positive equilibria is a saddle and the others are anti-saddles.
And we have $\mu\in \Sigma_2$, where
\begin{equation}\label{cond4.9}
\Sigma_2=\left\{\mu\in\mathbb{R}^3\left|K>1,\, a>0,\, -2\sqrt{a}<
b<\dfrac{1}{K}-1-2\sqrt{\dfrac{a}{K}},\, b\neq \frac{1}{K}-2-\frac{a}{K}\right.\right\}.
\end{equation}
Without loss of generality, assume that $E^{*}(1, 1)$ is one of the two positive anti-saddles.
If $E^{*}(1, 1)$ is center or focus type, then $s=s^{*}$ and $\mu\in \Lambda$.
Moreover, further assume that the value $1$ is the minimal root of  Eq. \eqref{dai2.13} in $(0, K)$, whereas the proof of the case that $1$ is the maximum  root of  Eq. \eqref{dai2.13} is
similar.
That is, we have $(1-Kb-K-\sqrt{\bar{\Delta}})/(2a)>1$, which implies $\mu\in \Sigma_3$, where
\begin{equation}\label{cond2.30}
\Sigma_3=\left\{\mu\in\mathbb{R}^3\left|K b+2 K+a-1>0, \, K b+K+2 a-1<0\right.\right\}.
\end{equation}

To prove the assertion (1), it suffices to prove the first three Lyapunov constants  in \eqref{dai2.19}  have no common real root, i.e., $\{\varphi_1, \varphi_2, \varphi_3\}$  have no common
real root,  such that $\mu\in \Lambda\cap\Sigma_2\cap\Sigma_3$.
Using the same arguments as the proof of Theorem \ref{T4.1}, we only need to check whether the parameter values $(K_2, a_4, b_2)$ locate in $\Lambda\cap\Sigma_2\cap\Sigma_3$.
Since $b_2>0$, we have $(K_2, a_4, b_2)\not\in\Sigma_2$.
This implies that the center or focus type equilibrium $E^{*}(1, 1)$ is not a center but a weak focus of order at most $3$. By Theorem 2.3.2 of \cite{Han2013}, at most $3$ limit cycles can
bifurcate from it.

To prove the assertion (2), we first find some parameter values such that $E^{*}(1, 1)$ is an order two weak focus.
To this end, we set $K=100$, $a=60$  and hence $\varphi_1$ is a polynomial of $b$. By command $\textsl{realroot}$ with accuracy $1/10^{20}$ in Maple and Theorem 2.3 of \cite{Lu2007}, we can
find $\varphi_1(100, 60, b)$ has a unique real root $b^{\ast}$ such that $(100, 60, b^{\ast})\in \Lambda\cap\Sigma_2\cap\Sigma_3$, where
\begin{align*}
b^{\ast}\in\left[-{\frac {1510607606947999582739}{590295810358705651712}}, -{\frac {
3021215213895999165477}{1180591620717411303424}}
\right]\triangleq \left[\underline{b}^{\ast}, \overline{b^{\ast}}\right].
\end{align*}
Furthermore, by Theorem 2.3 of \cite{Lu2007}, we have
\begin{equation*}
\varphi_2(100, 60, b^{\ast})\leq {\varphi}^+_{2}(100, 60, \overline{b^{\ast}})+{\varphi}^-_{2}(100, 60, \underline{b}^{\ast})\approx -{1.546063838\times 10^{19}}<0,
\end{equation*}
which yields $V_5(100, 60, b^{\ast})<0$.
We first perturb $b$ small near $b^{\ast}$ such that $V_3V_5<0$ and adjust $s$ such that $V_1=0$ holds.  A limit cycle bifurcates. For the second limit cycle, perturb $s$ so that $V_1$ is of the opposite sign of $V_3$. Therefore, two limit cycles can bifurcate from $E^{*}(1, 1)$ (see Fig. 3.2). This completes the proof.

\subsection{Proof of Theorem \ref{T4.4}}

By Theorem \ref{T4.3}, the conditions in Theorem \ref{T4.4} imply that system \eqref{norm} has three distinct positive equilibria,  one of them is a saddle and the others are anti-saddles.
To simplify the calculations, we introduce the new parameters as follows.
Under  the hypothesis that system \eqref{norm} has three distinct positive equilibria,  we assume that  Eq. \eqref{dai2.13} has three positive zeros $1, \alpha$ and $\beta$ in the interval $(0,
K)$  with $1<\alpha<\beta<K$.
Then $E^{*}_2(\alpha, \alpha)$ is a saddle, $E^{*}(1, 1)$ and $E^{*}_3(\beta, \beta)$ are both anti-saddle, respectively.
From Eq. \eqref{dai2.13} and Vieta's formulas for quadratic polynomial, we have
\begin{equation*}
\alpha\beta=\dfrac{K}{a}, \,\,
\alpha+\beta=-\dfrac{Kb+K-1}{a},
\end{equation*}
i.e.,
\begin{equation}\label{dai4.11}
a=\dfrac{K}{\alpha\beta}, \,\,
b=\dfrac{1}{K}-\dfrac{1}{\alpha}-\dfrac{1}{\beta}-1.
\end{equation}
Then, from Lemma \ref{L2.2}, the conditions that system \eqref{norm} has three distinct positive equilibria are equivalent to
\begin{equation}\label{cond2.31}
1<\alpha<\beta<K,\,\,  s>0,\,\, \alpha \beta -K \alpha-K \beta -K\alpha \beta+2 K \sqrt{K\alpha\beta}>0.
\end{equation}
With a time scaling transformation $d\tau=(K-\alpha) ( K-\beta) d t/(K^2 {\alpha}^{2}{\beta}^{2})$, system \eqref{norm} which has been replaced by Eq. \eqref{dai4.11} can be reduced to the following equivalent differential system (we will still use $t$ to denote $\tau$ for ease of
notation).
\begin{equation}\label{dai4.12}
\left\{
\begin{aligned}
\dot{x}=&x^2( K-x )\left[{K}
^{2}{x}^{2}-(K\alpha\,\beta+K\alpha+K\beta-\alpha\,
\beta)x+K\alpha\,\beta\right]\\
&-(K-1)(K-\alpha)(K-\beta) {x}^{3}y
, \\[2pt]
\dot{y}=& K s y (x-y)\left[{K}
^{2}{x}^{2}-(K\alpha\,\beta+K\alpha+K\beta-\alpha\,
\beta)x+K\alpha\,\beta\right],
\end{aligned}
\right.
\end{equation}
where the parameters $K$, $\alpha$, $\beta$ and $s$ satisfy \eqref{cond2.31}.
Next we only need to consider system \eqref{dai4.12}, instead of system \eqref{norm}.

For convenience, we denote the $n$-th Lyapunov constant of system \eqref{dai4.12} at  $E^{*}(1, 1)$ and $E_3^{*}(\beta, \beta)$  by $V_{2n+1}^{(1)}$ and $V_{2n+1}^{(3)}$, respectively, $n=0, 1,
\cdots$. Then  $V_{1}^{(1)}=\textrm{Tr}(J(E^{*}))$ and $V_{1}^{(3)}=\textrm{Tr}(J(E^{*}_3))$, where $J(E^{*})$ is the Jacobian matrix of system \eqref{dai4.12} at $E^{*}$.
From system \eqref{dai4.12}, it's not difficult to obtain
\begin{equation}\label{dai4.13}
\left\{
\begin{aligned}
V_{1}^{(1)}=&-{K}^{2}\alpha\beta+{K}^{3}+K\alpha\beta-2
{K}^{2}+K\alpha+K\beta-\alpha\beta-s K (K-\alpha) (K-\beta),\\
V_{1}^{(3)}=&(K^3\beta-2K^2{\beta}^{2}+K\alpha{\beta}^{2}-{K}^{2}\alpha+K\alpha\beta+
K{\beta}^2-\alpha{\beta}^{2}){\beta}^{2}-\\
&s K{\beta}^3(K-1)(K-\alpha).
\end{aligned}
\right.
\end{equation}
If both $E^{*}(1, 1)$ and $E_3^{*}(\beta, \beta)$ are center or focus type, then $V_{1}^{(1)}=V_{1}^{(3)}$, which implies
\begin{equation}\label{dai4.14}
\left\{
\begin{aligned}
&s={\dfrac {{K}^{3}\beta-2{K}^{2}{\beta}^{2}+{K}^{3}-2{K}^{2}\beta+3
K{\beta}^{2}-2{K}^{2}+3 K\beta-2\,{\beta}^{2}}{K \left( (\beta+1) {K}^{2}
-4\,\beta\,K+{\beta}^{2}+\beta \right) }}\triangleq s_0, \\
&\alpha={\dfrac {K \beta ( 3 {K}^{2}-2 K\beta-2 K+\beta ) }{
(K \beta+K-\beta) ( {K}^{2}-\beta ) }}\triangleq \alpha_0.
\end{aligned}
\right.
\end{equation}
Denote $l(K)\triangleq (\beta+1) {K}^{2}
-4\,\beta\,K+{\beta}^{2}+\beta$, where $l(\cdot)$ is a polynomial function. Since $l'(K)=2(\beta+1)K-4\beta
=2[\beta(K-1)+(K-\beta)]>0$ for $K>\beta>1$, we have $l(K)>l(1)=(\beta-1)^2>0$. Thus, the denominators of the right-hand side of Eq. \eqref{dai4.14} are both
positive.

We first prove the assertion (1). It suffices to prove $V_{3}^{(1)}>0$ and $V_{3}^{(3)}>0$.

Substituting \eqref{dai4.14}  to the system \eqref{dai4.12} and using the software Maple for symbolic calculate, we get
\begin{equation*}
V_{3}^{(1)}=\dfrac {2 {K}^{6}\beta(K-\beta)^{2}(\beta-1)^{2}F_1}{
({K}^{2}-\beta) ^{4}(K\beta+K-\beta)^{3}},\ \
V_{3}^{(3)}=\dfrac {2{K}^{6}{\beta}^{4}(K-1)^{2}(\beta-1)^{2}G_1}{({K}^{2}-\beta)^{4}
(K\beta+K-\beta)^{3}},
\end{equation*}
where $F_1$ and $G_1$ are both polynomials in $K$ and $\beta$ of degree $11$.  From \eqref{cond2.31}, it's easy to see that the denominators of $V_{3}^{(1)}$ and $V_{3}^{(3)}$ are both positive
and the signs of $V_{3}^{(1)}$ and $V_{3}^{(3)}$ are the same as that of $F_1$ and $G_1$, respectively. Hence we only need to prove $F_1>0$ and $G_1>0$.

Denote $\epsilon\triangleq (K-\beta)/(\beta-1)>0$, then $\beta=(\epsilon+K)/(\epsilon+1)$. Substituting it to $F_1$, we obtain
$F_1=(K-1)^{6}\overline{F}_1/(\epsilon+1)^{6}$, where $\overline{F}_1$ is collected by $\epsilon$ as follows.
\begin{eqnarray*}
\begin{aligned}
\overline{F}_1=&(24{K}^{3}-4K+4
){\epsilon}^{6}+(24{K}^{4}+45{K}^{3}-21{K}^{2}+7
K+5) {\epsilon}^{5}+(6{K}^{5}+75{K}^{4}-\\
&
4{K}^{3}-11{K}^{2}+12K+2){\epsilon}^{4}+(33{K}^{5}+42{K}^{4}-18{K}^{3}+3{K}^{2}+6K) {\epsilon}^{3}+\\
&(29{K}^{5}+7{K}^{4}-9{K}^{3}+6{K}^{2}) {\epsilon}^{2}+(11{K}^{5}-4{K}^{4}+2{K}^{3})\epsilon+{K}^{5}.
\end{aligned}
\end{eqnarray*}
It's easy to see that all the coefficients of the power of $\epsilon$ in $\overline{F}_1$ are positive since $K>1$. Noting that $\epsilon>0$, we have $\overline{F}_1>0$, which yields $F_1>0$.

Using the same arguments as the proof of $F_1>0$, we can prove $G_1>0$. This completes the proof of the assertion (1).

Next we will prove the assertion (2). For any given values of $K$ and $\beta$, let $(s, \alpha)=(s_0+\varepsilon, \alpha_0)$, where $s_0$, $\alpha_0$ are given by Eq. \eqref{dai4.14}, and
$\varepsilon$ is a perturbation parameter. From \eqref{dai4.13}, we have
\begin{equation}\label{dai4.15}
L_1=-K ( K-\beta)(K-\alpha_0)\varepsilon,\ \
T_1=-K{\beta}^{3} ( K-1)( K-\alpha_0)\varepsilon.
\end{equation}
By continuity, we can choose $\varepsilon$ small enough such that $V_{3}^{(1)}>0$ and $V_{3}^{(3)}>0$.
From \eqref{dai4.15},  it follows from $K>\beta>\alpha_0>1$ that $V_{1}^{(1)}<0$ and $V_{1}^{(3)}<0$ hold as long as $\varepsilon>0$.
That is to say, there exist parameter values $(s, \alpha)$ near $(s_0, \alpha_0)$ such that $L_1V_{3}^{(1)}<0$ and $T_1V_{3}^{(3)}<0$. This implies that one limit cycle can bifurcate from
$E^{*}(1, 1)$ and $E_3^{*}(\beta, \beta)$, respectively. This ends the proof.

\section{Discussion}

In this paper, we studied the Hopf bifurcation and global dynamics of system \eqref{model} with generalized Holling functional response of type III when it has no degenerate positive equilibrium.
We not only rigorously proved the existence of some phenomena that have been observed in \cite{Hsu1995,Huang2014}, but also observed some new and complicated dynamical behaviors,
such as the coexistence of four limit cycles including two stable cycles, three limit cycles including two stable cycles, one  big stable limit cycle enclosing two unstable limit cycles, or three hyperbolic positive equilibria and two limit cycles surrounding one of them.
Therefore, our results can be a complement to the works by Hsu et al \cite{Hsu1995} and  Huang et al \cite{Huang2014} for this model, and show that the maximum number of limit cycles in the
model of this kind is at least $4$, which improves the preceding results that this number is at least $2$.
Our results also indicate that the nonlinear dynamics of such biological and epidemiological models not only depend on more bifurcation parameters but also are very sensitive to parameter perturbations, which are important for the control of biological species or infectious diseases.
Furthermore, the coexistence of bistable states (one stable limit cycle and one stable equilibrium, or two stable limit cycles) and tristable states (two stable limit cycles and one stable equilibrium, or two stable equilibria  and one stable limit cycle) show that the model of this kind is also highly sensitive to initial values.

In this work, the coexistence of two stable limit cycles gives a positive answer to one of the almost impossible projects proposed in \cite{Coleman1983} which is given in the first section.
This shows the ecological relevance of the coexistence of multiple limit cycles in
predator-prey systems and the importance of our results, which should serve for the outlined
problem to be actually feasible in a biological lab with appropriate little creatures \cite{Coleman1983}. It will be interesting to see if two stable limit cycles occurs in realistic predator-prey systems.

This paper provided some available methods for studying the existence, stability, number and distribution of limit cycles in a given predator-prey systems.
It is worth mentioning that these methods can be applied to other models to study the Hopf bifurcation at positive equilibria with complicated expressions or no explicit expressions. However, there are two interesting problems that remain open. One is to prove the existence of the big stable limit cycle for some given parameter values, see Fig 1.1. If we can rule out the possibility  that there are  two homoclinic loop connecting with the saddle $E_2^*$, called eight-loop in general, then Poincar\'{e}-Bendixson Theorem implies the existence of the big stable limit cycle. We conjecture that this
problem may be solved by constructing an appropriate inner boundary of an annular region. The other is whether $3$ limit cycles can bifurcate from an arbitrary anti-saddle positive equilibrium if system has three distinct positive equilibria, see Remark \ref{rem3}. This problem may be very challenging since both the conditions and the expressions of $V_3$ and $V_5$ are extremely complicated.
Also, is there a natural (or laboratory) predator-prey system with at least three ecologically stable cycles?

\section*{Appendix A}

\subsection{Resultant elimination theory}

Let $\textbf{K}$ be an algebraically closed field. Given two polynomials $A(x_1, x_2, \cdots, x_n)$, $B(x_1, x_2, \cdots, x_n)\in \textbf{K}[x_1, x_2, \cdots, x_n]$, $(x_1, x_2, \cdots, x_n)\in
{\textbf{K}}^{n}$, of the forms
\begin{equation*}
A(x_1, x_2, \cdots, x_n)=\sum_{i=1}^{k}A_i(x_1, x_2, \cdots, x_{n-1})x_n^{i},\,
B(x_1, x_2, \cdots, x_n)=\sum_{i=1}^{l}B_i(x_1, x_2, \cdots, x_{n-1})x_n^{i},
\end{equation*}
where both $k$ and $l$ are positive integers.
Denote the Sylvester resultant of $A$ and $B$ with respect to $x_n$, as defined in \cite{GelfandIM1994},  by $\textrm{Res}(A, B, x_n)$.
Then the following lemma holds (see Theorem 5 in \cite{CollinsGE1971}).
\begin{lemm}\label{L4.2}  Denote $C(x_1, x_2, \cdots, x_{n-1})\triangleq \textrm{Res}(A, B, x_n)$.
If $(a_1, a_2, \cdots, a_{n})$ is a common zero of $A$ and $B$, then $C(a_1, a_2, \cdots, a_{n-1})=0$.
Conversely, if $C(a_1, a_2, \cdots, a_{n-1})=0$, then at least one of the following holds:

(a)\ \ $A_k(a_1, a_2, \cdots, a_{n-1})=\cdots=A_0(a_1, a_2, \cdots, a_{n-1})=0$,

(b)\ \ $B_l(a_1, a_2, \cdots, a_{n-1})=\cdots=B_0(a_1, a_2, \cdots, a_{n-1})=0$,

(c)\ \ $A_k(a_1, a_2, \cdots, a_{n-1})=B_l(a_1, a_2, \cdots, a_{n-1})=0$,

(d)\ \ For some $a_{n}\in \textbf{K}$,  $(a_1, a_2, \cdots, a_{n})$ is a common zero of $A$ and $B$.
\end{lemm}
Clearly, $C=0$ is a necessary condition of  $A=B=0$, but not sufficient. This fact not only gives a criterion for existence of common zeros, but also provides a method of finding the common
zeros of multivariate polynomial systems.
Let $f_1, f_2, \cdots, f_m$ be (finitely many) elements of $\textbf{K}[x_1, \cdots, x_n]$.
Denote the algebraic variety of  $f_1, f_2, \cdots, f_m$, the set of common zeros of $f_1, f_2, \cdots, f_m$, by $\textbf{V}(f_1, \cdots, f_m)$.
Then it follows from Lemma \ref{L4.2} that $\textbf{V}(A, B)=\textbf{V}(A, B, C)$.
This can be generalized to the case of multiple polynomials.
Taking $m=n=3$ as an example, it's not difficult to get the following theorem.
\begin{theo}\label{L4.3}
Denote $r_i^{(k)}(j)\triangleq \textrm{Res}(f_i, f_j, x_k)$, $i, j, k=1, 2, 3$. Then  the following equality holds for any $k\in \{1, 2, 3\}$.
\begin{equation*}
\textbf{V}(f_1, f_2, f_3)=\textbf{V}(f_1, f_2, f_3, r_1^{(k)}(2), r_1^{(k)}(3), r_2^{(k)}(3)).
\end{equation*}
\end{theo}

\section*{Appendix B}
\begin{scriptsize}\quad\ \ The expressions of $\varphi_1$,  $\varphi_2$, $\varphi_3$, $\varphi_4$ in Eq. \eqref{dai2.19}  are listed as follows.
\end{scriptsize}
$$
\begin{scriptsize}
\begin{aligned}
\varphi_1=&-K^2 a b^3+K^2 a^2 b-6 K^2 a b^2+3 K a^3 b+6 K a^2 b^2+6 K a b^3+K b^4-18 K^2 a b-2 K^2 b^2+12 K a^3+27 K a^2b+34K a b^2+6K b^3+a^2 b^2-\\
&16K^2a-3K^2b+16Ka^2+57Kab+12Kb^2-4a^3-3a^2b-6ab^2-b^3+36Ka+9Kb-14ab-
3b^2-12a-3b,\\
\varphi_2=&-24 K^6 a^4 b^5-78 K^6 a^3 b^6-15 K^6 a^2 b^7+24 K^6 a^5 b^3-173 K^6 a^4 b^4-959 K^6 a^3 b^5-153 K^6 a^2 b^6+15 K^6 a b^7+66 K^5 a^6 b^3+308 K^5 a^5 b^4+\\
&628 K^5 a^4 b^5+923 K^5 a^3 b^6+333 K^5 a^2 b^7+30 K^5 a b^8+115 K^6 a^5 b^2-685 K^6 a^4 b^3-5559 K^6 a^3 b^4-818 K^6 a^2 b^5+144 K^6 a b^6+595 K^5 a^6 b^2\\
&+3120 K^5 a^5 b^3+6488 K^5 a^4 b^4+11089 K^5 a^3 b^5+3931 K^5 a^2 b^6+237 K^5 a b^7-15 K^5 b^8+114 K^4 a^7 b^2+116 K^4 a^6 b^3-917 K^4 a^5 b^4-\\
&2489 K^4 a^4 b^5-3666 K^4 a^3 b^6-1824 K^4 a^2 b^7-328 K^4 a b^8-15 K^4 b^9+139 K^6 a^5 b-2002 K^6 a^4 b^2-18450 K^6 a^3 b^3-2357 K^6 a^2 b^4+\\
&788 K^6 a b^5+34 K^6 b^6+1642 K^5 a^6 b+11249 K^5 a^5 b^2+26254 K^5 a^4 b^3+57646 K^5 a^3 b^4+20645 K^5 a^2 b^5+720 K^5 a b^6-144 K^5 b^7+\\
&670 K^4 a^7 b-220 K^4 a^6 b^2-12726 K^4 a^5 b^3-30525 K^4 a^4 b^4-46438 K^4 a^3 b^5-23310 K^4 a^2 b^6-4044 K^4 a b^7-158 K^4 b^8+63 K^3 a^8 b-\\
&274 K^3 a^7 b^2-1428 K^3 a^6 b^3-1476 K^3 a^5 b^4+393 K^3 a^4 b^5+3655 K^3 a^3 b^6+2805 K^3 a^2 b^7+781 K^3 a b^8+73 K^3 b^9-3385 K^6 a^4 b-\\
&34982 K^6 a^3 b^2-1992 K^6 a^2 b^3+2999 K^6 a b^4+271 K^6 b^5+1444 K^5 a^6+17630 K^5 a^5 b+52608 K^5 a^4 b^2+166428 K^5 a^3 b^3+60158 K^5 a^2 b^4+\\
&105 K^5 a b^5-696 K^5 b^6+892 K^4 a^7-3576 K^4 a^6 b-54462 K^4 a^5 b^2-132793 K^4 a^4 b^3-234587 K^4 a^3 b^4-122607 K^4 a^2 b^5-20717 K^4 a b^6-\\
&672 K^4 b^7+220 K^3 a^8-2588 K^3 a^7 b-11078 K^3 a^6 b^2-4835 K^3 a^5 b^3+18789 K^3 a^4 b^4+56096 K^3 a^3 b^5+39785 K^3 a^2 b^6 +10733 K^3 a b^7+\\
&980 K^3 b^8-168 K^2 a^8 b+88 K^2 a^7 b^2+2344 K^2 a^6 b^3+5658 K^2 a^5 b^4+7297 K^2 a^4 b^5+3985 K^2 a^3 b^6+896 K^2 a^2 b^7+49 K^2 a b^8-4 K^2 b^9-\\
&2192 K^6 a^4-34418 K^6 a^3 b+5944 K^6 a^2 b^2+7162 K^6 a b^3+850 K^6 b^4+10188 K^5 a^5+51256 K^5 a^4 b+277194 K^5 a^3 b^2+93250 K^5 a^2 b^3-\\
&8702 K^5 a b^4-2455 K^5 b^5-5044 K^4 a^6-90753 K^4 a^5 b-257833 K^4 a^4 b^2-616634 K^4 a^3 b^3-342099 K^4 a^2 b^4-55211 K^4 a b^5-1158 K^4 b^6-\\
&5204 K^3 a^7-22452 K^3 a^6 b+20842 K^3 a^5 b^2+124319 K^3 a^4 b^3+314732 K^3 a^3 b^4+221974 K^3 a^2 b^5+59517 K^3 a b^6+5382 K^3 b^7-676 K^2 a^8\\
&+2847 K^2 a^7 b+23457 K^2 a^6 b^2+49262 K^2 a^5 b^3+60483 K^2 a^4 b^4+25252 K^2 a^3 b^5+257 K^2 a^2 b^6-1843 K^2 a b^7-263 K^2 b^8+243 K a^8 b+\\
&1023 K a^7 b^2+2047 K a^6 b^3+2481 K a^5 b^4+1152 K a^4 b^5+153 K a^3 b^6-23 K a^2 b^7-4 K a b^8-13360 K^6 a^3+14494 K^6 a^2 b+9035 K^6 a b^2+\\
&1221 K^6 b^3+18760 K^5 a^4+249164 K^5 a^3 b+52523 K^5 a^2 b^2-33164 K^5 a b^3-5978 K^5 b^4-52432 K^4 a^5-221851 K^4 a^4 b-913686 K^4 a^3 b^2-\\
&532240 K^4 a^2 b^3-72716 K^4 a b^4+1295 K^4 b^5-8116 K^3 a^6+88516 K^3 a^5 b+298058 K^3 a^4 b^2+867498 K^3 a^3 b^3+640510 K^3 a^2 b^4+\\
&173416 K^3 a b^5+15647 K^3 b^6+8608 K^2 a^7+59471 K^2 a^6 b+118483 K^2 a^5 b^2+154271 K^2 a^4 b^3+12211 K^2 a^3 b^4-52055 K^2 a^2 b^5-\\
&21961 K^2 a b^6-2485 K^2 b^7+1164 K a^8+4006 K a^7 b+4568 K a^6 b^2+2411 K a^5 b^3-7470 K a^4 b^4-7039 K a^3 b^5-2005 K a^2 b^6-134 K a b^7+\\
&8 K b^8+81 a^7 b^2+203 a^6 b^3+90 a^5 b^4+3 a^4 b^5-4 a^3 b^6+8784 K^6 a^2+4471 K^6 a b+690 K^6 b^2+93816 K^5 a^3-32118 K^5 a^2 b-51403 K^5 a b^2\\
&-8274 K^5 b^3-62552 K^4 a^4-744012 K^4 a^3 b-416826 K^4 a^2 b^2-16448 K^4 a b^3+10044 K^4 b^4+72780 K^3 a^5+273482 K^3 a^4 b+\\
&1290930 K^3 a^3 b^2+1031122 K^3 a^2 b^3+281328 K^3 a b^4+25009 K^3 b^5+32276 K^2 a^6+66224 K^2 a^5 b+132509 K^2 a^4 b^2-205610 K^2 a^3 b^3-\\
&267131 K^2 a^2 b^4-92725 K^2 a b^5-10036 K^2 b^6-1612 K a^7-15436 K a^6 b-29574 K a^5 b^2-61652 K a^4 b^3-36456 K a^3 b^4-4671 K a^2 b^5+\\
&1496 K a b^6+287 K b^7-324 a^8-1047 a^7 b-1434 a^6 b^2-2144 a^5 b^3-1152 a^4 b^4-177 a^3 b^5+23 a^2 b^6+4 a b^7-16 K^6 a+43 K^6 b-\\
&38636 K^5 a^2-30922 K^5 a b-5134 K^5 b^2-267892 K^4 a^3-81694 K^4 a^2 b+69450 K^4 a b^2+18293 K^4 b^3+64596 K^3 a^4+1036780 K^3 a^3 b+\\
&891166 K^3 a^2 b^2+232121 K^3 a b^3+18373 K^3 b^4-17112 K^2 a^5+39060 K^2 a^4 b-537454 K^2 a^3 b^2-572436 K^2 a^2 b^3-192853 K^2 a b^4-\\
&20979 K^2 b^5-14008 K a^6-29560 K a^5 b-104733 K a^4 b^2-35758 K a^3 b^3+25170 K a^2 b^4+14933 K a b^5+1922 K b^6+100 a^7+1543 a^6 b+\\
&91 a^5 b^2+3386 a^4 b^3+3373 a^3 b^4+1036 a^2 b^5+73 a b^6-4 b^7-2628 K^5 a-716 K^5 b+56396 K^4 a^2+64799 K^4 a b+13223 K^4 b^2+\\
&376100 K^3 a^3+314572 K^3 a^2 b+49878 K^3 a b^2-1565 K^3 b^3+27348 K^2 a^4-558241 K^2 a^3 b-602033 K^2 a^2 b^2-204548 K^2 a b^3-22576 K^2 b^4\\
&+3992 K a^5-80052 K a^4 b+47139 K a^3 b^2+108523 K a^2 b^3+44907 K a b^4+5445 K b^5+2088 a^6+460 a^5 b+11473 a^4 b^2+8846 a^3 b^3+\\
&1189 a^2 b^4-506 a b^5-97 b^6+8712 K^4 a+2465 K^4 b-20732 K^3 a^2-39764 K^3 a b-9282 K^3 b^2-242688 K^2 a^3-269921 K^2 a^2 b-\\
&89173 K^2 a b^2-9411 K^2 b^3-39364 K a^4+112054 K a^3 b+156126 K a^2 b^2+61283 K a b^3+7528 K b^4-2536 a^5+13310 a^4 b+5379 a^3 b^2-\\
&4967 a^2 b^3-3123 a b^4-426 b^5-8348 K^3 a-2401 K^3 b-7428 K^2 a^2+3650 K^2 a b+1499 K^2 b^2+71396 K a^3+89064 K a^2 b+35540 K a b^2+\\
&4554 K b^3+8764 a^4-5119 a^3 b-13028 a^2 b^2-5918 a b^3-777 b^4+3192 K^2 a+918 K^2 b+5744 K a^2+3196 K a b+471 K b^2-7868 a^3-\\
&10229 a^2 b-4439 a b^2-606 b^3-432 K a-123 K b-864 a^2-678 a b-123 b^2,\\
\varphi_3=&-18951 K^{10} a^7 b^7-112782 K^{10} a^6 b^8-187937 K^{10} a^5 b^9-62993 K^{10} a^4 b^{10}-5289 K^{10} a^3 b^{11}+18951 K^{10} a^8 b^5-178690 K^{10} a^7 b^6-\\
&1747491 K^{10} a^6 b^7-34{10}665 K^{10} a^5 b^8-{10}17161 K^{10} a^4 b^9-5479 K^{10} a^3 b^{10}+{10}578 K^{10} a^2 b^{11}+53529 K^9 a^9 b^5+372544 K^9 a^8 b^6+\\
&116{10}95 K^9 a^7 b^7+25673{10} K^9 a^6 b^8+3466028 K^9 a^5 b^9+1662059 K^9 a^4 b^{10}+301148 K^9 a^3 b^{11}+15867 K^9 a^2 b^{12}+182390 K^{10} a^8 b^4-\\
&757867 K^{10} a^7 b^5-12896093 K^{10} a^6 b^6-29198695 K^{10} a^5 b^7-7568997 K^{10} a^4 b^8+744785 K^{10} a^3 b^9+152676 K^{10} a^2 b^{10}-5289 K^{10} a b^{11}+
\end{aligned}
\end{scriptsize}
$$
$$
\begin{scriptsize}
\begin{aligned}
\qquad &765482 K^9 a^9 b^4+5725549 K^9 a^8 b^5+18193304 K^9 a^7 b^6+42183074 K^9 a^6 b^7+62143124 K^9 a^5 b^8+28951570 K^9 a^4 b^9+4132412 K^9 a^3 b^{10}-\\
&60759 K^9 a^2 b^{11}-21156 K^9 a b^{12}+145080 K^8 a^{10} b^4+515180 K^8 a^9 b^5-9{10}394 K^8 a^8 b^6-6950294 K^8 a^7 b^7-17323887 K^8 a^6 b^8-\\&24653771 K^8 a^5 b^9-15040180 K^8 a^4 b^{10}-4140723 K^8 a^3 b^{11}-478004 K^8 a^2 b^{12}-15867 K^8 a b^{13}+661797 K^{10} a^8 b^3-2520579 K^{10} a^7 b^4-\\
&59341174 K^{10} a^6 b^5-152577257 K^{10} a^5 b^6-32039308 K^{10} a^4 b^7+8874305 K^{10} a^3 b^8+1{10}2508 K^{10} a^2 b^9-72082 K^{10} a b^{10}+\\
&4160645 K^9 a^9 b^3+35322569 K^9 a^8 b^4+119977423 K^9 a^7 b^5+305155875 K^9 a^6 b^6+5066802{10} K^9 a^5 b^7+228601263 K^9 a^4 b^8+\\
&21876137 K^9 a^3 b^9-3479965 K^9 a^2 b^{10}-277332 K^9 a b^{11}+5289 K^9 b^{12}+1618096 K^8 a^{10} b^3+5263047 K^8 a^9 b^4-22333894 K^8 a^8 b^5-\\
&128447547 K^8 a^7 b^6-309055978 K^8 a^6 b^7-455379133 K^8 a^5 b^8-276246221 K^8 a^4 b^9-70782186 K^8 a^3 b^{10}-6237660 K^8 a^2 b^{11}+\\
&64147 K^8 a b^{12}+{10}578 K^8 b^{13}+150390 K^7 a^{11} b^3-472851 K^7 a^{10} b^4-6194624 K^7 a^9 b^5-12942012 K^7 a^8 b^6-2777985 K^7 a^7 b^7+\\
&33566450 K^7 a^6 b^8+743463{10} K^7 a^5 b^9+59782633 K^7 a^4 b^{10}+22560525 K^7 a^3 b^{11}+4081403 K^7 a^2 b^{12}+304536 K^7 a b^{13}+5289 K^7 b^{14}+\\
&1072680 K^{10} a^8 b^2-8022168 K^{10} a^7 b^3-182434147 K^{10} a^6 b^4-525040369 K^{10} a^5 b^5-70406042 K^{10} a^4 b^6+55477095 K^{10} a^3 b^7+\\
&5074351 K^{10}a^2 b^8-5303{10} K^{10}a b^9-12122 K^{10} b^{10}+{10}935654 K^9a^9 b^2+113290280 K^9 a^8 b^3+434385404 K^9 a^7 b^4+1277108633 K^9a^6b^5\\
&+2468307382K^9a^5b^6+1059752603 K^9 a^4 b^7+30569360 K^9 a^3 b^8-34648505 K^9 a^2 b^9-1649796 K^9 a b^{10}+72082 K^9 b^{11}+\\
&6519204 K^8 a^{10} b^2+19470408 K^8 a^9 b^3-180679066 K^8 a^8 b^4-95020{10}26 K^8 a^7 b^5-2331899848 K^8 a^6 b^6-3703404311 K^8 a^5 b^7-\\
&2251712546 K^8 a^4 b^8-526444840 K^8 a^3 b^9-26507683 K^8 a^2 b^{10}+3673022 K^8 a b^{11}+150982 K^8 b^{12}+1263990 K^7 a^{11} b^2-\\
&7592445 K^7 a^{10} b^3-81915792 K^7 a^9 b^4-145506578 K^7 a^8 b^5+93019176 K^7 a^7 b^6+750466370 K^7 a^6 b^7+1493120848 K^7 a^5 b^8+\\
&1171979681 K^7 a^4 b^9+427983918 K^7 a^3 b^{10}+71899593 K^7 a^2 b^{11}+4340284 K^7 a b^{12}+2091 K^7 b^{13}+77016 K^6 a^{12} b^2-{10}79632 K^6 a^{11} b^3-\\
&3725992 K^6 a^{10} b^4+11621385 K^6 a^9 b^5+60492346 K^6 a^8 b^6+{10}0995253 K^6 a^7 b^7+65462435 K^6 a^6 b^8-49520115 K^6 a^5 b^9-\\
&91640181 K^6 a^4 b^{10}-50628885 K^6 a^3 b^{11}-13144159 K^6 a^2 b^{12}-1575648 K^6 a b^{13}-64687 K^6 b^{14}+652226 K^{10} a^8 b-18088648 K^{10} a^7 b^2-\\
&371799604 K^{10} a^6 b^3-1206919179 K^{10} a^5 b^4-1706354 K^{10} a^4 b^5+221127292 K^{10} a^3 b^6+13766465 K^{10} a^2 b^7-3023576 K^{10} a b^8-\\
&174671 K^{10} b^9+14027606 K^9 a^9 b+200500458 K^9 a^8 b^2+936027430 K^9 a^7 b^3+3392951626 K^9 a^6 b^4+7884024147 K^9 a^5 b^5+\\
&3063477107 K^9 a^4 b^6-275876151 K^9 a^3 b^7-190686054 K^9 a^2 b^8-5573757 K^9 a b^9+504836 K^9 b^{10}+11429217 K^8 a^{10} b+31178346 K^8 a^9 b^2\\
&-695252176 K^8 a^8 b^3-3708670063 K^8 a^7 b^4-9775027861 K^8 a^6 b^5-17496562494 K^8 a^5 b^6-{10}674133418 K^8 a^4 b^7-2165016431 K^8 a^3 b^8\\
&+29232350 K^8 a^2 b^9+35946315 K^8 a b^{10}+9642{10} K^8 b^{11}+3346489 K^7 a^{11} b-39740308 K^7 a^{10} b^2-408763718 K^7 a^9 b^3-\\
&577357247 K^7 a^8 b^4+1465968811 K^7 a^7 b^5+6532167983 K^7 a^6 b^6+12775882514 K^7 a^5 b^7+9995023957 K^7 a^4 b^8+3538000384 K^7 a^3 b^9+\\
&539375663 K^7 a^2 b^{10}+21543128 K^7 a b^{11}-925928 K^7 b^{12}+472459 K^6 a^{12} b-{10}430336 K^6 a^{11} b^2-26195364 K^6 a^{10} b^3+221129811 K^6 a^9 b^4\\
&+907423797 K^6 a^8 b^5+1324950314 K^6 a^7 b^6+521615832 K^6 a^6 b^7-1538163246 K^6 a^5 b^8-2072587533 K^6 a^4 b^9-{10}67884084 K^6 a^3 b^{10}-\\
&264922922 K^6 a^2 b^{11}-30044319 K^6 a b^{12}-1112571 K^6 b^{13}+12627 K^5 a^{13} b-662504 K^5 a^{12} b^2+1579791 K^5 a^{11} b^3+20564380 K^5 a^{10} b^4\\
&+31217841 K^5 a^9 b^5-39929197 K^5 a^8 b^6-191860563 K^5 a^7 b^7-286484322 K^5 a^6 b^8-173630763 K^5 a^5 b^9-20840709 K^5 a^4 b^{10}+\\
&2504{10}28 K^5 a^3 b^{11}+127672{10} K^5 a^2 b^{12}+2389047 K^5 a b^{13}+160846 K^5 b^{14}-21778648 K^{10} a^7 b-476082304 K^{10} a^6 b^2-\\
&1816311584 K^{10} a^5 b^3+455693073 K^{10} a^4 b^4+576918459 K^{10} a^3 b^5+9141863 K^{10} a^2 b^6-13518921 K^{10} a b^7-{10}88356 K^{10} b^8+\\
&7021592 K^9 a^9+185944114 K^9 a^8 b+1201817852 K^9 a^7 b^2+58579788{10} K^9 a^6 b^3+16977116896 K^9 a^5 b^4+5204333033 K^9 a^4 b^5-\\
&1973461776 K^9 a^3 b^6-674366490 K^9 a^2 b^7-5878844 K^9 a b^8+2761743 K^9 b^9+7506304 K^8 a^{10}+18631942 K^8 a^9 b-1402049932 K^8 a^8 b^2\\
&-8315423190 K^8 a^7 b^3-24995847401 K^8 a^6 b^4-53149164822 K^8 a^5 b^5-32127558447 K^8 a^4 b^6-4869168738 K^8 a^3 b^7+\\
&823982378K^8 a^2b^8+189773157 K^8 a b^9+3523422 K^8 b^{10}+2774316 K^7 a^{11}-85080239 K^7a^{10} b-969450372 K^7 a^9 b^2-831603519 K^7 a^8 b^3\\
&+8078583446 K^7 a^7 b^4+29958479263 K^7 a^6 b^5+61737245660 K^7 a^5 b^6+48860797479 K^7 a^4 b^7+16638415602 K^7 a^3 b^8+\\
&2169295478 K^7 a^2 b^9+6715000 K^7 a b^{10}-10357050 K^7 b^{11}+698584 K^6 a^{12}-31726798 K^6 a^{11} b-35688427 K^6 a^{10} b^2+\\
&1407961369 K^6 a^9 b^3+5148241867 K^6 a^8 b^4+6454254947 K^6 a^7 b^5-458138623 K^6 a^6 b^6-16776{10}4266 K^6 a^5 b^7-19629838876 K^6 a^4 b^8-\\
&9673219495 K^6 a^3 b^9-2302503210 K^6 a^2 b^{10}-244297969 K^6 a b^{11}-7713683 K^6 b^{12}+42764 K^5 a^{13}-4388263 K^5 a^{12} b+24404974 K^5 a^{11} b^2\\
&+216799265 K^5 a^{10} b^3+160706884 K^5 a^9 b^4-113142{10}48 K^5 a^8 b^5-3367316200 K^5 a^7 b^6-4459841848 K^5 a^6 b^7-2148593680 K^5 a^5 b^8\\
&+378300472 K^5 a^4 b^9+787988582 K^5 a^3 b^{10}+314703403 K^5 a^2 b^{11}+53893864 K^5 a b^{12}+3433663 K^5 b^{13}-1{10}958 K^4 a^{13} b+\\
&2139423 K^4 a^{12} b^2+5648686 K^4 a^{11} b^3-24567169 K^4 a^{10} b^4-111703770 K^4 a^9 b^5-173176552 K^4 a^8 b^6-83808373 K^4 a^7 b^7+\\
&130743703 K^4 a^6 b^8+212231573 K^4 a^5 b^9+128225180 K^4 a^4 b^{10}+38815717 K^4 a^3 b^{11}+5597852 K^4 a^2 b^{12}+231967 K^4 a b^{13}-13175 K^4 b^{14}\\
&-10321696 K^{10} a^7-342612864 K^{10} a^6 b-1695091288 K^{10} a^5 b^2+1323354886 K^{10} a^4 b^3+948364815 K^{10} a^3 b^4-76113038 K^{10} a^2 b^5-\\
&42799073 K^{10} a b^6-3770530 K^{10} b^7+70746992 K^9 a^8+852575664 K^9 a^7 b+6400331852 K^9 a^6 b^2+24383776864 K^9 a^5 b^3+\\
&3358754828 K^9 a^4 b^4-6471484199 K^9 a^3 b^5-1527491459 K^9 a^2 b^6+49655782 K^9 a b^7+12554817 K^9 b^8-129688 K^8 a^9-\\
&1439239465 K^8 a^8 b-10750919404 K^8 a^7 b^2-40059670222 K^8 a^6 b^3-108007660369 K^8 a^5 b^4-62061751909 K^8 a^4 b^5-
\end{aligned}
\end{scriptsize}
$$
$$
\begin{scriptsize}
\begin{aligned}
\qquad &3100998747 K^8 a^3 b^6+4256652876 K^8 a^2 b^7+633359916 K^8 a b^8+5994315 K^8 b^9-63989120 K^7 a^{10}-1091681637 K^7 a^9 b+\\
&332035370 K^7 a^8 b^2+22350863780 K^7 a^7 b^3+80477877764 K^7 a^6 b^4+187363816955 K^7 a^5 b^5+151078086043 K^7 a^4 b^6+\\&48176025115 K^7 a^3 b^7+4442150366 K^7 a^2 b^8-443629110K^7 a b^9-58493374 K^7 b^{10}-30993812 K^6 a^{11}+71824623 K^6 a^{10} b+\\
&3995306486 K^6 a^9 b^2+14057807156 K^6 a^8 b^3+13504847452 K^6 a^7 b^4-17974793785 K^6 a^6 b^5-94564626560 K^6 a^5 b^6-\\
&103893797765 K^6 a^4 b^7-49650729916 K^6 a^3 b^8-11285813670 K^6 a^2 b^9-1085904724 K^6 a b^{10}-24445161 K^6 b^{11}-6838544 K^5 a^{12}+\\
&101887564 K^5 a^{11} b+748706074 K^5 a^{10} b^2-4920{10}649 K^5 a^9 b^3-9017342465 K^5 a^8 b^4-22053596025 K^5 a^7 b^5-27055947259 K^5 a^6 b^6-\\
&7881854170 K^5 a^5 b^7+9654992781 K^5 a^4 b^8+9189909717 K^5 a^3 b^9+3228017099 K^5 a^2 b^{10}+518993041 K^5 a b^{11}+31533704 K^5 b^{12}-\\
&416796 K^4 a^{13}+15809328 K^4 a^{12} b+15062180 K^4 a^{11} b^2-396341375 K^4 a^{10} b^3-1326572897 K^4 a^9 b^4-1517167164 K^4 a^8 b^5+\\
&382396680 K^4 a^7 b^6+3576722360 K^4 a^6 b^7+4061439986 K^4 a^5 b^8+2059560230 K^4 a^4 b^9+482656074 K^4 a^3 b^{10}+26310415 K^4 a^2 b^{11}-\\
&8141275 K^4 a b^{12}-1036160 K^4 b^{13}+455493 K^3 a^{13}b-2205742 K^3 a^{12} b^2-18362054 K^3 a^{11} b^3-27230933 K^3 a^{10} b^4+35602625 K^3 a^9 b^5+\\
&191188504 K^3 a^8 b^6+327619137 K^3 a^7 b^7+282786618 K^3 a^6 b^8+133462941 K^3 a^5 b^9+33346705 K^3 a^4 b^{10}+3391315 K^3 a^3 b^{11}-\\
&116351 K^3 a^2 b^{12}-32273 K^3 a b^{13}+143 K^3 b^{14}-{10}4617664 K^{10} a^6-874775320 K^{10} a^5 b+1841624544 K^{10} a^4 b^2+869493928 K^{10} a^3 b^3-\\
&294253345 K^{10} a^2 b^4-88228415 K^{10} a b^5-7668256 K^{10} b^6+258351152 K^9 a^7+4025027776 K^9 a^6 b+22270411304 K^9 a^5 b^2-\\
&5074157214 K^9 a^4 b^3-12422659420 K^9 a^3 b^4-1936747265 K^9 a^2 b^5+304836674 K^9 a b^6+41992629 K^9 b^7-595776368 K^8 a^8-\\
&7419548648 K^8 a^7 b-39438869764 K^8 a^6 b^2-147137747334 K^8 a^5 b^3-71561639549 K^8 a^4 b^4+14726259914 K^8 a^3 b^5+\\
&12166099734 K^8 a^2 b^6+1326450839 K^8 a b^7-12974625 K^8 b^8-463592012 K^7 a^9+1991560907 K^7 a^8 b+33305664788 K^7 a^7 b^2+\\
&130645186502 K^7 a^6 b^3+373912931214 K^7 a^5 b^4+305128895315 K^7 a^4 b^5+84343238008 K^7 a^3 b^6+627681736 K^7 a^2 b^7-\\
&2525432434 K^7 a b^8-206213839 K^7 b^9+137720424 K^6 a^{10}+5168076124 K^6 a^9 b+18865615809 K^6 a^8 b^2+6261016788 K^6 a^7 b^3-\\
&79252833344 K^6 a^6 b^4-317465574507 K^6 a^5 b^5-342398821553 K^6 a^4 b^6-159236390967 K^6 a^3 b^7-33842904957 K^6 a^2 b^8-\\
&2698442348 K^6 a b^9-5376505 K^6 b^{10}+126779696 K^5 a^{11}+938589108 K^5 a^{10} b-4307789878 K^5 a^9 b^2-31406348469 K^5 a^8 b^3-\\
&70121312722 K^5 a^7 b^4-81522424365 K^5 a^6 b^5+7687873468 K^5 a^5 b^6+76597823663 K^5 a^4 b^7+56695269366 K^5 a^3 b^8+\\
&18410797315 K^5 a^2 b^9+2813652946 K^5 a b^{10}+163054402 K^5 b^{11}+27273568 K^4 a^{12}-86304276 K^4 a^{11} b-1864849380 K^4 a^{10} b^2-\\
&5028450076 K^4 a^9 b^3-2436810064 K^4 a^8 b^4+11586898069 K^4 a^7 b^5+31266041893 K^4 a^6 b^6+29455973188 K^4 a^5 b^7+\\
&12056882536 K^4 a^4 b^8+12934{10}493 K^4 a^3 b^9-575243281 K^4 a^2 b^{10}-188220324 K^4 a b^{11}-16028906 K^4 b^{12}+1838812 K^3 a^{13}-\\
&24071345 K^3 a^{12} b-138571850 K^3 a^{11} b^2-13896513 K^3 a^{10} b^3+{10}39963432 K^3 a^9 b^4+2921627395 K^3 a^8 b^5+4129155168 K^3 a^7 b^6+\\
&2877716262 K^3 a^6 b^7+825806372 K^3 a^5 b^8-78396600 K^3 a^4 b^9-110304130 K^3 a^3 b^{10}-24334819 K^3 a^2 b^{11}-1426868 K^3 a b^{12}+\\
&54988 K^3 b^{13}-804804 K^2 a^{13} b-1178591 K^2 a^{12} b^2+13039306 K^2 a^{11} b^3+6092{10}44 K^2 a^{10} b^4+133817745 K^2 a^9 b^5+173657177 K^2 a^8 b^6\\
&+128001570 K^2 a^7 b^7+51741503 K^2 a^6 b^8+9625992 K^2 a^5 b^9-77007 K^2 a^4 b^{10}-278116 K^2 a^3 b^{11}-24057 K^2 a^2 b^{12}+286 K^2 a b^{13}-\\
&186052320 K^{10} a^5+1309063796 K^{10} a^4 b+257713160 K^{10} a^3 b^2-497942500 K^{10} a^2 b^3-108469633 K^{10} a b^4-8583025 K^{10} b^5+\\
&1107174960 K^9 a^6+11603982508 K^9 a^5 b-13174470848 K^9 a^4 b^2-13874508230 K^9 a^3 b^3-494633114 K^9 a^2 b^4+841349356 K^9 a b^5+\\
&93361133 K^9 b^6-2098315472 K^8 a^7-21909825734 K^8 a^6 b-129622799272 K^8 a^5 b^2-33857464270 K^8 a^4 b^3+48658689571 K^8 a^3 b^4+\\
&20999075801 K^8 a^2 b^5+1406316388 K^8 a b^6-117805888 K^8 b^7+1408749408 K^7 a^8+25295155930 K^7 a^7 b+125426108016 K^7 a^6 b^2+\\
&495961372268 K^7 a^5 b^3+395145948044 K^7 a^4 b^4+71426796041 K^7 a^3 b^5-21842234555 K^7 a^2 b^6-7421240124 K^7 a b^7-469030229 K^7 b^8\\
&+2457141924 K^6 a^9+10650202010 K^6 a^8 b-20076097976 K^6 a^7b^2-164024945800 K^6 a^6 b^3-675460745070 K^6 a^5 b^4-\\
&733878946511 K^6 a^4 b^5-328020673342 K^6 a^3 b^6-61169089178 K^6 a^2 b^7-26874523{10} K^6 a b^8+257995625 K^6 b^9+250679056 K^5 a^{10}-\\
&8402553049 K^5 a^9 b-52500296230 K^5 a^8 b^2-113554287382 K^5 a^7 b^3-125280214254 K^5 a^6 b^4+143799098061 K^5 a^5 b^5+\\
&320937013595 K^5 a^4 b^6+211112466301 K^5 a^3 b^7+64866918830 K^5 a^2 b^8+9428975610 K^5 a b^9+514357482 K^5 b^{10}-224125168 K^4 a^{11}-\\
&3138945957 K^4 a^{10} b-6317653166 K^4 a^9 b^2+10898304790 K^4 a^8 b^3+61482916780 K^4 a^7 b^4+129396633117 K^4 a^6 b^5+\\
&105018872030 K^4 a^5 b^6+26737994015 K^4 a^4 b^7-9203865400 K^4 a^3 b^8-7139923296 K^4 a^2 b^9-1559299266 K^4 a b^{10}-116879633 K^4 b^{11}\\
&-55240416 K^3 a^{12}-211734737 K^3 a^{11} b+{10}83880120 K^3 a^{10} b^2+6442225218 K^3 a^9 b^3+14143313305 K^3 a^8 b^4+17119524577 K^3 a^7 b^5\\
&+7003089837 K^3 a^6 b^6-3327049966 K^3 a^5 b^7-4248235473 K^3 a^4 b^8-1519382556 K^3 a^3 b^9-189{10}8895 K^3 a^2 b^{10}+7712424 K^3 a b^{11}+\\
&2426482 K^3 b^{12}-3646828 K^2 a^{13}+8834973 K^2 a^{12} b+163846{10}4 K^2 a^{11} b^2+550382784 K^2 a^{10} b^3+957644117 K^2 a^9 b^4+\\
&917328068 K^2 a^8 b^5+199070648 K^2 a^7 b^6-343064469 K^2 a^6 b^7-304811917 K^2 a^5 b^8-{10}3022786 K^2 a^4 b^9-13799005 K^2 a^3 b^{10}+\\
&31188 K^2 a^2 b^{11}+{10}0823 K^2 a b^{12}-429 K^2 b^{13}+661626 K a^{13} b+4292000 K a^{12} b^2+13404369 K a^{11} b^3+26830204 K a^{10} b^4+\\
&33406639 K a^9 b^5+22790019 K a^8 b^6+7631372 K a^7 b^7+529544 K a^6 b^8-391470 K a^5 b^9-982{10} K a^4 b^{10}-4816 K a^3 b^{11}+143 K a^2 b^{12}+\\
&379371904 K^{10} a^4-201717576 K^{10} a^3 b-428240864 K^{10} a^2 b^2-64549872 K^{10} a b^3-3437187 K^{10} b^4+2598256096 K^9 a^5-
\end{aligned}
\end{scriptsize}
$$
$$
\begin{scriptsize}
\begin{aligned}
\qquad &11371984172 K^9 a^4 b-7677701980 K^9 a^3 b^2+2246197422 K^9 a^2 b^3+1293688150 K^9 a b^4+127980818 K^9 b^5-5290721520 K^8 a^6-\\
&66902178764 K^8 a^5 b+22293409468 K^8 a^4 b^2+68871290894 K^8 a^3 b^3+20196161493 K^8 a^2 b^4-439011774 K^8 a b^5-361359733 K^8 b^6+\\&7533055960 K^7 a^7+65048777066 K^7 a^6 b+426553142112 K^7 a^5 b^2+295565803962 K^7 a^4 b^3-25835865950 K^7 a^3 b^4-\\
&60959412027 K^7 a^2 b^5-12988277500 K^7 a b^6-627643809 K^7 b^7+1066064592 K^6 a^8-30592131524 K^6 a^7 b-1773093567{10} K^6 a^6 b^2-\\
&930492203314 K^6 a^5 b^3-1024439224682 K^6 a^4 b^4-419217515623 K^6 a^3 b^5-53718738281 K^6 a^2 b^6+4639056076 K^6 a b^7+\\
&1078256918 K^6 b^8-4995156068 K^5 a^9-39224732211 K^5 a^8 b-84697176188 K^5 a^7 b^2-82429723658 K^5 a^6 b^3+483620651{10}2 K^5 a^5 b^4+\\
&811315315919 K^5 a^4 b^5+500447924608 K^5 a^3 b^6+146116150202 K^5 a^2 b^7+19825513642 K^5 a b^8+968983129 K^5 b^9-1434169760 K^4 a^{10}\\
&+726362688 K^4 a^9 b+40638154041 K^4 a^8 b^2+139088518168 K^4 a^7 b^3+287189635986 K^4 a^6 b^4+186855178541 K^4 a^5 b^5-\\
&25480447333 K^4 a^4 b^6-80908019385 K^4 a^3 b^7-36862270373 K^4 a^2 b^8-7010849536 K^4 a b^9-492459963 K^4 b^{10}+100236292 K^3 a^{11}+\\
&3211646571 K^3 a^{10} b+13854841908 K^3 a^9 b^2+25872748953 K^3 a^8 b^3+23857735346 K^3 a^7 b^4-16901758867 K^3 a^6 b^5-\\
&41583360012 K^3 a^5 b^6-25725987443 K^3 a^4 b^7-5973240842 K^3 a^3 b^8+138699538 K^3 a^2 b^9+265888600 K^3 a b^{10}+29234676 K^3 b^{11}+\\
&53028760 K^2 a^{12}+441447126 K^2 a^{11} b+{10}79880479 K^2 a^{10} b^2+938579447 K^2 a^9 b^3-{10}08792573 K^2 a^8 b^4-4652146508 K^2 a^7 b^5-\\
&5034091784 K^2 a^6 b^6-2275842399 K^2 a^5 b^7-293627023 K^2 a^4 b^8+{10}4771626 K^2 a^3 b^9+36797550 K^2 a^2 b^{10}+2671626 K^2 a b^{11}-\\
&84913 K^2 b^{12}+3657496 K a^{13}+18458906 K a^{12} b+30327470 K a^{11} b^2+9199325 K a^{10} b^3-73943158 K a^9 b^4-198302548 K a^8 b^5-\\
&203106376 K a^7 b^6-{10}1406834 K a^6 b^7-23356184 K a^5 b^8-861629 K a^4 b^9+505082 K a^3 b^{10}+49830 K a^2 b^{11}-572 K a b^{12}+220542 a^{12} b^2\\
&+1089352 a^{11} b^3+1885293 a^{10} b^4+1365690 a^9 b^5+390653 a^8 b^6-28613 a^7 b^7-39516 a^6 b^8-5716 a^5 b^9+143 a^4 b^{10}-137781984 K^{10} a^3\\
&-159230448 K^{10} a^2 b-235240 K^{10} a b^2+2552557 K^{10} b^3-3653758384 K^9 a^4-665407680 K^9 a^3 b+3181208820 K^9 a^2 b^2+\\
&1069525307 K^9 a b^3+95248251 K^9 b^4-15364126272 K^8 a^5+37715820414 K^8 a^4 b+48434436076 K^8 a^3 b^2+6402307126 K^8 a^2 b^3-\\
&3462791210 K^8 a b^4-608480633 K^8 b^5+13960657376 K^7 a^6+218250642070 K^7 a^5 b+76423304220 K^7 a^4 b^2-13{10}21655482 K^7 a^3 b^3-\\
&81440748745 K^7 a^2 b^4-12617944047 K^7 a b^5-261783507 K^7 b^6-12138814760 K^6 a^7-93622862314 K^6 a^6 b-818702429916 K^6 a^5 b^2-\\
&885846925096 K^6 a^4 b^3-2771257{10}372 K^6 a^3 b^4+18601124805 K^6 a^2 b^5+20465700468 K^6 a b^6+2287167207 K^6 b^7-91732{10}272 K^5 a^8\\
&-15997757680 K^5 a^7 b+1681366684 K^5 a^6 b^2+841224790642 K^5 a^5 b^3+1289661364440 K^5 a^4 b^4+758350267181 K^5 a^3 b^5+\\
&205519168591 K^5 a^2 b^6+24149401136 K^5 a b^7+839712823 K^5 b^8+3724860804 K^4 a^9+43070757607 K^4 a^8 b+144509229032 K^4 a^7 b^2+\\
&355695337562 K^4 a^6 b^3+{10}6779430458 K^4 a^5 b^4-278649973543 K^4 a^4 b^5-280076874460 K^4 a^3 b^6-{10}8212814070 K^4 a^2 b^7-\\
&19136731258 K^4 a b^8-1283455360 K^4 b^9+1933883648 K^3 a^{10}+91943{10}170 K^3 a^9 b+{10}992967952 K^3 a^8 b^2-{10}120442020 K^3 a^7 b^3-\\
&111616333268 K^3 a^6 b^4-140409574869 K^3 a^5 b^5-62145113231 K^3 a^4 b^6-2123461347 K^3 a^3 b^7+6650684782 K^3 a^2 b^8+\\
&1957302439 K^3 a b^9+171558960 K^3 b^{10}+128413844 K^2 a^{11}-362742847 K^2 a^{10} b-36313439{10} K^2 a^9 b^2-{10}192455111 K^2 a^8 b^3-\\
&19077380097 K^2 a^7 b^4-13784118957 K^2 a^6 b^5-1223519927 K^2 a^5 b^6+3245319889 K^2 a^4 b^7+1637066083 K^2 a^3 b^8+264604486 K^2 a^2 b^9\\
&-1212314 K^2 a b^{10}-2657361 K^2 b^{11}-12261392 K a^{12}-128550960 K a^{11} b-405791666 K a^{10} b^2-755840016 K a^9 b^3-999964634 K a^8 b^4\\
&-516041770 K a^7 b^5+116457124 K a^6 b^6+241767604 K a^5 b^7+{10}0898704 K a^4 b^8+15709218 K a^3 b^9+1673{10} K a^2 b^{10}-103397 K a b^{11}\\
&+429 K b^{12}-882168 a^{13}-5003514 a^{12} b-12175684 a^{11} b^2-21009933 a^{10} b^3-27826206 a^9 b^4-21060848 a^8 b^5-7945629 a^7 b^6-\\
&826702 a^6 b^7+346895 a^5 b^8+99497 a^4 b^9+4816 a^3 b^{10}-143 a^2 b^{11}-8695488 K^{10} a^2+16435768 K^{10} a b+2985048 K^{10} b^2+\\
&807572560 K^9 a^3+1544361552 K^9 a^2 b+364712914 K^9 a b^2+23552398 K^9 b^3+14679238960 K^8 a^4+12400154040 K^8 a^3 b-\\
&5668187376 K^8 a^2 b^2-4335404922 K^8 a b^3-572426893 K^8 b^4+5{10}03446152 K^7 a^5-48538565802 K^7 a^4 b-127953667130 K^7 a^3 b^2-\\
&53628344473 K^7 a^2 b^3-3772382134 K^7 a b^4+591603803 K^7 b^5-18104773840 K^6 a^6-428051584880 K^6 a^5 b-\\
&397945620366 K^6 a^4 b^2+5807498132 K^6 a^3 b^3+102520956368 K^6 a^2 b^4+31709218662 K^6 a b^5+2751254783 K^6 b^6+5551235520 K^5 a^7+\\
&18521232944 K^5 a^6 b+851878423148 K^5 a^5 b^2+1269467124226 K^5 a^4 b^3+694250009194 K^5 a^3 b^4+156690026921 K^5 a^2 b^5+\\
&10059112428 K^5 a b^6-543878513 K^5 b^7+12977675136 K^4 a^8+59477248464 K^4 a^7 b+250497046160 K^4 a^6 b^2-178320782684 K^4 a^5 b^3\\
&-683122609028 K^4 a^4 b^4-544416617441 K^4 a^3 b^5-193696223005 K^4 a^2 b^6-32326468486 K^4 a b^7-2050639962 K^4 b^8+\\
&250880600 K^3 a^9-1{10}24335072 K^3 a^8 b-46936064524 K^3 a^7 b^2-215496755538 K^3 a^6 b^3-218815218582 K^3 a^5 b^4-\\
&40585025549 K^3 a^4 b^5+50524905264 K^3 a^3 b^6+32365819288 K^3 a^2 b^7+7287015318 K^3 a b^8+581566630 K^3 b^9-826414840 K^2 a^{10}-\\
&4902877504 K^2 a^9 b-12009911356 K^2 a^8 b^2-23774814592 K^2 a^7 b^3-2648781597 K^2 a^6 b^4+23127532908 K^2 a^5 b^5+\\
&19449068741 K^2 a^4 b^6+5786245104 K^2 a^3 b^7+228879962 K^2 a^2 b^8-191261918 K^2 a b^9-25262728 K^2 b^{10}-72710944 K a^{11}-\\
&188847692 K a^{10} b-124671656 K a^9 b^2+122253403 K a^8 b^3+2012567950 K a^7 b^4+285{10}51616 K a^6 b^5+1550663982 K a^5 b^6+\\
&283479090 K a^4 b^7-49686830 K a^3 b^8-24308238 K a^2 b^9-1974214 K a b^{10}+57562 K b^{11}+689200 a^{12}+11629768 a^{11} b+33658560 a^{10} b^2\\
&+54227052 a^9 b^3+92919518 a^8 b^4+94130249 a^7 b^5+49662687 a^6 b^6+12408403 a^5 b^7+664954 a^4 b^8-24{10}69 a^3 b^9-25344 a^2 b^{10}+
\end{aligned}
\end{scriptsize}
$$
$$
\begin{scriptsize}
\begin{aligned}
\qquad &286 a b^{11}+3804384 K^{10} a+718394 K^{10} b+162034704 K^9 a^2-24056050 K^9 a b-9942122 K^9 b^2-1346037168 K^8 a^3-5017955563 K^8 a^2 b\\
&-2197699522 K^8 a b^2-261508090 K^8 b^3-31185459168 K^7 a^4-46089910595 K^7 a^3 b-8965510988 K^7 a^2 b^2+4538493306 K^7 a b^3\\&+1076459403 K^7 b^4-103493689368 K^6 a^5-24487788665 K^6 a^4 b+142801503752 K^6 a^3 b^2+105215940748 K^6 a^2 b^3+\\
&23662759523 K^6 a b^4+1575728910 K^6 b^5+827287840 K^5 a^6+491740481169 K^5 a^5 b+700801183844 K^5 a^4 b^2+306485979191 K^5 a^3 b^3+\\
&18626826322 K^5 a^2 b^4-14770958002 K^5 a b^5-2320990082 K^5 b^6+4809584384 K^4 a^7+{10}7093970646 K^4 a^6 b-387670055500 K^4 a^5 b^2-\\
&862244626332 K^4 a^4 b^3-623168949140 K^4 a^3 b^4-206561699273 K^4 a^2 b^5-31513724522 K^4 a b^6-1750651217 K^4 b^7-6242646304 K^3 a^8\\
&-27828670466 K^3 a^7 b-206292065768 K^3 a^6 b^2-144273813172 K^3 a^5 b^3+109799467608 K^3 a^4 b^4+169719135171 K^3 a^3 b^5+\\
&78191327423 K^3 a^2 b^6+15837833888 K^3 a b^7+1203308057 K^3 b^8-821859248 K^2 a^9-1604082188 K^2 a^8 b-9319357680 K^2 a^7 b^2+\\
&37051860860 K^2 a^6 b^3+69441641153 K^2 a^5 b^4+38869702661 K^2 a^4 b^5+4440359554 K^2 a^3 b^6-3322805591 K^2 a^2 b^7-1207929718 K^2 a b^8\\
&-117555558 K^2 b^9+153251696 K a^{10}+829184396 K a^9 b+1725254512 K a^8 b^2+5846350384 K a^7 b^3+5630497524 K a^6 b^4+\\
&980935446 K a^5 b^5-1283387742 K a^4 b^6-766022622 K a^3 b^7-138639316 K a^2 b^8-608872 K a b^9+1377974 K b^{10}+11335088 a^{11}+\\
&33466160 a^{10} b+47461308 a^9 b^2+127207936 a^8 b^3+70899563 a^7 b^4-57700996 a^6 b^5-79279344 a^5 b^6-33290085 a^4 b^7-5411457 a^3 b^8\\
&-78286 a^2 b^9+34847 a b^{10}-143 b^{11}-20899448 K^9 a-4263446 K^9 b-739058256 K^8 a^2-279368714 K^8 a b-29523168 K^8 b^2-\\
&960054084 K^7 a^3+6152303445 K^7 a^2 b+4375514580 K^7 a b^2+703044181 K^7 b^3+35949459288 K^6 a^4+76332881858 K^6 a^3 b+\\
&40125138961 K^6 a^2 b^2+5367070849 K^6 a b^3-119032265 K^6 b^4+129462040196 K^5 a^5+146068693859 K^5 a^4 b-182055{10}578 K^5 a^3 b^2-\\
&67783693943 K^5 a^2 b^3-24000533306 K^5 a b^4-2502256211 K^5 b^5+27082613760 K^4 a^6-305318000470 K^4 a^5 b-583765754799 K^4 a^4 b^2\\
&-387289195590 K^4 a^3 b^3-1{10}661795017 K^4 a^2 b^4-12122077575 K^4 a b^5-193134939 K^4 b^6-2509392152 K^3 a^7-111115982842 K^3 a^6 b\\
&+22088499184 K^3 a^5 b^2+272131236666 K^3 a^4 b^3+261348922678 K^3 a^3 b^4+106963238663 K^3 a^2 b^5+20409548924 K^3 a b^6+\\
&1489751205 K^3 b^7+1213587344 K^2 a^8-2227919252 K^2 a^7 b+59836058582 K^2 a^6 b^2+82523343706 K^2 a^5 b^3+23043055881 K^2 a^4 b^4-\\
&18202928288 K^2 a^3 b^5-14239096754 K^2 a^2 b^6-357443{10}75 K^2 a b^7-311303714 K^2 b^8+188448000 K a^9+345822820 K a^8 b+\\
&5318098632 K a^7 b^2+1618861436 K a^6 b^3-6283005462 K a^5 b^4-6162282386 K a^4 b^5-2034971720 K a^3 b^6-102759488 K a^2 b^7+\\
&71751886 K a b^8+{10}182599 K b^9-15089488 a^{10}-68194176 a^9 b-3431144 a^8 b^2-264735038 a^7 b^3-462228329 a^6 b^4-273373251 a^5 b^5\\
&-47260578 a^4 b^6+14091355 a^3 b^7+6092591 a^2 b^8+499205 a b^9-14462 b^{10}+31331160 K^8 a+6577403 K^8 b+1434154016 K^7 a^2+\\
&976293119 K^7 a b+163832618 K^7 b^2+5754782972 K^6 a^3-14280{10}61 K^6 a^2 b-2772977642 K^6 a b^2-614729964 K^6 b^3-\\
&19781926736 K^5 a^4-60235306076 K^5 a^3 b-45284353430 K^5 a^2 b^2-11902345481 K^5 a b^3-{10}13413079 K^5 b^4-97509741716 K^4 a^5-\\
&169109331434 K^4 a^4 b-85007814556 K^4 a^3 b^2-3808156979 K^4 a^2 b^3+6272328847 K^4 a b^4+1052283535 K^4 b^5-30420448992 K^3 a^6+\\
&92389010569 K^3 a^5 b+250819916666 K^3 a^4 b^2+208850955498 K^3 a^3 b^3+79539422225 K^3 a^2 b^4+14140466800 K^3 a b^5+949185291 K^3 b^6\\
&-1810835736 K^2 a^7+43548042346 K^2 a^6 b+37993593524 K^2 a^5 b^2-29208861698 K^2 a^4 b^3-51849722275 K^2 a^3 b^4-\\
&26416997611 K^2 a^2 b^5-5877080979 K^2 a b^6-488161397 K^2 b^7-23{10}35152 K a^8+3227908192 K a^7 b-5311865576 K a^6 b^2-\\
&13899808634 K a^5 b^3-8679682520 K a^4 b^4-988146970 K a^3 b^5+931990392 K a^2 b^6+354254678 K a b^7+36351572 K b^8-15096864 a^9+\\
&73502980 a^8 b-431781596 a^7 b^2-560343532 a^6 b^3-52791889 a^5 b^4+218776604 a^4 b^5+125267367 a^3 b^6+22758298 a^2 b^7-123566 a b^8\\
&-271781 b^9+7726340 K^7 a+3997135 K^7 b-13{10}056600 K^6 a^2-1175418652 K^6 a b-225918603 K^6 b^2-7250556144 K^5 a^3-\\
&5580421412 K^5 a^2 b-620891638 K^5 a b^2+89937427 K^5 b^3+1726318880 K^4 a^4+20381243444 K^4 a^3 b+21558525316 K^4 a^2 b^2+\\
&7452673784 K^4 a b^3+838406412 K^4 b^4+43757262124 K^3 a^5+94041216223 K^3 a^4 b+72064595174 K^3 a^3 b^2+23657735459 K^3 a^2 b^3+\\
&3027514886 K^3 a b^4+79031190 K^3 b^5+14508573040 K^2 a^6-7968254692 K^2 a^5 b-54510642447 K^2 a^4 b^2-56141586338 K^2 a^3 b^3-\\
&25281216571 K^2 a^2 b^4-5344264241 K^2 a b^5-433190740 K^2 b^6+1275398096 K a^7-7126058408 K a^6 b-11143385332 K a^5 b^2-\\
&2589833256 K a^4 b^3+38969860{10} K a^3 b^4+2953473124 K a^2 b^5+778988322 K a b^6+71893864 K b^7+53966224 a^8-422434864 a^7 b-\\
&129680196 a^6 b^2+630782922 a^5 b^3+634902625 a^4 b^4+209420075 a^3 b^5+3502579 a^2 b^6-10900737 a b^7-1492901 b^8-55491788 K^6 a\\
&-16842508 K^6 b+536241360 K^5 a^2+634{10}4821 K^5 a b+134365066 K^5 b^2+4331378544 K^4 a^3+4838507983 K^4 a^2 b+\\
&1542870538 K^4 a b^2+152963422 K^4 b^3+3558020992 K^3 a^4+457309955 K^3 a^3 b-3180664320 K^3 a^2 b^2-1708276494 K^3 a b^3-\\
&239325145 K^3 b^4-11379085132 K^2 a^5-27864836499 K^2 a^4 b-25674673336 K^2 a^3 b^2-1{10}9{10}06964 K^2 a^2 b^3-2250913394 K^2 a b^4-\\
&174138705 K^2 b^5-3271883344 K a^6-1905198062 K a^5 b+4991748876 K a^4 b^2+7154137167 K a^3 b^3+3764535184 K a^2 b^4+\\
&895915007 K a b^5+80070975 K b^6-201439056 a^7+352725176 a^6 b+899521668 a^5 b^2+534599006 a^4 b^3+3461289 a^3 b^4-\\
&106475700 a^2 b^5-37297098 a b^6-3886319 b^7+53259476 K^5 a+15626167 K^5 b-25579616 K^4 a^2-124848492 K^4 a b-32098009 K^4 b^2\\
&-1318632876 K^3 a^3-1772671041 K^3 a^2 b-715205300 K^3 a b^2-9{10}83619 K^3 b^3-1843845032 K^2 a^4-2342717202 K^2 a^3 b-\\
&903938133 K^2 a^2 b^2-88867705 K^2 a b^3+5112696 K^2 b^4+1574648168 K a^5+4213437650 K a^4 b+4345819594 K a^3 b^2+
\end{aligned}
\end{scriptsize}
$$
$$
\begin{scriptsize}
\begin{aligned}
\qquad &2145367199 K a^2 b^3+506{10}7492 K a b^4+45726114 K b^5+287105808 a^6+341632320 a^5 b-65741238 a^4 b^2-319986130 a^3 b^3-\\
&208731920 a^2 b^4-55835147 a b^5-5383037 b^6-23585364 K^4 a-6837549 K^4 b-55474592 K^3 a^2-19711000 K^3 a b-1488684 K^3 b^2+\\&171990084 K^2 a^3+274683685 K^2 a^2 b+128872514 K^2 a b^2+18284289 K^2 b^3+370848096 K a^4+620610976 K a^3 b+375638554 K a^2 b^2\\
&+97456570 K a b^3+9264306 K b^4-89381608 a^5-254264410 a^4 b-284950416 a^3 b^2-153933335 a^2 b^3-39668261 a b^4-3859476 b^5+\\
&5145696 K^3 a+1478274 K^3 b+18065928 K^2 a^2+12061368 K^2 a b+20{10}762 K^2 b^2+1789392 K a^3-3779676 K a^2 b-3733668 K a b^2-\\
&685995 K b^3-27519360 a^4-52704600 a^3 b-36987996 a^2 b^2-11194734 a b^3-1218483 b^4-448344 K^2 a-127566 K^2 b-1793376 K a^2\\
&-1406952 K a b-255132 K b^2-1793376 a^3-2303640 a^2 b-958608 a b^2-127566 b^3,\\
\varphi_4=&-16476027\,{K}^{14}{a}^{10}{b}^{9}-139967169\,{K}^{14}{a}^{9}{b}^{10}-
422050748\,{K}^{14}{a}^{8}{b}^{11}-522260634\,{K}^{14}{a}^{7}{b}^{12}-
228175353\,{K}^{14}{a}^{6}{b}^{13}-\\
&37081436\,{K}^{14}{a}^{5}{b}^{14}-
1477905\,{K}^{14}{a}^{4}{b}^{15}+16476027\,{K}^{14}{a}^{11}{b}^{7}-
191978596\,{K}^{14}{a}^{10}{b}^{8}-2636867711\,{K}^{14}{a}^{9}{b}^{9}-
9308222958\,{K}^{14}{a}^{8}{b}^{10}\\
&-12304207236\,{K}^{14}{a}^{7}{b}^{
11}-4958673572\,{K}^{14}{a}^{6}{b}^{12}-481040343\,{K}^{14}{a}^{5}{b}^
{13}+55773259\,{K}^{14}{a}^{4}{b}^{14}+4433715\,{K}^{14}{a}^{3}{b}^{15
}+\\
&47398173\,{K}^{13}{a}^{12}{b}^{7}+431584237\,{K}^{13}{a}^{11}{b}^{8}
+1814971927\,{K}^{13}{a}^{10}{b}^{9}+5271634350\,{K}^{13}{a}^{9}{b}^{
10}+10954315659\,{K}^{13}{a}^{8}{b}^{11}+\\
&13165140278\,{K}^{13}{a}^{7}{
b}^{12}+7301254470\,{K}^{13}{a}^{6}{b}^{13}+1883006375\,{K}^{13}{a}^{5
}{b}^{14}+203527747\,{K}^{13}{a}^{4}{b}^{15}+5911620\,{K}^{13}{a}^{3}{
b}^{16}+\\
&235912691\,{K}^{14}{a}^{11}{b}^{6}-834558458\,{K}^{14}{a}^{10}
{b}^{7}-23616625459\,{K}^{14}{a}^{9}{b}^{8}-97317878205\,{K}^{14}{a}^{
8}{b}^{9}-137234709255\,{K}^{14}{a}^{7}{b}^{10}-\\
&49561296183\,{K}^{14}{
a}^{6}{b}^{11}-321517043\,{K}^{14}{a}^{5}{b}^{12}+1561432181\,{K}^{14}
{a}^{4}{b}^{13}+31999960\,{K}^{14}{a}^{3}{b}^{14}-4433715\,{K}^{14}{a}
^{2}{b}^{15}+\\
&922803821\,{K}^{13}{a}^{12}{b}^{6}+8867408496\,{K}^{13}{a
}^{11}{b}^{7}+37979320676\,{K}^{13}{a}^{10}{b}^{8}+112408459137\,{K}^{
13}{a}^{9}{b}^{9}+246171845543\,{K}^{13}{a}^{8}{b}^{10}+\\
&309867162399\,
{K}^{13}{a}^{7}{b}^{11}+166957434408\,{K}^{13}{a}^{6}{b}^{12}+
36662071061\,{K}^{13}{a}^{5}{b}^{13}+1763896903\,{K}^{13}{a}^{4}{b}^{
14}-268587307\,{K}^{13}{a}^{3}{b}^{15}-\\
&13301145\,{K}^{13}{a}^{2}{b}^{
16}+184636032\,{K}^{12}{a}^{13}{b}^{6}+1071876948\,{K}^{12}{a}^{12}{b}
^{7}+459473552\,{K}^{12}{a}^{11}{b}^{8}-11072099055\,{K}^{12}{a}^{10}{
b}^{9}-\\
&45486661982\,{K}^{12}{a}^{9}{b}^{10}-102527905248\,{K}^{12}{a}^
{8}{b}^{11}-134342866618\,{K}^{12}{a}^{7}{b}^{12}-90182950149\,{K}^{12
}{a}^{6}{b}^{13}-31493119151\,{K}^{12}{a}^{5}{b}^{14}-\\
&5552938279\,{K}^
{12}{a}^{4}{b}^{15}-427102792\,{K}^{12}{a}^{3}{b}^{16}-8867430\,{K}^{
12}{a}^{2}{b}^{17}+1412808955\,{K}^{14}{a}^{11}{b}^{5}-1762318940\,{K}
^{14}{a}^{10}{b}^{6}-\\
&134828850755\,{K}^{14}{a}^{9}{b}^{7}-636199934819
\,{K}^{14}{a}^{8}{b}^{8}-952503311604\,{K}^{14}{a}^{7}{b}^{9}-
291921439939\,{K}^{14}{a}^{6}{b}^{10}+40239330792\,{K}^{14}{a}^{5}{b}^
{11}+\\
&16932790595\,{K}^{14}{a}^{4}{b}^{12}-458369333\,{K}^{14}{a}^{3}{b
}^{13}-79581271\,{K}^{14}{a}^{2}{b}^{14}+1477905\,{K}^{14}a{b}^{15}+
7405876880\,{K}^{13}{a}^{12}{b}^{5}+\\
&77289865114\,{K}^{13}{a}^{11}{b}^{
6}+346832804602\,{K}^{13}{a}^{10}{b}^{7}+1074425832882\,{K}^{13}{a}^{9
}{b}^{8}+2535676821952\,{K}^{13}{a}^{8}{b}^{9}+3383978800035\,{K}^{13}
{a}^{7}{b}^{10}\\
&+1755398577505\,{K}^{13}{a}^{6}{b}^{11}+298000402641\,{
K}^{13}{a}^{5}{b}^{12}-17947671943\,{K}^{13}{a}^{4}{b}^{13}-6649478595
\,{K}^{13}{a}^{3}{b}^{14}-59373039\,{K}^{13}{a}^{2}{b}^{15}+\\
&8867430\,{
K}^{13}a{b}^{16}+2988664332\,{K}^{12}{a}^{13}{b}^{5}+17316498186\,{K}^
{12}{a}^{12}{b}^{6}-7606247340\,{K}^{12}{a}^{11}{b}^{7}-282970051276\,
{K}^{12}{a}^{10}{b}^{8}-\\
&1072353557580\,{K}^{12}{a}^{9}{b}^{9}-
2415809443148\,{K}^{12}{a}^{8}{b}^{10}-3238617191747\,{K}^{12}{a}^{7}{
b}^{11}-2149623062689\,{K}^{12}{a}^{6}{b}^{12}-\\
&700869325697\,{K}^{12}{
a}^{5}{b}^{13}-100608812041\,{K}^{12}{a}^{4}{b}^{14}-2800347361\,{K}^{
12}{a}^{3}{b}^{15}+431408289\,{K}^{12}{a}^{2}{b}^{16}+13301145\,{K}^{
12}a{b}^{17}+\\
&299240529\,{K}^{11}{a}^{14}{b}^{5}-3450786\,{K}^{11}{a}^{
13}{b}^{6}-12256753273\,{K}^{11}{a}^{12}{b}^{7}-47111290779\,{K}^{11}{
a}^{11}{b}^{8}-55707410368\,{K}^{11}{a}^{10}{b}^{9}+\\
&69569636328\,{K}^{
11}{a}^{9}{b}^{10}+383965220987\,{K}^{11}{a}^{8}{b}^{11}+665410293717
\,{K}^{11}{a}^{7}{b}^{12}+555335901667\,{K}^{11}{a}^{6}{b}^{13}+
248940984736\,{K}^{11}{a}^{5}{b}^{14}+\\
&61218952540\,{K}^{11}{a}^{4}{b}^
{15}+7805054484\,{K}^{11}{a}^{3}{b}^{16}+430956254\,{K}^{11}{a}^{2}{b}
^{17}+5911620\,{K}^{11}a{b}^{18}+4528095185\,{K}^{14}{a}^{11}{b}^{4}-\\
&
4386712314\,{K}^{14}{a}^{10}{b}^{5}-545417021157\,{K}^{14}{a}^{9}{b}^{
6}-2869986979586\,{K}^{14}{a}^{8}{b}^{7}-4535251107782\,{K}^{14}{a}^{7
}{b}^{8}-1043685062428\,{K}^{14}{a}^{6}{b}^{9}+\\
&454016417691\,{K}^{14}{
a}^{5}{b}^{10}+106738930973\,{K}^{14}{a}^{4}{b}^{11}-9567887975\,{K}^{
14}{a}^{3}{b}^{12}-693389291\,{K}^{14}{a}^{2}{b}^{13}+25476806\,{K}^{
14}a{b}^{14}+\\
&32136160673\,{K}^{13}{a}^{12}{b}^{4}+376355926120\,{K}^{
13}{a}^{11}{b}^{5}+1822087856040\,{K}^{13}{a}^{10}{b}^{6}+
6082291214804\,{K}^{13}{a}^{9}{b}^{7}+15846038683195\,{K}^{13}{a}^{8}{
b}^{8}\\
&+22675137792518\,{K}^{13}{a}^{7}{b}^{9}+11116934352519\,{K}^{13}
{a}^{6}{b}^{10}+1102775794337\,{K}^{13}{a}^{5}{b}^{11}-413839894922\,{
K}^{13}{a}^{4}{b}^{12}-\\
&62529593856\,{K}^{13}{a}^{3}{b}^{13}+2133840961
\,{K}^{13}{a}^{2}{b}^{14}+149909405\,{K}^{13}a{b}^{15}-1477905\,{K}^{
13}{b}^{16}+19664978948\,{K}^{12}{a}^{13}{b}^{4}+\\
&116365658489\,{K}^{12
}{a}^{12}{b}^{5}-191116560787\,{K}^{12}{a}^{11}{b}^{6}-2989385494494\,
{K}^{12}{a}^{10}{b}^{7}-11028146773945\,{K}^{12}{a}^{9}{b}^{8}-\\
&
25528195701294\,{K}^{12}{a}^{8}{b}^{9}-35616311672195\,{K}^{12}{a}^{7}
{b}^{10}-23383615923040\,{K}^{12}{a}^{6}{b}^{11}-6947697210631\,{K}^{
12}{a}^{5}{b}^{12}-\\
&675519118416\,{K}^{12}{a}^{4}{b}^{13}+58144469048\,
{K}^{12}{a}^{3}{b}^{14}+10281007984\,{K}^{12}{a}^{2}{b}^{15}+53592846
\,{K}^{12}a{b}^{16}-4433715\,{K}^{12}{b}^{17}+\\
&4013960610\,{K}^{11}{a}^
{14}{b}^{4}-5232123109\,{K}^{11}{a}^{13}{b}^{5}-233383274114\,{K}^{11}
{a}^{12}{b}^{6}-841544583216\,{K}^{11}{a}^{11}{b}^{7}-686134498710\,{K
}^{11}{a}^{10}{b}^{8}+\\
&2560697228941\,{K}^{11}{a}^{9}{b}^{9}+
10107330115159\,{K}^{11}{a}^{8}{b}^{10}+16943634262468\,{K}^{11}{a}^{7
}{b}^{11}+13904030519727\,{K}^{11}{a}^{6}{b}^{12}+\\
&6004362513241\,{K}^{
11}{a}^{5}{b}^{13}+1356340075268\,{K}^{11}{a}^{4}{b}^{14}+139665966288
\,{K}^{11}{a}^{3}{b}^{15}+2878073223\,{K}^{11}{a}^{2}{b}^{16}-
280147693\,{K}^{11}a{b}^{17}-\\
&4433715\,{K}^{11}{b}^{18}+263692056\,{K}^
{10}{a}^{15}{b}^{4}-2305542749\,{K}^{10}{a}^{14}{b}^{5}-17846922269\,{
K}^{10}{a}^{13}{b}^{6}+3914270076\,{K}^{10}{a}^{12}{b}^{7}+\\
&
230788654125\,{K}^{10}{a}^{11}{b}^{8}+650011642748\,{K}^{10}{a}^{10}{b
}^{9}+773784608800\,{K}^{10}{a}^{9}{b}^{10}-28544883981\,{K}^{10}{a}^{
8}{b}^{11}-1378800154125\,{K}^{10}{a}^{7}{b}^{12}-\\
&1724160213904\,{K}^{
10}{a}^{6}{b}^{13}-1016942812306\,{K}^{10}{a}^{5}{b}^{14}-331144447385
\,{K}^{10}{a}^{4}{b}^{15}-59995407176\,{K}^{10}{a}^{3}{b}^{16}-
5566660529\,{K}^{10}{a}^{2}{b}^{17}-\\
&209307940\,{K}^{10}a{b}^{18}-
1477905\,{K}^{10}{b}^{19}+8186722346\,{K}^{14}{a}^{11}{b}^{3}-
23418962733\,{K}^{14}{a}^{10}{b}^{4}-1611264978515\,{K}^{14}{a}^{9}{b}
^{5}-\\
&9278969969422\,{K}^{14}{a}^{8}{b}^{6}-15406452354550\,{K}^{14}{a}
^{7}{b}^{7}-1790771351586\,{K}^{14}{a}^{6}{b}^{8}+2773060348534\,{K}^{
14}{a}^{5}{b}^{9}+\\
&416986240965\,{K}^{14}{a}^{4}{b}^{10}-85995166290\,{
K}^{14}{a}^{3}{b}^{11}-3954557611\,{K}^{14}{a}^{2}{b}^{12}+229696728\,
{K}^{14}a{b}^{13}+3412682\,{K}^{14}{b}^{14}+
\end{aligned}
\end{scriptsize}
$$
$$
\begin{scriptsize}
\begin{aligned}
\qquad &81973128457\,{K}^{13}{a}^{
12}{b}^{3}+1125346817056\,{K}^{13}{a}^{11}{b}^{4}+6087140431497\,{K}^{
13}{a}^{10}{b}^{5}+22650236507688\,{K}^{13}{a}^{9}{b}^{6}+\\
&
66855435709993\,{K}^{13}{a}^{8}{b}^{7}+103591224029566\,{K}^{13}{a}^{7
}{b}^{8}+46251562642798\,{K}^{13}{a}^{6}{b}^{9}-618200336573\,{K}^{13}
{a}^{5}{b}^{10}-\\
&3534474620094\,{K}^{13}{a}^{4}{b}^{11}-309445679352\,{
K}^{13}{a}^{3}{b}^{12}+35105861795\,{K}^{13}{a}^{2}{b}^{13}+1155035022
\,{K}^{13}a{b}^{14}-25476806\,{K}^{13}{b}^{15}+\\
&67855108749\,{K}^{12}{a
}^{13}{b}^{3}+422619642450\,{K}^{12}{a}^{12}{b}^{4}-1460202901386\,{K}
^{12}{a}^{11}{b}^{5}-17537408327281\,{K}^{12}{a}^{10}{b}^{6}-\\
&
65520754533457\,{K}^{12}{a}^{9}{b}^{7}-160188024383486\,{K}^{12}{a}^{8
}{b}^{8}-236723458980588\,{K}^{12}{a}^{7}{b}^{9}-153086206360753\,{K}^
{12}{a}^{6}{b}^{10}-\\
&39468501733430\,{K}^{12}{a}^{5}{b}^{11}-
854849645158\,{K}^{12}{a}^{4}{b}^{12}+1052295414152\,{K}^{12}{a}^{3}{b
}^{13}+90131926552\,{K}^{12}{a}^{2}{b}^{14}-2325978232\,{K}^{12}a{b}^{
15}-\\
&77654540\,{K}^{12}{b}^{16}+21008972482\,{K}^{11}{a}^{14}{b}^{3}-
66106180497\,{K}^{11}{a}^{13}{b}^{4}-1827452822051\,{K}^{11}{a}^{12}{b
}^{5}-6295513439691\,{K}^{11}{a}^{11}{b}^{6}-\\
&1997222542511\,{K}^{11}{a
}^{10}{b}^{7}+34021931927349\,{K}^{11}{a}^{9}{b}^{8}+115774370012772\,
{K}^{11}{a}^{8}{b}^{9}+193628172115795\,{K}^{11}{a}^{7}{b}^{10}+\\
&
157460282836813\,{K}^{11}{a}^{6}{b}^{11}+65130927553865\,{K}^{11}{a}^{
5}{b}^{12}+13055979168008\,{K}^{11}{a}^{4}{b}^{13}+865732150888\,{K}^{
11}{a}^{3}{b}^{14}-\\
&60852459218\,{K}^{11}{a}^{2}{b}^{15}-6974371057\,{K
}^{11}a{b}^{16}-26219767\,{K}^{11}{b}^{17}+2858370047\,{K}^{10}{a}^{15
}{b}^{3}-35976884599\,{K}^{10}{a}^{14}{b}^{4}-\\
&251387922082\,{K}^{10}{a
}^{13}{b}^{5}+361760982114\,{K}^{10}{a}^{12}{b}^{6}+4981123948124\,{K}
^{10}{a}^{11}{b}^{7}+12633563354145\,{K}^{10}{a}^{10}{b}^{8}+\\
&
12584574043677\,{K}^{10}{a}^{9}{b}^{9}-8259273138683\,{K}^{10}{a}^{8}{
b}^{10}-40808881281669\,{K}^{10}{a}^{7}{b}^{11}-46767928086816\,{K}^{
10}{a}^{6}{b}^{12}-\\
&26449351650057\,{K}^{10}{a}^{5}{b}^{13}-
8198628321896\,{K}^{10}{a}^{4}{b}^{14}-1365011957704\,{K}^{10}{a}^{3}{
b}^{15}-105070216072\,{K}^{10}{a}^{2}{b}^{16}-1899838528\,{K}^{10}a{b}
^{17}\\
&+61553452\,{K}^{10}{b}^{18}+134201703\,{K}^{9}{a}^{16}{b}^{3}-
3022953480\,{K}^{9}{a}^{15}{b}^{4}-3326954271\,{K}^{9}{a}^{14}{b}^{5}+
97874034062\,{K}^{9}{a}^{13}{b}^{6}+\\
&323785037364\,{K}^{9}{a}^{12}{b}^{
7}+32430971515\,{K}^{9}{a}^{11}{b}^{8}-1586165082958\,{K}^{9}{a}^{10}{
b}^{9}-3799021992098\,{K}^{9}{a}^{9}{b}^{10}-3935350218551\,{K}^{9}{a}
^{8}{b}^{11}-\\
&830081599694\,{K}^{9}{a}^{7}{b}^{12}+1907173041898\,{K}^{
9}{a}^{6}{b}^{13}+1968753527342\,{K}^{9}{a}^{5}{b}^{14}+890329955990\,
{K}^{9}{a}^{4}{b}^{15}+219529606096\,{K}^{9}{a}^{3}{b}^{16}+\\
&
29405864954\,{K}^{9}{a}^{2}{b}^{17}+1883930965\,{K}^{9}a{b}^{18}+
39008167\,{K}^{9}{b}^{19}+7908343648\,{K}^{14}{a}^{11}{b}^{2}-
81806177066\,{K}^{14}{a}^{10}{b}^{3}-\\
&3437749874761\,{K}^{14}{a}^{9}{b}
^{4}-21670560688842\,{K}^{14}{a}^{8}{b}^{5}-37740603940762\,{K}^{14}{a
}^{7}{b}^{6}+2438194984096\,{K}^{14}{a}^{6}{b}^{7}+\\
&11153028352916\,{K}
^{14}{a}^{5}{b}^{8}+850989159768\,{K}^{14}{a}^{4}{b}^{9}-499750328669
\,{K}^{14}{a}^{3}{b}^{10}-14876483127\,{K}^{14}{a}^{2}{b}^{11}+
1684533640\,{K}^{14}a{b}^{12}+\\
&69845411\,{K}^{14}{b}^{13}+123358247996
\,{K}^{13}{a}^{12}{b}^{2}+2122085110746\,{K}^{13}{a}^{11}{b}^{3}+
13440000815788\,{K}^{13}{a}^{10}{b}^{4}+58120917116440\,{K}^{13}{a}^{9
}{b}^{5}+\\
&199930974940797\,{K}^{13}{a}^{8}{b}^{6}+338266864358876\,{K}^
{13}{a}^{7}{b}^{7}+127094098448576\,{K}^{13}{a}^{6}{b}^{8}-
28928179154424\,{K}^{13}{a}^{5}{b}^{9}-\\
&18202154868479\,{K}^{13}{a}^{4}
{b}^{10}-640724597485\,{K}^{13}{a}^{3}{b}^{11}+283088743162\,{K}^{13}{
a}^{2}{b}^{12}+5382435299\,{K}^{13}a{b}^{13}-225585704\,{K}^{13}{b}^{
14}+\\
&130200374019\,{K}^{12}{a}^{13}{b}^{2}+896027322312\,{K}^{12}{a}^{
12}{b}^{3}-5775350019921\,{K}^{12}{a}^{11}{b}^{4}-63735890609267\,{K}^
{12}{a}^{10}{b}^{5}-\\
&249947721911140\,{K}^{12}{a}^{9}{b}^{6}-
665030463045448\,{K}^{12}{a}^{8}{b}^{7}-1060273602682070\,{K}^{12}{a}^
{7}{b}^{8}-667051103974848\,{K}^{12}{a}^{6}{b}^{9}-\\
&132926058451814\,{K
}^{12}{a}^{5}{b}^{10}+18175248648586\,{K}^{12}{a}^{4}{b}^{11}+
8135528030811\,{K}^{12}{a}^{3}{b}^{12}+374773360510\,{K}^{12}{a}^{2}{b
}^{13}-\\
&36928362338\,{K}^{12}a{b}^{14}-613706829\,{K}^{12}{b}^{15}+
54029957207\,{K}^{11}{a}^{14}{b}^{2}-331262807997\,{K}^{11}{a}^{13}{b}
^{3}-7696714243352\,{K}^{11}{a}^{12}{b}^{4}-\\
&25669087278272\,{K}^{11}{a
}^{11}{b}^{5}+11399205267054\,{K}^{11}{a}^{10}{b}^{6}+239789396169294
\,{K}^{11}{a}^{9}{b}^{7}+768738764880098\,{K}^{11}{a}^{8}{b}^{8}+\\
&
1317698051315836\,{K}^{11}{a}^{7}{b}^{9}+1066875566299538\,{K}^{11}{a}
^{6}{b}^{10}+416917615833117\,{K}^{11}{a}^{5}{b}^{11}+68718244721814\,
{K}^{11}{a}^{4}{b}^{12}-\\
&45474887656\,{K}^{11}{a}^{3}{b}^{13}-
1100413818399\,{K}^{11}{a}^{2}{b}^{14}-63159615150\,{K}^{11}a{b}^{15}+
655701568\,{K}^{11}{b}^{16}+11270438209\,{K}^{10}{a}^{15}{b}^{2}-\\
&
214415643170\,{K}^{10}{a}^{14}{b}^{3}-1387868573929\,{K}^{10}{a}^{13}{
b}^{4}+4748596841060\,{K}^{10}{a}^{12}{b}^{5}+44186660894713\,{K}^{10}
{a}^{11}{b}^{6}+\\
&103004239627341\,{K}^{10}{a}^{10}{b}^{7}+
77261052904064\,{K}^{10}{a}^{9}{b}^{8}-160182144215907\,{K}^{10}{a}^{8
}{b}^{9}-520442507701807\,{K}^{10}{a}^{7}{b}^{10}-\\
&565633913540164\,{K}
^{10}{a}^{6}{b}^{11}-308523872791728\,{K}^{10}{a}^{5}{b}^{12}-
90367829132546\,{K}^{10}{a}^{4}{b}^{13}-13365849495985\,{K}^{10}{a}^{3
}{b}^{14}-\\
&710600830431\,{K}^{10}{a}^{2}{b}^{15}+21233754939\,{K}^{10}a
{b}^{16}+1737451051\,{K}^{10}{b}^{17}+1128989577\,{K}^{9}{a}^{16}{b}^{
2}-35763166670\,{K}^{9}{a}^{15}{b}^{3}+\\
&9949826083\,{K}^{9}{a}^{14}{b}^
{4}+1666635376306\,{K}^{9}{a}^{13}{b}^{5}+4627472143796\,{K}^{9}{a}^{
12}{b}^{6}-4214041096869\,{K}^{9}{a}^{11}{b}^{7}-39078033416597\,{K}^{
9}{a}^{10}{b}^{8}-\\
&81130798384666\,{K}^{9}{a}^{9}{b}^{9}-73674805220287
\,{K}^{9}{a}^{8}{b}^{10}+5586523317171\,{K}^{9}{a}^{7}{b}^{11}+
66477627072877\,{K}^{9}{a}^{6}{b}^{12}+\\
&57740094344503\,{K}^{9}{a}^{5}{
b}^{13}+24438363853061\,{K}^{9}{a}^{4}{b}^{14}+5701167218956\,{K}^{9}{
a}^{3}{b}^{15}+709582428478\,{K}^{9}{a}^{2}{b}^{16}+39670256615\,{K}^{
9}a{b}^{17}+\\
&539256106\,{K}^{9}{b}^{18}+36449472\,{K}^{8}{a}^{17}{b}^{2
}-1814408259\,{K}^{8}{a}^{16}{b}^{3}+10444703251\,{K}^{8}{a}^{15}{b}^{
4}+84612651182\,{K}^{8}{a}^{14}{b}^{5}-\\
&61408548310\,{K}^{8}{a}^{13}{b}
^{6}-1184282204628\,{K}^{8}{a}^{12}{b}^{7}-2663411958703\,{K}^{8}{a}^{
11}{b}^{8}-1459642600038\,{K}^{8}{a}^{10}{b}^{9}+4207090872254\,{K}^{8
}{a}^{9}{b}^{10}+\\
&9829633012420\,{K}^{8}{a}^{8}{b}^{11}+8499274828070\,
{K}^{8}{a}^{7}{b}^{12}+2790709141625\,{K}^{8}{a}^{6}{b}^{13}-
658790024888\,{K}^{8}{a}^{5}{b}^{14}-923476667307\,{K}^{8}{a}^{4}{b}^{
15}-\\
&351721082422\,{K}^{8}{a}^{3}{b}^{16}-66806148045\,{K}^{8}{a}^{2}{b
}^{17}-6287685055\,{K}^{8}a{b}^{18}-224217663\,{K}^{8}{b}^{19}+
3184423556\,{K}^{14}{a}^{11}b-\\
&149305254732\,{K}^{14}{a}^{10}{b}^{2}-
5108388303854\,{K}^{14}{a}^{9}{b}^{3}-36133337768679\,{K}^{14}{a}^{8}{
b}^{4}-66103825069338\,{K}^{14}{a}^{7}{b}^{5}+\\
&24580631042072\,{K}^{14}
{a}^{6}{b}^{6}+30917882042366\,{K}^{14}{a}^{5}{b}^{7}-563005818808\,{K
}^{14}{a}^{4}{b}^{8}-2042132140784\,{K}^{14}{a}^{3}{b}^{9}-21877525669
\,{K}^{14}{a}^{2}{b}^{10}\\
&+10982950994\,{K}^{14}a{b}^{11}+637374812\,{K
}^{14}{b}^{12}+101533289300\,{K}^{13}{a}^{12}b+2468518923252\,{K}^{13}
{a}^{11}{b}^{2}+19645081280936\,{K}^{13}{a}^{10}{b}^{3}+\\
&
104256435812040\,{K}^{13}{a}^{9}{b}^{4}+431621919901755\,{K}^{13}{a}^{
8}{b}^{5}+803649493132164\,{K}^{13}{a}^{7}{b}^{6}+206914158576816\,{K}
^{13}{a}^{6}{b}^{7}-\\
&165670708718152\,{K}^{13}{a}^{5}{b}^{8}-
61861611640964\,{K}^{13}{a}^{4}{b}^{9}+2114584512813\,{K}^{13}{a}^{3}{
b}^{10}+1502746239563\,{K}^{13}{a}^{2}{b}^{11}+\\
&13440664057\,{K}^{13}a{
b}^{12}-1640392690\,{K}^{13}{b}^{13}+132251393870\,{K}^{12}{a}^{13}b+
1106608375499\,{K}^{12}{a}^{12}{b}^{2}-13289269965768\,{K}^{12}{a}^{11
}{b}^{3}-
\end{aligned}
\end{scriptsize}
$$
$$
\begin{scriptsize}
\begin{aligned}
\qquad &149857991866957\,{K}^{12}{a}^{10}{b}^{4}-641265672041080\,{K}
^{12}{a}^{9}{b}^{5}-1922880137781898\,{K}^{12}{a}^{8}{b}^{6}-
3371704702510544\,{K}^{12}{a}^{7}{b}^{7}-\\
&2006984668499370\,{K}^{12}{a}
^{6}{b}^{8}-205803185415322\,{K}^{12}{a}^{5}{b}^{9}+156456660211060\,{
K}^{12}{a}^{4}{b}^{10}+37758127417531\,{K}^{12}{a}^{3}{b}^{11}+\\
&
188449106462\,{K}^{12}{a}^{2}{b}^{12}-287808442945\,{K}^{12}a{b}^{13}-
2912043318\,{K}^{12}{b}^{14}+68650137499\,{K}^{11}{a}^{14}b-
822980538068\,{K}^{11}{a}^{13}{b}^{2}-\\
&18929167889770\,{K}^{11}{a}^{12}
{b}^{3}-61727589646732\,{K}^{11}{a}^{11}{b}^{4}+108879126750776\,{K}^{
11}{a}^{10}{b}^{5}+1032679153527805\,{K}^{11}{a}^{9}{b}^{6}+\\
&
3307689165757705\,{K}^{11}{a}^{8}{b}^{7}+5964984163315722\,{K}^{11}{a}
^{7}{b}^{8}+4809223142000353\,{K}^{11}{a}^{6}{b}^{9}+1723907426137727
\,{K}^{11}{a}^{5}{b}^{10}+\\
&186831230313931\,{K}^{11}{a}^{4}{b}^{11}-
34952024239795\,{K}^{11}{a}^{3}{b}^{12}-8395307721920\,{K}^{11}{a}^{2}
{b}^{13}-274159856016\,{K}^{11}a{b}^{14}+\\
&11646788639\,{K}^{11}{b}^{15}
+19290102976\,{K}^{10}{a}^{15}b-618051914572\,{K}^{10}{a}^{14}{b}^{2}-
3794937784271\,{K}^{10}{a}^{13}{b}^{3}+27137808072996\,{K}^{10}{a}^{12
}{b}^{4}+\\
&212703854630908\,{K}^{10}{a}^{11}{b}^{5}+458312813210486\,{K}
^{10}{a}^{10}{b}^{6}+175612669112832\,{K}^{10}{a}^{9}{b}^{7}-
1407151915680044\,{K}^{10}{a}^{8}{b}^{8}-\\
&3842435792013398\,{K}^{10}{a}
^{7}{b}^{9}-4049385673364058\,{K}^{10}{a}^{6}{b}^{10}-2131793811363353
\,{K}^{10}{a}^{5}{b}^{11}-579913882979471\,{K}^{10}{a}^{4}{b}^{12}-\\
&
70501554776744\,{K}^{10}{a}^{3}{b}^{13}-612252037236\,{K}^{10}{a}^{2}{
b}^{14}+480955683658\,{K}^{10}a{b}^{15}+17806369871\,{K}^{10}{b}^{16}+
3053066878\,{K}^{9}{a}^{16}b-\\
&153010661471\,{K}^{9}{a}^{15}{b}^{2}+
364828413706\,{K}^{9}{a}^{14}{b}^{3}+11111309092358\,{K}^{9}{a}^{13}{b
}^{4}+24829307259346\,{K}^{9}{a}^{12}{b}^{5}-\\
&71923908317843\,{K}^{9}{a
}^{11}{b}^{6}-388793461656790\,{K}^{9}{a}^{10}{b}^{7}-729747481416026
\,{K}^{9}{a}^{9}{b}^{8}-556617902197867\,{K}^{9}{a}^{8}{b}^{9}+\\
&
345469437963777\,{K}^{9}{a}^{7}{b}^{10}+955773539226550\,{K}^{9}{a}^{6
}{b}^{11}+745602357043282\,{K}^{9}{a}^{5}{b}^{12}+298436362767623\,{K}
^{9}{a}^{4}{b}^{13}+\\
&65591422646745\,{K}^{9}{a}^{3}{b}^{14}+
7411028447246\,{K}^{9}{a}^{2}{b}^{15}+327717764433\,{K}^{9}a{b}^{16}-
225405999\,{K}^{9}{b}^{17}+217227834\,{K}^{8}{a}^{17}b-\\
&16312007569\,{K
}^{8}{a}^{16}{b}^{2}+151454659051\,{K}^{8}{a}^{15}{b}^{3}+974170374900
\,{K}^{8}{a}^{14}{b}^{4}-2529274516037\,{K}^{8}{a}^{13}{b}^{5}-
21933270549165\,{K}^{8}{a}^{12}{b}^{6}-\\
&40192444297599\,{K}^{8}{a}^{11}
{b}^{7}+1967352915714\,{K}^{8}{a}^{10}{b}^{8}+128361755938788\,{K}^{8}
{a}^{9}{b}^{9}+235011662367697\,{K}^{8}{a}^{8}{b}^{10}+\\
&176164041982608
\,{K}^{8}{a}^{7}{b}^{11}+32278066307968\,{K}^{8}{a}^{6}{b}^{12}-
38758166923258\,{K}^{8}{a}^{5}{b}^{13}-31211816542749\,{K}^{8}{a}^{4}{
b}^{14}-\\
&10556861947296\,{K}^{8}{a}^{3}{b}^{15}-1873921010612\,{K}^{8}{
a}^{2}{b}^{16}-165770984303\,{K}^{8}a{b}^{17}-5448509151\,{K}^{8}{b}^{
18}+4138335\,{K}^{7}{a}^{18}b-\\
&529222631\,{K}^{7}{a}^{17}{b}^{2}+
9534299332\,{K}^{7}{a}^{16}{b}^{3}+8744650524\,{K}^{7}{a}^{15}{b}^{4}-
287794398492\,{K}^{7}{a}^{14}{b}^{5}-829389057239\,{K}^{7}{a}^{13}{b}^
{6}+\\
&436666511291\,{K}^{7}{a}^{12}{b}^{7}+5535592534957\,{K}^{7}{a}^{11
}{b}^{8}+10990952579700\,{K}^{7}{a}^{10}{b}^{9}+8694928776403\,{K}^{7}
{a}^{9}{b}^{10}-2215364207542\,{K}^{7}{a}^{8}{b}^{11}\\
&-9810055820574\,{
K}^{7}{a}^{7}{b}^{12}-8035296747722\,{K}^{7}{a}^{6}{b}^{13}-
3074004167155\,{K}^{7}{a}^{5}{b}^{14}-428798038681\,{K}^{7}{a}^{4}{b}^
{15}+90764751961\,{K}^{7}{a}^{3}{b}^{16}\\
&+45842745167\,{K}^{7}{a}^{2}{b
}^{17}+6902270186\,{K}^{7}a{b}^{18}+372400052\,{K}^{7}{b}^{19}-
137154360084\,{K}^{14}{a}^{10}b-4968549151736\,{K}^{14}{a}^{9}{b}^{2}-\\
&
41727938288138\,{K}^{14}{a}^{8}{b}^{3}-80542004636918\,{K}^{14}{a}^{7}
{b}^{4}+75838153694900\,{K}^{14}{a}^{6}{b}^{5}+58512932932734\,{K}^{14
}{a}^{5}{b}^{6}-\\
&10207065668376\,{K}^{14}{a}^{4}{b}^{7}-5920093710218\,
{K}^{14}{a}^{3}{b}^{8}+128874157835\,{K}^{14}{a}^{2}{b}^{9}+
56466275698\,{K}^{14}a{b}^{10}+3372848337\,{K}^{14}{b}^{11}+\\
&
35240176304\,{K}^{13}{a}^{12}+1621172890224\,{K}^{13}{a}^{11}b+
18359480852396\,{K}^{13}{a}^{10}{b}^{2}+129206634381942\,{K}^{13}{a}^{
9}{b}^{3}+\\
&670343237044533\,{K}^{13}{a}^{8}{b}^{4}+1386189982085088\,{K
}^{13}{a}^{7}{b}^{5}+72954287650256\,{K}^{13}{a}^{6}{b}^{6}-
548012781374044\,{K}^{13}{a}^{5}{b}^{7}-\\
&135882473739226\,{K}^{13}{a}^{
4}{b}^{8}+22665077968282\,{K}^{13}{a}^{3}{b}^{9}+5610536484497\,{K}^{
13}{a}^{2}{b}^{10}-28155736123\,{K}^{13}a{b}^{11}-10973587320\,{K}^{13
}{b}^{12}+\\
&55680988512\,{K}^{12}{a}^{13}+732524948886\,{K}^{12}{a}^{12}
b-18014141298174\,{K}^{12}{a}^{11}{b}^{2}-228941178837870\,{K}^{12}{a}
^{10}{b}^{3}-\\
&1122540651768463\,{K}^{12}{a}^{9}{b}^{4}-3960859617752802
\,{K}^{12}{a}^{8}{b}^{5}-7795831974670394\,{K}^{12}{a}^{7}{b}^{6}-
4131595052975000\,{K}^{12}{a}^{6}{b}^{7}+\\
&335469147375104\,{K}^{12}{a}^
{5}{b}^{8}+688618638024828\,{K}^{12}{a}^{4}{b}^{9}+110066534545953\,{K
}^{12}{a}^{3}{b}^{10}-7063667700674\,{K}^{12}{a}^{2}{b}^{11}-\\
&
1467153573562\,{K}^{12}a{b}^{12}-7919352127\,{K}^{12}{b}^{13}+
34704712872\,{K}^{11}{a}^{14}-1014661375436\,{K}^{11}{a}^{13}b-
27258553730623\,{K}^{11}{a}^{12}{b}^{2}-\\
&87465080840340\,{K}^{11}{a}^{
11}{b}^{3}+389309963095134\,{K}^{11}{a}^{10}{b}^{4}+2891180567527341\,
{K}^{11}{a}^{9}{b}^{5}+9727840726746332\,{K}^{11}{a}^{8}{b}^{6}+\\
&
18971463985569484\,{K}^{11}{a}^{7}{b}^{7}+15120097005244352\,{K}^{11}{
a}^{6}{b}^{8}+4652728771142606\,{K}^{11}{a}^{5}{b}^{9}+7935754643638\,
{K}^{11}{a}^{4}{b}^{10}-\\
&260964207429384\,{K}^{11}{a}^{3}{b}^{11}-
37777297138001\,{K}^{11}{a}^{2}{b}^{12}-235725834511\,{K}^{11}a{b}^{13
}+94580704227\,{K}^{11}{b}^{14}+12088173212\,{K}^{10}{a}^{15}-\\
&
869400676882\,{K}^{10}{a}^{14}b-5234564483080\,{K}^{10}{a}^{13}{b}^{2}
+82642094368055\,{K}^{10}{a}^{12}{b}^{3}+610249066221783\,{K}^{10}{a}^
{11}{b}^{4}+\\
&1199638153546702\,{K}^{10}{a}^{10}{b}^{5}-322166187460065
\,{K}^{10}{a}^{9}{b}^{6}-7231757136619890\,{K}^{10}{a}^{8}{b}^{7}-
18472682571952828\,{K}^{10}{a}^{7}{b}^{8}-\\
&19155901829629422\,{K}^{10}{
a}^{6}{b}^{9}-9672650542511106\,{K}^{10}{a}^{5}{b}^{10}-
2352733921764997\,{K}^{10}{a}^{4}{b}^{11}-184172502649356\,{K}^{10}{a}
^{3}{b}^{12}+\\
&22844226469695\,{K}^{10}{a}^{2}{b}^{13}+3960450269698\,{K
}^{10}a{b}^{14}+96064837755\,{K}^{10}{b}^{15}+2704152732\,{K}^{9}{a}^{
16}-282783442869\,{K}^{9}{a}^{15}b+\\
&1812808291250\,{K}^{9}{a}^{14}{b}^{
2}+37372347671457\,{K}^{9}{a}^{13}{b}^{3}+57532245928867\,{K}^{9}{a}^{
12}{b}^{4}-487467696159861\,{K}^{9}{a}^{11}{b}^{5}-\\
&2088714049512858\,{
K}^{9}{a}^{10}{b}^{6}-3618461852447637\,{K}^{9}{a}^{9}{b}^{7}-
2036890862050872\,{K}^{9}{a}^{8}{b}^{8}+4094899631764359\,{K}^{9}{a}^{
7}{b}^{9}+\\
&7785389063664478\,{K}^{9}{a}^{6}{b}^{10}+5642368005380278\,{
K}^{9}{a}^{5}{b}^{11}+2142168484274233\,{K}^{9}{a}^{4}{b}^{12}+
437610639273969\,{K}^{9}{a}^{3}{b}^{13}+\\
&42584265682256\,{K}^{9}{a}^{2}
{b}^{14}+1023654223028\,{K}^{9}a{b}^{15}-57766868229\,{K}^{9}{b}^{16}+
301036356\,{K}^{8}{a}^{17}-47026489315\,{K}^{8}{a}^{16}b+\\
&764302957569
\,{K}^{8}{a}^{15}{b}^{2}+3896753915869\,{K}^{8}{a}^{14}{b}^{3}-
25368506291608\,{K}^{8}{a}^{13}{b}^{4}-159004603220094\,{K}^{8}{a}^{12
}{b}^{5}-\\
&221496407923719\,{K}^{8}{a}^{11}{b}^{6}+285862139326960\,{K}^
{8}{a}^{10}{b}^{7}+1488256455679034\,{K}^{8}{a}^{9}{b}^{8}+
2359199946333475\,{K}^{8}{a}^{8}{b}^{9}+\\
&1485763917976251\,{K}^{8}{a}^{
7}{b}^{10}-102477804006590\,{K}^{8}{a}^{6}{b}^{11}-707333999227617\,{K
}^{8}{a}^{5}{b}^{12}-448742084763560\,{K}^{8}{a}^{4}{b}^{13}-\\
&
139234926029893\,{K}^{8}{a}^{3}{b}^{14}-23213148220661\,{K}^{8}{a}^{2}
{b}^{15}-1913714942359\,{K}^{8}a{b}^{16}-56112741840\,{K}^{8}{b}^{17}+
15026492\,{K}^{7}{a}^{18}-
\end{aligned}
\end{scriptsize}
$$
$$
\begin{scriptsize}
\begin{aligned}
\qquad &3330650891\,{K}^{7}{a}^{17}b+94548704291\,{K
}^{7}{a}^{16}{b}^{2}-89895307737\,{K}^{7}{a}^{15}{b}^{3}-4201366552153
\,{K}^{7}{a}^{14}{b}^{4}-8617112925025\,{K}^{7}{a}^{13}{b}^{5}+\\
&
23103484460009\,{K}^{7}{a}^{12}{b}^{6}+115808205057689\,{K}^{7}{a}^{11
}{b}^{7}+188576504258332\,{K}^{7}{a}^{10}{b}^{8}+94084539811435\,{K}^{
7}{a}^{9}{b}^{9}-\\
&152024750229991\,{K}^{7}{a}^{8}{b}^{10}-
279312474974775\,{K}^{7}{a}^{7}{b}^{11}-190198038196173\,{K}^{7}{a}^{6
}{b}^{12}-57579842917955\,{K}^{7}{a}^{5}{b}^{13}-\\
&137347680635\,{K}^{7}
{a}^{4}{b}^{14}+5458868289623\,{K}^{7}{a}^{3}{b}^{15}+1697950750626\,{
K}^{7}{a}^{2}{b}^{16}+222977338916\,{K}^{7}a{b}^{17}+11138876354\,{K}^
{7}{b}^{18}-\\
&62979693\,{K}^{6}{a}^{18}b+3234386947\,{K}^{6}{a}^{17}{b}^
{2}-21179612990\,{K}^{6}{a}^{16}{b}^{3}-145745814802\,{K}^{6}{a}^{15}{
b}^{4}+163658229840\,{K}^{6}{a}^{14}{b}^{5}+\\
&2180091092192\,{K}^{6}{a}^
{13}{b}^{6}+4592090089990\,{K}^{6}{a}^{12}{b}^{7}+1712026492259\,{K}^{
6}{a}^{11}{b}^{8}-9422615169243\,{K}^{6}{a}^{10}{b}^{9}-\\
&20016176469785
\,{K}^{6}{a}^{9}{b}^{10}-17620636656661\,{K}^{6}{a}^{8}{b}^{11}-
6272706330282\,{K}^{6}{a}^{7}{b}^{12}+1380891533427\,{K}^{6}{a}^{6}{b}
^{13}+\\
&2339255518064\,{K}^{6}{a}^{5}{b}^{14}+1022900194683\,{K}^{6}{a}^
{4}{b}^{15}+224940754767\,{K}^{6}{a}^{3}{b}^{16}+24120761757\,{K}^{6}{
a}^{2}{b}^{17}+724739915\,{K}^{6}a{b}^{18}-\\
&41016217\,{K}^{6}{b}^{19}-
50531308480\,{K}^{14}{a}^{10}-2821685568804\,{K}^{14}{a}^{9}b-
31513378801532\,{K}^{14}{a}^{8}{b}^{2}-64588568915756\,{K}^{14}{a}^{7}
{b}^{3}+\\
&138579749895342\,{K}^{14}{a}^{6}{b}^{4}+70411322834802\,{K}^{
14}{a}^{5}{b}^{5}-36083423105196\,{K}^{14}{a}^{4}{b}^{6}-
11778103677589\,{K}^{14}{a}^{3}{b}^{7}+\\
&1024228207926\,{K}^{14}{a}^{2}{
b}^{8}+208675192125\,{K}^{14}a{b}^{9}+11142340601\,{K}^{14}{b}^{10}+
460739873392\,{K}^{13}{a}^{11}+9972648850152\,{K}^{13}{a}^{10}b+\\
&
105868199989516\,{K}^{13}{a}^{9}{b}^{2}+730565381291033\,{K}^{13}{a}^{
8}{b}^{3}+1701480951898336\,{K}^{13}{a}^{7}{b}^{4}-547235313584102\,{K
}^{13}{a}^{6}{b}^{5}-\\
&1194065459423552\,{K}^{13}{a}^{5}{b}^{6}-
156002287416429\,{K}^{13}{a}^{4}{b}^{7}+94574111005629\,{K}^{13}{a}^{3
}{b}^{8}+14494325390747\,{K}^{13}{a}^{2}{b}^{9}-\\
&495681391729\,{K}^{13}
a{b}^{10}-58809378689\,{K}^{13}{b}^{11}+197936999696\,{K}^{12}{a}^{12}
-13394278195872\,{K}^{12}{a}^{11}b-219749310979130\,{K}^{12}{a}^{10}{b
}^{2}-\\
&1325797832601407\,{K}^{12}{a}^{9}{b}^{3}-5825504170769009\,{K}^{
12}{a}^{8}{b}^{4}-13164655642040532\,{K}^{12}{a}^{7}{b}^{5}-
5341463033615334\,{K}^{12}{a}^{6}{b}^{6}+\\
&2844019586671302\,{K}^{12}{a}
^{5}{b}^{7}+1939299678494474\,{K}^{12}{a}^{4}{b}^{8}+171506229771836\,
{K}^{12}{a}^{3}{b}^{9}-46746249130013\,{K}^{12}{a}^{2}{b}^{10}-\\
&
5236276942358\,{K}^{12}a{b}^{11}+6320650157\,{K}^{12}{b}^{12}-
497604274216\,{K}^{11}{a}^{13}-21299351023273\,{K}^{11}{a}^{12}b-
67036121674344\,{K}^{11}{a}^{11}{b}^{2}+\\
&771463775013592\,{K}^{11}{a}^{
10}{b}^{3}+5377463906004063\,{K}^{11}{a}^{9}{b}^{4}+20033719379721480
\,{K}^{11}{a}^{8}{b}^{5}+43565342979914178\,{K}^{11}{a}^{7}{b}^{6}+\\
&
33603361264374378\,{K}^{11}{a}^{6}{b}^{7}+7392694075630402\,{K}^{11}{a
}^{5}{b}^{8}-2080199975372483\,{K}^{11}{a}^{4}{b}^{9}-1067553308666791
\,{K}^{11}{a}^{3}{b}^{10}-\\
&103166219689987\,{K}^{11}{a}^{2}{b}^{11}+
4320351829561\,{K}^{11}a{b}^{12}+494202108539\,{K}^{11}{b}^{13}-
480400205840\,{K}^{10}{a}^{14}-3111278168772\,{K}^{10}{a}^{13}b+\\
&
140110012871789\,{K}^{10}{a}^{12}{b}^{2}+1067656975977167\,{K}^{10}{a}
^{11}{b}^{3}+1813687367178376\,{K}^{10}{a}^{10}{b}^{4}-
3135035065525276\,{K}^{10}{a}^{9}{b}^{5}-\\
&23982060626990151\,{K}^{10}{a
}^{8}{b}^{6}-61353544399212576\,{K}^{10}{a}^{7}{b}^{7}-
63110719547246584\,{K}^{10}{a}^{6}{b}^{8}-30012833409703380\,{K}^{10}{
a}^{5}{b}^{9}-\\
&5947792266282533\,{K}^{10}{a}^{4}{b}^{10}+70871891959793
\,{K}^{10}{a}^{3}{b}^{11}+184990521937143\,{K}^{10}{a}^{2}{b}^{12}+
18913665317331\,{K}^{10}a{b}^{13}+\\
&253634969371\,{K}^{10}{b}^{14}-
191327033800\,{K}^{9}{a}^{15}+3497637416479\,{K}^{9}{a}^{14}b+
66855651522437\,{K}^{9}{a}^{13}{b}^{2}+26280559181030\,{K}^{9}{a}^{12}
{b}^{3}-\\
&1761587138378454\,{K}^{9}{a}^{11}{b}^{4}-6710669712566067\,{K}
^{9}{a}^{10}{b}^{5}-10772208838896729\,{K}^{9}{a}^{9}{b}^{6}-
2481083150752946\,{K}^{9}{a}^{8}{b}^{7}+\\
&25558188139580602\,{K}^{9}{a}^
{7}{b}^{8}+40778426455627381\,{K}^{9}{a}^{6}{b}^{9}+27952240445058689
\,{K}^{9}{a}^{5}{b}^{10}+10011485640331752\,{K}^{9}{a}^{4}{b}^{11}+\\
&
1841382772484480\,{K}^{9}{a}^{3}{b}^{12}+133099707909710\,{K}^{9}{a}^{
2}{b}^{13}-3372393272333\,{K}^{9}a{b}^{14}-599625492324\,{K}^{9}{b}^{
15}-44073097048\,{K}^{8}{a}^{16}+\\
&1624119054178\,{K}^{8}{a}^{15}b+
5924430671197\,{K}^{8}{a}^{14}{b}^{2}-111301242128600\,{K}^{8}{a}^{13}
{b}^{3}-584078203496881\,{K}^{8}{a}^{12}{b}^{4}-\\
&457833244829843\,{K}^{
8}{a}^{11}{b}^{5}+2670864282012514\,{K}^{8}{a}^{10}{b}^{6}+
9129556194928068\,{K}^{8}{a}^{9}{b}^{7}+13131216222926796\,{K}^{8}{a}^
{8}{b}^{8}+\\
&6222950633254828\,{K}^{8}{a}^{7}{b}^{9}-4069329104151875\,{
K}^{8}{a}^{6}{b}^{10}-6706271725577709\,{K}^{8}{a}^{5}{b}^{11}-
3703814134952727\,{K}^{8}{a}^{4}{b}^{12}-\\
&1070685697713437\,{K}^{8}{a}^
{3}{b}^{13}-167089941574188\,{K}^{8}{a}^{2}{b}^{14}-12554824012517\,{K
}^{8}a{b}^{15}-302336359126\,{K}^{8}{b}^{16}-4894701920\,{K}^{7}{a}^{
17}+\\
&298486596657\,{K}^{7}{a}^{16}b-1278805245788\,{K}^{7}{a}^{15}{b}^{
2}-22118903159513\,{K}^{7}{a}^{14}{b}^{3}-22249961745479\,{K}^{7}{a}^{
13}{b}^{4}+\\
&261142339281241\,{K}^{7}{a}^{12}{b}^{5}+929731277649429\,{K
}^{7}{a}^{11}{b}^{6}+1225839957611770\,{K}^{7}{a}^{10}{b}^{7}-
7931936642419\,{K}^{7}{a}^{9}{b}^{8}-\\
&2393031788253814\,{K}^{7}{a}^{8}{
b}^{9}-3211903742011160\,{K}^{7}{a}^{7}{b}^{10}-1840069229181703\,{K}^
{7}{a}^{6}{b}^{11}-350205954697131\,{K}^{7}{a}^{5}{b}^{12}+\\
&
144835387201941\,{K}^{7}{a}^{4}{b}^{13}+103537152930075\,{K}^{7}{a}^{3
}{b}^{14}+26145330211318\,{K}^{7}{a}^{2}{b}^{15}+3126366617657\,{K}^{7
}a{b}^{16}+\\
&146321603727\,{K}^{7}{b}^{17}-239367016\,{K}^{6}{a}^{18}+
21819296858\,{K}^{6}{a}^{17}b-259508909190\,{K}^{6}{a}^{16}{b}^{2}-
1246580539917\,{K}^{6}{a}^{15}{b}^{3}+\\
&5269741227796\,{K}^{6}{a}^{14}{b
}^{4}+34089680143641\,{K}^{6}{a}^{13}{b}^{5}+51053051501289\,{K}^{6}{a
}^{12}{b}^{6}-39322994324066\,{K}^{6}{a}^{11}{b}^{7}-\\
&256091743030663\,
{K}^{6}{a}^{10}{b}^{8}-409589675682442\,{K}^{6}{a}^{9}{b}^{9}-
285640900279431\,{K}^{6}{a}^{8}{b}^{10}-28266868958843\,{K}^{6}{a}^{7}
{b}^{11}+\\
&96301413528594\,{K}^{6}{a}^{6}{b}^{12}+73680755326038\,{K}^{6
}{a}^{5}{b}^{13}+25444972454228\,{K}^{6}{a}^{4}{b}^{14}+4336215540792
\,{K}^{6}{a}^{3}{b}^{15}+\\
&218615071229\,{K}^{6}{a}^{2}{b}^{16}-
32146227033\,{K}^{6}a{b}^{17}-3589796829\,{K}^{6}{b}^{18}+414611097\,{
K}^{5}{a}^{18}b-10274999835\,{K}^{5}{a}^{17}{b}^{2}+\\
&120610835\,{K}^{5}
{a}^{16}{b}^{3}+333248544656\,{K}^{5}{a}^{15}{b}^{4}+902845282456\,{K}
^{5}{a}^{14}{b}^{5}-509778413929\,{K}^{5}{a}^{13}{b}^{6}-6453913661296
\,{K}^{5}{a}^{12}{b}^{7}-\\
&13789941431262\,{K}^{5}{a}^{11}{b}^{8}-
13205006679481\,{K}^{5}{a}^{10}{b}^{9}-1617906055033\,{K}^{5}{a}^{9}{b
}^{10}+9512476766469\,{K}^{5}{a}^{8}{b}^{11}+\\
&10755237658918\,{K}^{5}{a
}^{7}{b}^{12}+5897951583654\,{K}^{5}{a}^{6}{b}^{13}+1856553357578\,{K}
^{5}{a}^{5}{b}^{14}+324079954224\,{K}^{5}{a}^{4}{b}^{15}+\\
&23478140233\,
{K}^{5}{a}^{3}{b}^{16}-833160515\,{K}^{5}{a}^{2}{b}^{17}-160193642\,{K
}^{5}a{b}^{18}+629829\,{K}^{5}{b}^{19}-705066966464\,{K}^{14}{a}^{9}-
13900932276460\,{K}^{14}{a}^{8}b\\
&-30434950068832\,{K}^{14}{a}^{7}{b}^{2
}+162229918358796\,{K}^{14}{a}^{6}{b}^{3}+40273736264110\,{K}^{14}{a}^
{5}{b}^{4}-73095855910236\,{K}^{14}{a}^{4}{b}^{5}-\\
&14456064484953\,{K}^
{14}{a}^{3}{b}^{6}+3657046209882\,{K}^{14}{a}^{2}{b}^{7}+531196569199
\,{K}^{14}a{b}^{8}+22257744897\,{K}^{14}{b}^{9}+2404213864208\,{K}^{13
}{a}^{10}+
\end{aligned}
\end{scriptsize}
$$
$$
\begin{scriptsize}
\begin{aligned}
\qquad &51722699044464\,{K}^{13}{a}^{9}b+529218105447632\,{K}^{13}{a
}^{8}{b}^{2}+1424694338777028\,{K}^{13}{a}^{7}{b}^{3}-1520458464437824
\,{K}^{13}{a}^{6}{b}^{4}-\\
&1728045966209152\,{K}^{13}{a}^{5}{b}^{5}+
65827135137999\,{K}^{13}{a}^{4}{b}^{6}+243359731324900\,{K}^{13}{a}^{3
}{b}^{7}+23189078027684\,{K}^{13}{a}^{2}{b}^{8}-\\
&2676444151753\,{K}^{13
}a{b}^{9}-226842124944\,{K}^{13}{b}^{10}-4228126423984\,{K}^{12}{a}^{
11}-120476350001832\,{K}^{12}{a}^{10}b-\\
&1012459338862255\,{K}^{12}{a}^{
9}{b}^{2}-6006720839156582\,{K}^{12}{a}^{8}{b}^{3}-16028707041015834\,
{K}^{12}{a}^{7}{b}^{4}-2685308693524342\,{K}^{12}{a}^{6}{b}^{5}+\\
&
8200720142235588\,{K}^{12}{a}^{5}{b}^{6}+3575509913060384\,{K}^{12}{a}
^{4}{b}^{7}-55290998748698\,{K}^{12}{a}^{3}{b}^{8}-169021374256057\,{K
}^{12}{a}^{2}{b}^{9}-\\
&12766127143490\,{K}^{12}a{b}^{10}+207865333348\,{
K}^{12}{b}^{11}-6961114276272\,{K}^{11}{a}^{12}-20312240806328\,{K}^{
11}{a}^{11}b+\\
&890469752830855\,{K}^{11}{a}^{10}{b}^{2}+6588454158677199
\,{K}^{11}{a}^{9}{b}^{3}+29009589329789864\,{K}^{11}{a}^{8}{b}^{4}+
72898610139265312\,{K}^{11}{a}^{7}{b}^{5}+\\
&51740151797958412\,{K}^{11}{
a}^{6}{b}^{6}+2528726611966056\,{K}^{11}{a}^{5}{b}^{7}-
9072616132124752\,{K}^{11}{a}^{4}{b}^{8}-2783181552276564\,{K}^{11}{a}
^{3}{b}^{9}-\\
&129098442897342\,{K}^{11}{a}^{2}{b}^{10}+29827996870975\,{
K}^{11}a{b}^{11}+1811630767080\,{K}^{11}{b}^{12}-333769827284\,{K}^{10
}{a}^{13}+\\
&124959878608274\,{K}^{10}{a}^{12}b+1105418870732087\,{K}^{10
}{a}^{11}{b}^{2}+1318048303860316\,{K}^{10}{a}^{10}{b}^{3}-
9042382223470068\,{K}^{10}{a}^{9}{b}^{4}-\\
&53555594064770185\,{K}^{10}{a
}^{8}{b}^{5}-145188811083781898\,{K}^{10}{a}^{7}{b}^{6}-
148102384287331002\,{K}^{10}{a}^{6}{b}^{7}-63506059964587446\,{K}^{10}
{a}^{5}{b}^{8}\\
&-7378095632118436\,{K}^{10}{a}^{4}{b}^{9}+
2525046257767373\,{K}^{10}{a}^{3}{b}^{10}+777623288108454\,{K}^{10}{a}
^{2}{b}^{11}+55928161091201\,{K}^{10}a{b}^{12}-\\
&175237596275\,{K}^{10}{
b}^{13}+2396679773404\,{K}^{9}{a}^{14}+59773343857971\,{K}^{9}{a}^{13}
b-128766379736466\,{K}^{9}{a}^{12}{b}^{2}-\\
&3680685339835972\,{K}^{9}{a}
^{11}{b}^{3}-13299878426005787\,{K}^{9}{a}^{10}{b}^{4}-
19359977322342343\,{K}^{9}{a}^{9}{b}^{5}+8824037075736828\,{K}^{9}{a}^
{8}{b}^{6}+\\
&100542606530609358\,{K}^{9}{a}^{7}{b}^{7}+
146235539984107614\,{K}^{9}{a}^{6}{b}^{8}+95423253411417435\,{K}^{9}{a
}^{5}{b}^{9}+31634950501498244\,{K}^{9}{a}^{4}{b}^{10}+\\
&
4836162252973567\,{K}^{9}{a}^{3}{b}^{11}+104420406447081\,{K}^{9}{a}^{
2}{b}^{12}-48698927448058\,{K}^{9}a{b}^{13}-3328311272599\,{K}^{9}{b}^
{14}+\\
&1244603122484\,{K}^{8}{a}^{15}+985798971664\,{K}^{8}{a}^{14}b-
243390761794607\,{K}^{8}{a}^{13}{b}^{2}-1150560788673712\,{K}^{8}{a}^{
12}{b}^{3}+\\
&383445995602464\,{K}^{8}{a}^{11}{b}^{4}+11742095034627488\,
{K}^{8}{a}^{10}{b}^{5}+33440678088265027\,{K}^{8}{a}^{9}{b}^{6}+
44630172755794683\,{K}^{8}{a}^{8}{b}^{7}+\\
&10032295343683980\,{K}^{8}{a}
^{7}{b}^{8}-34230098542741248\,{K}^{8}{a}^{6}{b}^{9}-39290503623069237
\,{K}^{8}{a}^{5}{b}^{10}-19735029846443613\,{K}^{8}{a}^{4}{b}^{11}-\\
&
5335180617989225\,{K}^{8}{a}^{3}{b}^{12}-767220165313105\,{K}^{8}{a}^{
2}{b}^{13}-49700956825817\,{K}^{8}a{b}^{14}-717563493139\,{K}^{8}{b}^{
15}+302714083332\,{K}^{7}{a}^{16}\\
&-4146905542495\,{K}^{7}{a}^{15}b-
51750855766725\,{K}^{7}{a}^{14}{b}^{2}+43298384783439\,{K}^{7}{a}^{13}
{b}^{3}+1290900010963874\,{K}^{7}{a}^{12}{b}^{4}+\\
&3740083442730037\,{K}
^{7}{a}^{11}{b}^{5}+3500702541907372\,{K}^{7}{a}^{10}{b}^{6}-
4690226605461151\,{K}^{7}{a}^{9}{b}^{7}-17965843377290868\,{K}^{7}{a}^
{8}{b}^{8}-\\
&20053895525106181\,{K}^{7}{a}^{7}{b}^{9}-9177918267457711\,
{K}^{7}{a}^{6}{b}^{10}+205655323559909\,{K}^{7}{a}^{5}{b}^{11}+
2152558185369836\,{K}^{7}{a}^{4}{b}^{12}+\\
&1022373582987439\,{K}^{7}{a}^
{3}{b}^{13}+228108542372100\,{K}^{7}{a}^{2}{b}^{14}+25324871922019\,{K
}^{7}a{b}^{15}+1112933941774\,{K}^{7}{b}^{16}+34625733648\,{K}^{6}{a}^
{17}\\
&-970798991951\,{K}^{6}{a}^{16}b-2597489672872\,{K}^{6}{a}^{15}{b}^
{2}+41317494911156\,{K}^{6}{a}^{14}{b}^{3}+188928393294581\,{K}^{6}{a}
^{13}{b}^{4}+\\
&131060677315968\,{K}^{6}{a}^{12}{b}^{5}-786095184230832\,
{K}^{6}{a}^{11}{b}^{6}-2468014341749097\,{K}^{6}{a}^{10}{b}^{7}-
3240596490801268\,{K}^{6}{a}^{9}{b}^{8}-\\
&1508398576948736\,{K}^{6}{a}^{
8}{b}^{9}+907525680302711\,{K}^{6}{a}^{7}{b}^{10}+1577854506839307\,{K
}^{6}{a}^{6}{b}^{11}+888754064285080\,{K}^{6}{a}^{5}{b}^{12}+\\
&
242460751203109\,{K}^{6}{a}^{4}{b}^{13}+23278245578285\,{K}^{6}{a}^{3}
{b}^{14}-3705137648126\,{K}^{6}{a}^{2}{b}^{15}-1101527020655\,{K}^{6}a
{b}^{16}-76246237494\,{K}^{6}{b}^{17}+\\
&1643378400\,{K}^{5}{a}^{18}-
76714156068\,{K}^{5}{a}^{17}b+234281860362\,{K}^{5}{a}^{16}{b}^{2}+
3946662148693\,{K}^{5}{a}^{15}{b}^{3}+6082447778789\,{K}^{5}{a}^{14}{b
}^{4}-\\
&27527586471853\,{K}^{5}{a}^{13}{b}^{5}-119962333970100\,{K}^{5}{
a}^{12}{b}^{6}-196248396610823\,{K}^{5}{a}^{11}{b}^{7}-112055660874378
\,{K}^{5}{a}^{10}{b}^{8}+\\
&131880891121946\,{K}^{5}{a}^{9}{b}^{9}+
292122114717556\,{K}^{5}{a}^{8}{b}^{10}+237982626942796\,{K}^{5}{a}^{7
}{b}^{11}+102319035910323\,{K}^{5}{a}^{6}{b}^{12}+\\
&21395652192352\,{K}^
{5}{a}^{5}{b}^{13}+158183151036\,{K}^{5}{a}^{4}{b}^{14}-892758402630\,
{K}^{5}{a}^{3}{b}^{15}-163862960209\,{K}^{5}{a}^{2}{b}^{16}-7836924831
\,{K}^{5}a{b}^{17}+\\
&260425239\,{K}^{5}{b}^{18}-1495052703\,{K}^{4}{a}^{
18}b+16395610065\,{K}^{4}{a}^{17}{b}^{2}+88285680725\,{K}^{4}{a}^{16}{
b}^{3}-137349431135\,{K}^{4}{a}^{15}{b}^{4}-\\
&1553588794285\,{K}^{4}{a}^
{14}{b}^{5}-3762863376807\,{K}^{4}{a}^{13}{b}^{6}-3234415344731\,{K}^{
4}{a}^{12}{b}^{7}+3643227620716\,{K}^{4}{a}^{11}{b}^{8}+\\
&13514815580354
\,{K}^{4}{a}^{10}{b}^{9}+17215107072525\,{K}^{4}{a}^{9}{b}^{10}+
12406061872078\,{K}^{4}{a}^{8}{b}^{11}+5457337590899\,{K}^{4}{a}^{7}{b
}^{12}+\\
&1423855032488\,{K}^{4}{a}^{6}{b}^{13}+182007629522\,{K}^{4}{a}^
{5}{b}^{14}-1285693446\,{K}^{4}{a}^{4}{b}^{15}-3006700650\,{K}^{4}{a}^
{3}{b}^{16}-213691776\,{K}^{4}{a}^{2}{b}^{17}+\\
&3245540\,{K}^{4}a{b}^{18
}+7177\,{K}^{4}{b}^{19}-2692936560384\,{K}^{14}{a}^{8}-6231984045528\,
{K}^{14}{a}^{7}b+119363225302312\,{K}^{14}{a}^{6}{b}^{2}-\\
&
16937901068108\,{K}^{14}{a}^{5}{b}^{3}-93207723735006\,{K}^{14}{a}^{4}
{b}^{4}-6404331427633\,{K}^{14}{a}^{3}{b}^{5}+8036636065172\,{K}^{14}{
a}^{2}{b}^{6}+888012078549\,{K}^{14}a{b}^{7}\\
&+19388755947\,{K}^{14}{b}^
{8}+11411787649872\,{K}^{13}{a}^{9}+227980041989484\,{K}^{13}{a}^{8}b+
751498136522120\,{K}^{13}{a}^{7}{b}^{2}-\\
&2093318110224008\,{K}^{13}{a}^
{6}{b}^{3}-1548622090345744\,{K}^{13}{a}^{5}{b}^{4}+599357310554502\,{
K}^{13}{a}^{4}{b}^{5}+408291033630290\,{K}^{13}{a}^{3}{b}^{6}+\\
&
11784706451286\,{K}^{13}{a}^{2}{b}^{7}-8645932651538\,{K}^{13}a{b}^{8}
-603731490172\,{K}^{13}{b}^{9}-28779257999760\,{K}^{12}{a}^{10}-
452008661948386\,{K}^{12}{a}^{9}b-\\
&4140541755580987\,{K}^{12}{a}^{8}{b}
^{2}-13607251899857144\,{K}^{12}{a}^{7}{b}^{3}+4201414789839054\,{K}^{
12}{a}^{6}{b}^{4}+14095327596202524\,{K}^{12}{a}^{5}{b}^{5}+\\
&
3905243646637698\,{K}^{12}{a}^{4}{b}^{6}-1008569279179356\,{K}^{12}{a}
^{3}{b}^{7}-396009326213810\,{K}^{12}{a}^{2}{b}^{8}-17961750240666\,{K
}^{12}a{b}^{9}+\\
&1211350018334\,{K}^{12}{b}^{10}+1294867135152\,{K}^{11}
{a}^{11}+561178892108111\,{K}^{11}{a}^{10}b+5092868830116724\,{K}^{11}
{a}^{9}{b}^{2}+\\
&29084069814987456\,{K}^{11}{a}^{8}{b}^{3}+
88269967150237616\,{K}^{11}{a}^{7}{b}^{4}+50679200978712249\,{K}^{11}{
a}^{6}{b}^{5}-18683969569201448\,{K}^{11}{a}^{5}{b}^{6}-\\
&
21748033000740663\,{K}^{11}{a}^{4}{b}^{7}-4519184617749551\,{K}^{11}{a
}^{3}{b}^{8}+198654617348963\,{K}^{11}{a}^{2}{b}^{9}+108556860717837\,
{K}^{11}a{b}^{10}+
\end{aligned}
\end{scriptsize}
$$
$$
\begin{scriptsize}
\begin{aligned}
\qquad &4653590255176\,{K}^{11}{b}^{11}+45684070929472\,{K}^
{10}{a}^{12}+609954852665572\,{K}^{10}{a}^{11}b-69666015311194\,{K}^{
10}{a}^{10}{b}^{2}-\\
&14147424714484267\,{K}^{10}{a}^{9}{b}^{3}-
81569884556963387\,{K}^{10}{a}^{8}{b}^{4}-247962202750152072\,{K}^{10}
{a}^{7}{b}^{5}-247266317159753160\,{K}^{10}{a}^{6}{b}^{6}-\\
&
85343224605194656\,{K}^{10}{a}^{5}{b}^{7}+6390370105968982\,{K}^{10}{a
}^{4}{b}^{8}+10440398812946386\,{K}^{10}{a}^{3}{b}^{9}+
2048700401926028\,{K}^{10}{a}^{2}{b}^{10}+\\
&88580757309689\,{K}^{10}a{b}
^{11}-4487532858348\,{K}^{10}{b}^{12}+20386010028392\,{K}^{9}{a}^{13}-
231865914389097\,{K}^{9}{a}^{12}b-\\
&4432083858294279\,{K}^{9}{a}^{11}{b}
^{2}-15952015861462018\,{K}^{9}{a}^{10}{b}^{3}-19155588791293200\,{K}^
{9}{a}^{9}{b}^{4}+46238658114553168\,{K}^{9}{a}^{8}{b}^{5}+\\
&
267664001032388354\,{K}^{9}{a}^{7}{b}^{6}+370334153147232244\,{K}^{9}{
a}^{6}{b}^{7}+228614863749578936\,{K}^{9}{a}^{5}{b}^{8}+
66812266193388843\,{K}^{9}{a}^{4}{b}^{9}+\\
&6391481580951535\,{K}^{9}{a}^
{3}{b}^{10}-984440408335642\,{K}^{9}{a}^{2}{b}^{11}-241232987182592\,{
K}^{9}a{b}^{12}-11677707367669\,{K}^{9}{b}^{13}-\\
&3479355836552\,{K}^{8}
{a}^{14}-257742458381510\,{K}^{8}{a}^{13}b-1143406434061603\,{K}^{8}{a
}^{12}{b}^{2}+3642066685006909\,{K}^{8}{a}^{11}{b}^{3}+\\
&
29201596500092107\,{K}^{8}{a}^{10}{b}^{4}+76596923865512760\,{K}^{8}{a
}^{9}{b}^{5}+95087254463205111\,{K}^{8}{a}^{8}{b}^{6}-
24736114821587130\,{K}^{8}{a}^{7}{b}^{7}-\\
&160059786921675078\,{K}^{8}{a
}^{6}{b}^{8}-154001291113408680\,{K}^{8}{a}^{5}{b}^{9}-
71687845876395345\,{K}^{8}{a}^{4}{b}^{10}-17956700855491082\,{K}^{8}{a
}^{3}{b}^{11}-\\
&2281182870950745\,{K}^{8}{a}^{2}{b}^{12}-106603816767819
\,{K}^{8}a{b}^{13}+1353752090555\,{K}^{8}{b}^{14}-4118509669360\,{K}^{
7}{a}^{15}-\\
&50725924298409\,{K}^{7}{a}^{14}b+294102756038628\,{K}^{7}{a
}^{13}{b}^{2}+3275438513430389\,{K}^{7}{a}^{12}{b}^{3}+
7942714436575110\,{K}^{7}{a}^{11}{b}^{4}+\\
&1920689711445251\,{K}^{7}{a}^
{10}{b}^{5}-30523001973679832\,{K}^{7}{a}^{9}{b}^{6}-78153924244013403
\,{K}^{7}{a}^{8}{b}^{7}-75238398189333222\,{K}^{7}{a}^{7}{b}^{8}-\\
&
22008490768208907\,{K}^{7}{a}^{6}{b}^{9}+14885317358832578\,{K}^{7}{a}
^{5}{b}^{10}+15825765968950595\,{K}^{7}{a}^{4}{b}^{11}+
6229134241259728\,{K}^{7}{a}^{3}{b}^{12}+\\
&1271473277237975\,{K}^{7}{a}^
{2}{b}^{13}+131773249589440\,{K}^{7}a{b}^{14}+5388151596071\,{K}^{7}{b
}^{15}-1130177466944\,{K}^{6}{a}^{16}+1659111043648\,{K}^{6}{a}^{15}b\\
&
+128763133458158\,{K}^{6}{a}^{14}{b}^{2}+447969305895330\,{K}^{6}{a}^{
13}{b}^{3}-436099149109350\,{K}^{6}{a}^{12}{b}^{4}-4927671655635497\,{
K}^{6}{a}^{11}{b}^{5}-\\
&11672475605896403\,{K}^{6}{a}^{10}{b}^{6}-
12498248625509607\,{K}^{6}{a}^{9}{b}^{7}-427568340711393\,{K}^{6}{a}^{
8}{b}^{8}+12125065702837560\,{K}^{6}{a}^{7}{b}^{9}+\\
&12279845992039297\,
{K}^{6}{a}^{6}{b}^{10}+5504155367580718\,{K}^{6}{a}^{5}{b}^{11}+
1003323832986341\,{K}^{6}{a}^{4}{b}^{12}-101955080750521\,{K}^{6}{a}^{
3}{b}^{13}-\\
&80830480891686\,{K}^{6}{a}^{2}{b}^{14}-13724706059495\,{K}^
{6}a{b}^{15}-805480953826\,{K}^{6}{b}^{16}-135191066488\,{K}^{5}{a}^{
17}+1539093492927\,{K}^{5}{a}^{16}b+\\
&14369624442333\,{K}^{5}{a}^{15}{b}
^{2}-5669255112152\,{K}^{5}{a}^{14}{b}^{3}-254359279060994\,{K}^{5}{a}
^{13}{b}^{4}-761809782836692\,{K}^{5}{a}^{12}{b}^{5}-\\
&890879122175846\,
{K}^{5}{a}^{11}{b}^{6}+253425994538859\,{K}^{5}{a}^{10}{b}^{7}+
2283349521672712\,{K}^{5}{a}^{9}{b}^{8}+3048414237654825\,{K}^{5}{a}^{
8}{b}^{9}+\\
&1926757638954412\,{K}^{5}{a}^{7}{b}^{10}+538423977328159\,{K
}^{5}{a}^{6}{b}^{11}-37927188854755\,{K}^{5}{a}^{5}{b}^{12}-
72978990077681\,{K}^{5}{a}^{4}{b}^{13}-\\
&20527221440877\,{K}^{5}{a}^{3}{
b}^{14}-2113399892306\,{K}^{5}{a}^{2}{b}^{15}+24195389120\,{K}^{5}a{b}
^{16}+13289468812\,{K}^{5}{b}^{17}-6243359208\,{K}^{4}{a}^{18}+\\
&
148028398292\,{K}^{4}{a}^{17}b+439150576578\,{K}^{4}{a}^{16}{b}^{2}-
3780997705599\,{K}^{4}{a}^{15}{b}^{3}-20309492641311\,{K}^{4}{a}^{14}{
b}^{4}-32832390773437\,{K}^{4}{a}^{13}{b}^{5}\\
&+11372526277032\,{K}^{4}{
a}^{12}{b}^{6}+141051823652163\,{K}^{4}{a}^{11}{b}^{7}+273246363726338
\,{K}^{4}{a}^{10}{b}^{8}+269704120041771\,{K}^{4}{a}^{9}{b}^{9}+\\
&
145149873977302\,{K}^{4}{a}^{8}{b}^{10}+33351415166770\,{K}^{4}{a}^{7}
{b}^{11}-6205343893158\,{K}^{4}{a}^{6}{b}^{12}-6295905948018\,{K}^{4}{
a}^{5}{b}^{13}-\\
&1673106158673\,{K}^{4}{a}^{4}{b}^{14}-175796970966\,{K}
^{4}{a}^{3}{b}^{15}-771509\,{K}^{4}{a}^{2}{b}^{16}+867726416\,{K}^{4}a
{b}^{17}-2998428\,{K}^{4}{b}^{18}+3187232664\,{K}^{3}{a}^{18}b-\\
&
6456736798\,{K}^{3}{a}^{17}{b}^{2}-139944160178\,{K}^{3}{a}^{16}{b}^{3
}-413265939786\,{K}^{3}{a}^{15}{b}^{4}-100680550135\,{K}^{3}{a}^{14}{b
}^{5}+2230528301031\,{K}^{3}{a}^{13}{b}^{6}+\\
&6863902786680\,{K}^{3}{a}^
{12}{b}^{7}+10880910451348\,{K}^{3}{a}^{11}{b}^{8}+10523219978368\,{K}
^{3}{a}^{10}{b}^{9}+6430277198065\,{K}^{3}{a}^{9}{b}^{10}+\\
&
2424900980901\,{K}^{3}{a}^{8}{b}^{11}+493884284695\,{K}^{3}{a}^{7}{b}^
{12}+18820916425\,{K}^{3}{a}^{6}{b}^{13}-13580903883\,{K}^{3}{a}^{5}{b
}^{14}-2548392493\,{K}^{3}{a}^{4}{b}^{15}-\\
&101212275\,{K}^{3}{a}^{3}{b}
^{16}+5640712\,{K}^{3}{a}^{2}{b}^{17}+21531\,{K}^{3}a{b}^{18}+
100752862464\,{K}^{14}{a}^{7}+50170250902968\,{K}^{14}{a}^{6}b-\\
&
47964137620560\,{K}^{14}{a}^{5}{b}^{2}-72433068519684\,{K}^{14}{a}^{4}
{b}^{3}+10520341979541\,{K}^{14}{a}^{3}{b}^{4}+11352459794342\,{K}^{14
}{a}^{2}{b}^{5}+\\
&851436293991\,{K}^{14}a{b}^{6}-24261081958\,{K}^{14}{b
}^{7}+43969169658336\,{K}^{13}{a}^{8}+212069196853728\,{K}^{13}{a}^{7}
b-1688270416699240\,{K}^{13}{a}^{6}{b}^{2}-\\
&637922294838036\,{K}^{13}{a
}^{5}{b}^{3}+1098475632107739\,{K}^{13}{a}^{4}{b}^{4}+424551695160696
\,{K}^{13}{a}^{3}{b}^{5}-40237559237720\,{K}^{13}{a}^{2}{b}^{6}-\\
&
18198906568872\,{K}^{13}a{b}^{7}-1066035732905\,{K}^{13}{b}^{8}-
89676720796528\,{K}^{12}{a}^{9}-1717788158236594\,{K}^{12}{a}^{8}b-\\
&
7543866870318356\,{K}^{12}{a}^{7}{b}^{2}+10231779167220804\,{K}^{12}{a
}^{6}{b}^{3}+15207925801267222\,{K}^{12}{a}^{5}{b}^{4}+
1183704361758389\,{K}^{12}{a}^{4}{b}^{5}-\\
&2623883549555371\,{K}^{12}{a}
^{3}{b}^{6}-598200318371162\,{K}^{12}{a}^{2}{b}^{7}-28668624359\,{K}^{
12}a{b}^{8}+4120965507754\,{K}^{12}{b}^{9}+\\
&149256643511192\,{K}^{11}{a
}^{10}+2237640946866444\,{K}^{11}{a}^{9}b+19341015141994697\,{K}^{11}{
a}^{8}{b}^{2}+75338438931591088\,{K}^{11}{a}^{7}{b}^{3}+\\
&
21294743961172398\,{K}^{11}{a}^{6}{b}^{4}-50816031654256163\,{K}^{11}{
a}^{5}{b}^{5}-32618237795799194\,{K}^{11}{a}^{4}{b}^{6}-
3363575679865404\,{K}^{11}{a}^{3}{b}^{7}+\\
&1326797011202182\,{K}^{11}{a}
^{2}{b}^{8}+254534501295101\,{K}^{11}a{b}^{9}+7595696006975\,{K}^{11}{
b}^{10}+132658747826380\,{K}^{10}{a}^{11}-\\
&738746676845934\,{K}^{10}{a}
^{10}b-12714979032421558\,{K}^{10}{a}^{9}{b}^{2}-83813734152140243\,{K
}^{10}{a}^{8}{b}^{3}-304881920016878426\,{K}^{10}{a}^{7}{b}^{4}-\\
&
284579362575182985\,{K}^{10}{a}^{6}{b}^{5}-49954067048226019\,{K}^{10}
{a}^{5}{b}^{6}+48383691311278268\,{K}^{10}{a}^{4}{b}^{7}+
24207730808591071\,{K}^{10}{a}^{3}{b}^{8}+\\
&3315552363578566\,{K}^{10}{a
}^{2}{b}^{9}-22211268185700\,{K}^{10}a{b}^{10}-20299827019460\,{K}^{10
}{b}^{11}-119066387592516\,{K}^{9}{a}^{12}-\\
&2821953969383025\,{K}^{9}{a
}^{11}b-10628461707386274\,{K}^{9}{a}^{10}{b}^{2}-5947281572604297\,{K
}^{9}{a}^{9}{b}^{3}+99231591366398213\,{K}^{9}{a}^{8}{b}^{4}+\\
&
497592939560442758\,{K}^{9}{a}^{7}{b}^{5}+668398116479733244\,{K}^{9}{
a}^{6}{b}^{6}+379678670190497190\,{K}^{9}{a}^{5}{b}^{7}+
84781726302150240\,{K}^{9}{a}^{4}{b}^{8}-\\
&4310311065285835\,{K}^{9}{a}^
{3}{b}^{9}-5035262013292274\,{K}^{9}{a}^{2}{b}^{10}-711797978466912\,{
K}^{9}a{b}^{11}-25934980128699\,{K}^{9}{b}^{12}-
\end{aligned}
\end{scriptsize}
$$
$$
\begin{scriptsize}
\begin{aligned}
\qquad &103305978421148\,{K}^{
8}{a}^{13}-421217608783212\,{K}^{8}{a}^{12}b+6947779140205101\,{K}^{8}
{a}^{11}{b}^{2}+42625672111203831\,{K}^{8}{a}^{10}{b}^{3}+\\
&
109918579021491339\,{K}^{8}{a}^{9}{b}^{4}+122789123721845370\,{K}^{8}{
a}^{8}{b}^{5}-178552741178374062\,{K}^{8}{a}^{7}{b}^{6}-
485033453321789604\,{K}^{8}{a}^{6}{b}^{7}-\\
&419131182012134250\,{K}^{8}{
a}^{5}{b}^{8}-180973004117895586\,{K}^{8}{a}^{4}{b}^{9}-
40641095885988951\,{K}^{8}{a}^{3}{b}^{10}-4014872362986644\,{K}^{8}{a}
^{2}{b}^{11}-\\
&8357591190600\,{K}^{8}a{b}^{12}+17470123297978\,{K}^{8}{b
}^{13}-12929794517612\,{K}^{7}{a}^{14}+467050356730421\,{K}^{7}{a}^{13
}b+\\
&4323207249022227\,{K}^{7}{a}^{12}{b}^{2}+7939456212265215\,{K}^{7}{
a}^{11}{b}^{3}-14511275048776223\,{K}^{7}{a}^{10}{b}^{4}-
96468616515317634\,{K}^{7}{a}^{9}{b}^{5}-\\
&212985246177841059\,{K}^{7}{a
}^{8}{b}^{6}-172308812349866262\,{K}^{7}{a}^{7}{b}^{7}+
4182673582118460\,{K}^{7}{a}^{6}{b}^{8}+98614741820444325\,{K}^{7}{a}^
{5}{b}^{9}+\\
&72100713016481131\,{K}^{7}{a}^{4}{b}^{10}+25266442719233159
\,{K}^{7}{a}^{3}{b}^{11}+4774427950066559\,{K}^{7}{a}^{2}{b}^{12}+
457731011567340\,{K}^{7}a{b}^{13}+\\
&16873852622759\,{K}^{7}{b}^{14}+
6538254640340\,{K}^{6}{a}^{15}+165619185143066\,{K}^{6}{a}^{14}b+
364571288087372\,{K}^{6}{a}^{13}{b}^{2}-\\
&2933225279941775\,{K}^{6}{a}^{
12}{b}^{3}-14903103012069741\,{K}^{6}{a}^{11}{b}^{4}-29547212218185812
\,{K}^{6}{a}^{10}{b}^{5}-21881728544446730\,{K}^{6}{a}^{9}{b}^{6}+\\
&
29067033215054568\,{K}^{6}{a}^{8}{b}^{7}+69409282617743016\,{K}^{6}{a}
^{7}{b}^{8}+54305493082636782\,{K}^{6}{a}^{6}{b}^{9}+18217401805391401
\,{K}^{6}{a}^{5}{b}^{10}+\\
&288590329144\,{K}^{6}{a}^{4}{b}^{11}-
2062976889285423\,{K}^{6}{a}^{3}{b}^{12}-686726467214185\,{K}^{6}{a}^{
2}{b}^{13}-96412645682295\,{K}^{6}a{b}^{14}-\\
&5153197131658\,{K}^{6}{b}^
{15}+2420277000500\,{K}^{5}{a}^{16}+15953711816837\,{K}^{5}{a}^{15}b-
109357648174616\,{K}^{5}{a}^{14}{b}^{2}-908623228824223\,{K}^{5}{a}^{
13}{b}^{3}\\
&-2039746857885538\,{K}^{5}{a}^{12}{b}^{4}-746379967413935\,{
K}^{5}{a}^{11}{b}^{5}+5905618256794211\,{K}^{5}{a}^{10}{b}^{6}+
14818752175887789\,{K}^{5}{a}^{9}{b}^{7}+\\
&15201813248327444\,{K}^{5}{a}
^{8}{b}^{8}+6607355644237153\,{K}^{5}{a}^{7}{b}^{9}-503172391575267\,{
K}^{5}{a}^{6}{b}^{10}-1840032767502222\,{K}^{5}{a}^{5}{b}^{11}-\\
&
812740841175985\,{K}^{5}{a}^{4}{b}^{12}-147606144356421\,{K}^{5}{a}^{3
}{b}^{13}-3786016989280\,{K}^{5}{a}^{2}{b}^{14}+2287945690950\,{K}^{5}
a{b}^{15}+224603891766\,{K}^{5}{b}^{16}\\
&+309386717100\,{K}^{4}{a}^{17}-
450716134502\,{K}^{4}{a}^{16}b-21590094627859\,{K}^{4}{a}^{15}{b}^{2}-
77280620228751\,{K}^{4}{a}^{14}{b}^{3}-26156343086583\,{K}^{4}{a}^{13}
{b}^{4}\\
&+419684077689451\,{K}^{4}{a}^{12}{b}^{5}+1273519040824233\,{K}^
{4}{a}^{11}{b}^{6}+1851736514494502\,{K}^{4}{a}^{10}{b}^{7}+
1334684555782330\,{K}^{4}{a}^{9}{b}^{8}+\\
&239957947410597\,{K}^{4}{a}^{8
}{b}^{9}-341163178159078\,{K}^{4}{a}^{7}{b}^{10}-287106476400073\,{K}^
{4}{a}^{6}{b}^{11}-96555090819376\,{K}^{4}{a}^{5}{b}^{12}-\\
&
12695298078414\,{K}^{4}{a}^{4}{b}^{13}+877747493270\,{K}^{4}{a}^{3}{b}
^{14}+411200864516\,{K}^{4}{a}^{2}{b}^{15}+27001285531\,{K}^{4}a{b}^{
16}-678753772\,{K}^{4}{b}^{17}+\\
&14247859828\,{K}^{3}{a}^{18}-
131870001809\,{K}^{3}{a}^{17}b-1229053278230\,{K}^{3}{a}^{16}{b}^{2}-
1663556797311\,{K}^{3}{a}^{15}{b}^{3}+8941598037433\,{K}^{3}{a}^{14}{b
}^{4}+\\
&43124915755610\,{K}^{3}{a}^{13}{b}^{5}+91964818350375\,{K}^{3}{a
}^{12}{b}^{6}+113124099635591\,{K}^{3}{a}^{11}{b}^{7}+74612016730936\,
{K}^{3}{a}^{10}{b}^{8}+\\
&13158860998236\,{K}^{3}{a}^{9}{b}^{9}-
17604258251288\,{K}^{3}{a}^{8}{b}^{10}-15246072102725\,{K}^{3}{a}^{7}{
b}^{11}-5547942043988\,{K}^{3}{a}^{6}{b}^{12}-\\
&961459535626\,{K}^{3}{a}
^{5}{b}^{13}-37855246951\,{K}^{3}{a}^{4}{b}^{14}+10318923385\,{K}^{3}{
a}^{3}{b}^{15}+958687833\,{K}^{3}{a}^{2}{b}^{16}-12573071\,{K}^{3}a{b}
^{17}-28708\,{K}^{3}{b}^{18}-\\
&3753273984\,{K}^{2}{a}^{18}b-14924756508
\,{K}^{2}{a}^{17}{b}^{2}+33759063064\,{K}^{2}{a}^{16}{b}^{3}+
379291052214\,{K}^{2}{a}^{15}{b}^{4}+1309224262221\,{K}^{2}{a}^{14}{b}
^{5}+\\
&2681708580065\,{K}^{2}{a}^{13}{b}^{6}+3558587977766\,{K}^{2}{a}^{
12}{b}^{7}+3063295002904\,{K}^{2}{a}^{11}{b}^{8}+1671009250448\,{K}^{2
}{a}^{10}{b}^{9}+529081273292\,{K}^{2}{a}^{9}{b}^{10}\\
&+62361748238\,{K}
^{2}{a}^{8}{b}^{11}-16165890171\,{K}^{2}{a}^{7}{b}^{12}-6850187226\,{K
}^{2}{a}^{6}{b}^{13}-778650662\,{K}^{2}{a}^{5}{b}^{14}+4578119\,{K}^{2
}{a}^{4}{b}^{15}+\\
&4071297\,{K}^{2}{a}^{3}{b}^{16}+21531\,{K}^{2}{a}^{2}
{b}^{17}+9138749030272\,{K}^{14}{a}^{6}-33138272229576\,{K}^{14}{a}^{5
}b-29076320448536\,{K}^{14}{a}^{4}{b}^{2}+\\
&21168171513346\,{K}^{14}{a}^
{3}{b}^{3}+9894961992239\,{K}^{14}{a}^{2}{b}^{4}+171744247857\,{K}^{14
}a{b}^{5}-102455992982\,{K}^{14}{b}^{6}+20482494237536\,{K}^{13}{a}^{7
}-\\
&758400710402384\,{K}^{13}{a}^{6}b+178209203479960\,{K}^{13}{a}^{5}{b
}^{2}+1061940602987627\,{K}^{13}{a}^{4}{b}^{3}+210465628830400\,{K}^{
13}{a}^{3}{b}^{4}-\\
&111438133952979\,{K}^{13}{a}^{2}{b}^{5}-
24835811525976\,{K}^{13}a{b}^{6}-1121798744793\,{K}^{13}{b}^{7}-
324686197454688\,{K}^{12}{a}^{8}-\\
&2397451956546320\,{K}^{12}{a}^{7}b+
9940072593584660\,{K}^{12}{a}^{6}{b}^{2}+9470849289821487\,{K}^{12}{a}
^{5}{b}^{3}-3259678991589712\,{K}^{12}{a}^{4}{b}^{4}-\\
&3712925575521778
\,{K}^{12}{a}^{3}{b}^{5}-498530006710001\,{K}^{12}{a}^{2}{b}^{6}+
61073494092387\,{K}^{12}a{b}^{7}+9140797073436\,{K}^{12}{b}^{8}+\\
&
421682681673960\,{K}^{11}{a}^{9}+7734644587643243\,{K}^{11}{a}^{8}b+
42848094216083152\,{K}^{11}{a}^{7}{b}^{2}-15841418773246606\,{K}^{11}{
a}^{6}{b}^{3}-\\
&68295751688444307\,{K}^{11}{a}^{5}{b}^{4}-
28832646138959683\,{K}^{11}{a}^{4}{b}^{5}+2931303589594457\,{K}^{11}{a
}^{3}{b}^{6}+3094786531802573\,{K}^{11}{a}^{2}{b}^{7}+\\
&383396493488585
\,{K}^{11}a{b}^{8}+4400321111259\,{K}^{11}{b}^{9}-333972747374480\,{K}
^{10}{a}^{10}-6092257153312040\,{K}^{10}{a}^{9}b-\\
&55856999800713909\,{K
}^{10}{a}^{8}{b}^{2}-264439656875882683\,{K}^{10}{a}^{7}{b}^{3}-
204449072809237951\,{K}^{10}{a}^{6}{b}^{4}+52200612796383814\,{K}^{10}
{a}^{5}{b}^{5}+\\
&104276314687547848\,{K}^{10}{a}^{4}{b}^{6}+
34516154437621492\,{K}^{10}{a}^{3}{b}^{7}+2323341037153351\,{K}^{10}{a
}^{2}{b}^{8}-499813548520597\,{K}^{10}a{b}^{9}-\\
&53742443038522\,{K}^{10
}{b}^{10}-719394998092120\,{K}^{9}{a}^{11}-3048692754220813\,{K}^{9}{a
}^{10}b+6597647439711413\,{K}^{9}{a}^{9}{b}^{2}+\\
&122396932782443872\,{K
}^{9}{a}^{8}{b}^{3}+651518163746829724\,{K}^{9}{a}^{7}{b}^{4}+
848919988177360099\,{K}^{9}{a}^{6}{b}^{5}+409310727291880508\,{K}^{9}{
a}^{5}{b}^{6}+\\
&28937156823814016\,{K}^{9}{a}^{4}{b}^{7}-
39015154046506957\,{K}^{9}{a}^{3}{b}^{8}-12714465165055035\,{K}^{9}{a}
^{2}{b}^{9}-1321491583878159\,{K}^{9}a{b}^{10}-\\
&29271868054181\,{K}^{9}
{b}^{11}+35750947789208\,{K}^{8}{a}^{12}+5739929124336014\,{K}^{8}{a}^
{11}b+35134843492212069\,{K}^{8}{a}^{10}{b}^{2}+\\
&94854587151399354\,{K}
^{8}{a}^{9}{b}^{3}+81067836562201188\,{K}^{8}{a}^{8}{b}^{4}-
481338847321338246\,{K}^{8}{a}^{7}{b}^{5}-1004093369821590900\,{K}^{8}
{a}^{6}{b}^{6}-\\
&799678186161717180\,{K}^{8}{a}^{5}{b}^{7}-
312977024146832130\,{K}^{8}{a}^{4}{b}^{8}-56983389925934480\,{K}^{8}{a
}^{3}{b}^{9}-1918129547961933\,{K}^{8}{a}^{2}{b}^{10}+\\
&733875799585003
\,{K}^{8}a{b}^{11}+69664042548080\,{K}^{8}{b}^{12}+235861764249280\,{K
}^{7}{a}^{13}+2662845304497119\,{K}^{7}{a}^{12}b+\\
&993218050827540\,{K}^
{7}{a}^{11}{b}^{2}-40759368682247813\,{K}^{7}{a}^{10}{b}^{3}-
175450024743402753\,{K}^{7}{a}^{9}{b}^{4}-374375009674292010\,{K}^{7}{
a}^{8}{b}^{5}-\\
&217515244579153062\,{K}^{7}{a}^{7}{b}^{6}+
202372455799240932\,{K}^{7}{a}^{6}{b}^{7}+358776169740042882\,{K}^{7}{
a}^{5}{b}^{8}+219861344810594501\,{K}^{7}{a}^{4}{b}^{9}+
\end{aligned}
\end{scriptsize}
$$
$$
\begin{scriptsize}
\begin{aligned}
\qquad & 70393640905527632\,{K}^{7}{a}^{3}{b}^{10}+12221971076415763\,{K}^{7}{a
}^{2}{b}^{11}+1044907079850673\,{K}^{7}a{b}^{12}+31291508547622\,{K}^{
7}{b}^{13}+\\
&65252427280480\,{K}^{6}{a}^{14}-150330308601000\,{K}^{6}{a}
^{13}b-5785751646434702\,{K}^{6}{a}^{12}{b}^{2}-23299799296986533\,{K}
^{6}{a}^{11}{b}^{3}-\\
&38082035644723335\,{K}^{6}{a}^{10}{b}^{4}+
1320023503009350\,{K}^{6}{a}^{9}{b}^{5}+144100540166508354\,{K}^{6}{a}
^{8}{b}^{6}+227482287519303780\,{K}^{6}{a}^{7}{b}^{7}+\\
&
142157384047334856\,{K}^{6}{a}^{6}{b}^{8}+22342529081227260\,{K}^{6}{a
}^{5}{b}^{9}-20131953878589013\,{K}^{6}{a}^{4}{b}^{10}-
13156604209825053\,{K}^{6}{a}^{3}{b}^{11}-\\
&3439758938852907\,{K}^{6}{a}
^{2}{b}^{12}-433990265452981\,{K}^{6}a{b}^{13}-21625409989612\,{K}^{6}
{b}^{14}-1305982952248\,{K}^{5}{a}^{15}-\\
&224498362551255\,{K}^{5}{a}^{
14}b-1380789812767895\,{K}^{5}{a}^{13}{b}^{2}-1769928094065453\,{K}^{5
}{a}^{12}{b}^{3}+5554277135328160\,{K}^{5}{a}^{11}{b}^{4}+\\
&
26283564358823838\,{K}^{5}{a}^{10}{b}^{5}+48887380527984955\,{K}^{5}{a
}^{9}{b}^{6}+37983269049561989\,{K}^{5}{a}^{8}{b}^{7}+1692007648022372
\,{K}^{5}{a}^{7}{b}^{8}-\\
&17713282910012387\,{K}^{5}{a}^{6}{b}^{9}-
12868753933347934\,{K}^{5}{a}^{5}{b}^{10}-3887837832743402\,{K}^{5}{a}
^{4}{b}^{11}-332321213422364\,{K}^{5}{a}^{3}{b}^{12}+\\
&102067306573954\,
{K}^{5}{a}^{2}{b}^{13}+27720215610057\,{K}^{5}a{b}^{14}+1987157066858
\,{K}^{5}{b}^{15}-2805917499368\,{K}^{4}{a}^{16}-37545518125666\,{K}^{
4}{a}^{15}b-\\
&70258321835843\,{K}^{4}{a}^{14}{b}^{2}+413930136572183\,{K
}^{4}{a}^{13}{b}^{3}+2151366507275804\,{K}^{4}{a}^{12}{b}^{4}+
4636288677395658\,{K}^{4}{a}^{11}{b}^{5}+\\
&5092106524646667\,{K}^{4}{a}^
{10}{b}^{6}+998220016538121\,{K}^{4}{a}^{9}{b}^{7}-3646615041089044\,{
K}^{4}{a}^{8}{b}^{8}-4120895921425718\,{K}^{4}{a}^{7}{b}^{9}-\\
&
1914545156743972\,{K}^{4}{a}^{6}{b}^{10}-324737907212762\,{K}^{4}{a}^{
5}{b}^{11}+63031628775294\,{K}^{4}{a}^{4}{b}^{12}+36826384244964\,{K}^
{4}{a}^{3}{b}^{13}+\\
&5293294253606\,{K}^{4}{a}^{2}{b}^{14}+85768464196\,
{K}^{4}a{b}^{15}-25488417092\,{K}^{4}{b}^{16}-409930836952\,{K}^{3}{a}
^{17}-2133533635942\,{K}^{3}{a}^{16}b+\\
&7428746593498\,{K}^{3}{a}^{15}{b
}^{2}+73670450261685\,{K}^{3}{a}^{14}{b}^{3}+218924804314412\,{K}^{3}{
a}^{13}{b}^{4}+343144143406698\,{K}^{3}{a}^{12}{b}^{5}+\\
&237004584728004
\,{K}^{3}{a}^{11}{b}^{6}-143240215520573\,{K}^{3}{a}^{10}{b}^{7}-
426084388561338\,{K}^{3}{a}^{9}{b}^{8}-360576888077064\,{K}^{3}{a}^{8}
{b}^{9}-\\
&145835956270225\,{K}^{3}{a}^{7}{b}^{10}-19426909116396\,{K}^{3
}{a}^{6}{b}^{11}+6710227592598\,{K}^{3}{a}^{5}{b}^{12}+3117904909227\,
{K}^{3}{a}^{4}{b}^{13}+\\
&425656971047\,{K}^{3}{a}^{3}{b}^{14}+7049522006
\,{K}^{3}{a}^{2}{b}^{15}-1836356608\,{K}^{3}a{b}^{16}+5760015\,{K}^{3}
{b}^{17}-18814990944\,{K}^{2}{a}^{18}-3489515052\,{K}^{2}{a}^{17}b+\\
&
827688535441\,{K}^{2}{a}^{16}{b}^{2}+4215525579574\,{K}^{2}{a}^{15}{b}
^{3}+10457757597703\,{K}^{2}{a}^{14}{b}^{4}+15262168445017\,{K}^{2}{a}
^{13}{b}^{5}+\\
&10803295881005\,{K}^{2}{a}^{12}{b}^{6}-3610463346757\,{K}
^{2}{a}^{11}{b}^{7}-14788999162838\,{K}^{2}{a}^{10}{b}^{8}-
13817876558200\,{K}^{2}{a}^{9}{b}^{9}-\\
&6717126825642\,{K}^{2}{a}^{8}{b}
^{10}-1731676330004\,{K}^{2}{a}^{7}{b}^{11}-153519236210\,{K}^{2}{a}^{
6}{b}^{12}+29717307962\,{K}^{2}{a}^{5}{b}^{13}+7760448081\,{K}^{2}{a}^
{4}{b}^{14}+\\
&369897981\,{K}^{2}{a}^{3}{b}^{15}-16556109\,{K}^{2}{a}^{2}
{b}^{16}-64593\,{K}^{2}a{b}^{17}+2045351484\,K{a}^{18}b+17811440660\,K
{a}^{17}{b}^{2}+75580840020\,K{a}^{16}{b}^{3}+\\
&209265992198\,K{a}^{15}{
b}^{4}+399232201977\,K{a}^{14}{b}^{5}+509916554269\,K{a}^{13}{b}^{6}+
418273114533\,K{a}^{12}{b}^{7}+208653225840\,K{a}^{11}{b}^{8}+\\
&
52492565476\,K{a}^{10}{b}^{9}-1142791751\,K{a}^{9}{b}^{10}-4579062339
\,K{a}^{8}{b}^{11}-1083565756\,K{a}^{7}{b}^{12}-48764744\,K{a}^{6}{b}^
{13}+12246127\,K{a}^{5}{b}^{14}+\\
&1053473\,K{a}^{4}{b}^{15}+7177\,K{a}^{
3}{b}^{16}-8366748854400\,{K}^{14}{a}^{5}-2116251005144\,{K}^{14}{a}^{
4}b+16210780951232\,{K}^{14}{a}^{3}{b}^{2}+\\
&4480117880798\,{K}^{14}{a}^
{2}{b}^{3}-630695848917\,{K}^{14}a{b}^{4}-152455941434\,{K}^{14}{b}^{5
}-146541897225952\,{K}^{13}{a}^{6}+309118006658912\,{K}^{13}{a}^{5}b+\\
&
545337679295772\,{K}^{13}{a}^{4}{b}^{2}-52349595300926\,{K}^{13}{a}^{3
}{b}^{3}-136256003317228\,{K}^{13}{a}^{2}{b}^{4}-20209597010776\,{K}^{
13}a{b}^{5}-\\
&387111808473\,{K}^{13}{b}^{6}-315742065011680\,{K}^{12}{a}
^{7}+4940793788557920\,{K}^{12}{a}^{6}b+2175422434441109\,{K}^{12}{a}^
{5}{b}^{2}-\\
&5315716943636116\,{K}^{12}{a}^{4}{b}^{3}-3001415087297317\,
{K}^{12}{a}^{3}{b}^{4}-19117863087767\,{K}^{12}{a}^{2}{b}^{5}+
140293899704984\,{K}^{12}a{b}^{6}+\\
&13304021281188\,{K}^{12}{b}^{7}+
1421474467075936\,{K}^{11}{a}^{8}+14454455221835888\,{K}^{11}{a}^{7}b-
29580391851589347\,{K}^{11}{a}^{6}{b}^{2}-\\
&52522280994575503\,{K}^{11}{
a}^{5}{b}^{3}-9204277797212104\,{K}^{11}{a}^{4}{b}^{4}+
10875640337748816\,{K}^{11}{a}^{3}{b}^{5}+4115195257179398\,{K}^{11}{a
}^{2}{b}^{6}+\\
&310903254459242\,{K}^{11}a{b}^{7}-12554723530946\,{K}^{11
}{b}^{8}-1177828168813092\,{K}^{10}{a}^{9}-22011654386909474\,{K}^{10}
{a}^{8}b-\\
&154112460984303377\,{K}^{10}{a}^{7}{b}^{2}-59161960595274706
\,{K}^{10}{a}^{6}{b}^{3}+147595786252235791\,{K}^{10}{a}^{5}{b}^{4}+
125655350928421034\,{K}^{10}{a}^{4}{b}^{5}+\\
&26498483100469881\,{K}^{10}
{a}^{3}{b}^{6}-2715367502572491\,{K}^{10}{a}^{2}{b}^{7}-
1322965401614618\,{K}^{10}a{b}^{8}-91520154306013\,{K}^{10}{b}^{9}-\\
&
38554305401828\,{K}^{9}{a}^{10}+6600979950639959\,{K}^{9}{a}^{9}b+
89989529127614764\,{K}^{9}{a}^{8}{b}^{2}+593849653565552738\,{K}^{9}{a
}^{7}{b}^{3}+\\
&724094342011802389\,{K}^{9}{a}^{6}{b}^{4}+
218658589577952796\,{K}^{9}{a}^{5}{b}^{5}-110451344062183908\,{K}^{9}{
a}^{4}{b}^{6}-86724004663199209\,{K}^{9}{a}^{3}{b}^{7}-\\
&
19302337084572353\,{K}^{9}{a}^{2}{b}^{8}-1339910186024989\,{K}^{9}a{b}
^{9}+17636007777647\,{K}^{9}{b}^{10}+1752494790353076\,{K}^{8}{a}^{11}
+\\
&14229086143906700\,{K}^{8}{a}^{10}b+44180642260988501\,{K}^{8}{a}^{9}
{b}^{2}+745343463110076\,{K}^{8}{a}^{8}{b}^{3}-776040806402193429\,{K}
^{8}{a}^{7}{b}^{4}-\\
&1437803549840667906\,{K}^{8}{a}^{6}{b}^{5}-
1049123594148428808\,{K}^{8}{a}^{5}{b}^{6}-343358674437483282\,{K}^{8}
{a}^{4}{b}^{7}-30317747851594845\,{K}^{8}{a}^{3}{b}^{8}+\\
&
9933980722781458\,{K}^{8}{a}^{2}{b}^{9}+2546563239464043\,{K}^{8}a{b}^
{10}+163338988699983\,{K}^{8}{b}^{11}+519083783447092\,{K}^{7}{a}^{12}
-\\
&4070943108490139\,{K}^{7}{a}^{11}b-46371736892722249\,{K}^{7}{a}^{10}
{b}^{2}-188078510911418113\,{K}^{7}{a}^{9}{b}^{3}-424488962108604480\,
{K}^{7}{a}^{8}{b}^{4}-\\
&52028578250464806\,{K}^{7}{a}^{7}{b}^{5}+
690280096834985520\,{K}^{7}{a}^{6}{b}^{6}+841825171719319074\,{K}^{7}{
a}^{5}{b}^{7}+460112829568030812\,{K}^{7}{a}^{4}{b}^{8}+\\
&
134237801030983929\,{K}^{7}{a}^{3}{b}^{9}+20541601882032779\,{K}^{7}{a
}^{2}{b}^{10}+1373223356467419\,{K}^{7}a{b}^{11}+15985249590442\,{K}^{
7}{b}^{12}-\\
&238703163488284\,{K}^{6}{a}^{13}-4704244264311646\,{K}^{6}{
a}^{12}b-17323750084546248\,{K}^{6}{a}^{11}{b}^{2}-15752294278585752\,
{K}^{6}{a}^{10}{b}^{3}+\\
&75356772277609232\,{K}^{6}{a}^{9}{b}^{4}+
353320178036431068\,{K}^{6}{a}^{8}{b}^{5}+458020431618463944\,{K}^{6}{
a}^{7}{b}^{6}+198566807698104318\,{K}^{6}{a}^{6}{b}^{7}-\\
&
64552243299115518\,{K}^{6}{a}^{5}{b}^{8}-106645675418017548\,{K}^{6}{a
}^{4}{b}^{9}-49000839732121651\,{K}^{6}{a}^{3}{b}^{10}-
11279031589336739\,{K}^{6}{a}^{2}{b}^{11}-\\
&1312623662618883\,{K}^{6}a{b
}^{12}-60994131496839\,{K}^{6}{b}^{13}-116193725316636\,{K}^{5}{a}^{14
}-721878196344595\,{K}^{5}{a}^{13}b+
\end{aligned}
\end{scriptsize}
$$
$$
\begin{scriptsize}
\begin{aligned}
\qquad &1610835370041160\,{K}^{5}{a}^{12}{
b}^{2}+17605130966084267\,{K}^{5}{a}^{11}{b}^{3}+54280149614660941\,{K
}^{5}{a}^{10}{b}^{4}+86800845930588638\,{K}^{5}{a}^{9}{b}^{5}+\\
&
32287612722246073\,{K}^{5}{a}^{8}{b}^{6}-64739932153516130\,{K}^{5}{a}
^{7}{b}^{7}-87190925399159556\,{K}^{5}{a}^{6}{b}^{8}-43913987835106271
\,{K}^{5}{a}^{5}{b}^{9}-\\
&7974713564514131\,{K}^{5}{a}^{4}{b}^{10}+
1524852527121546\,{K}^{5}{a}^{3}{b}^{11}+991358268709255\,{K}^{5}{a}^{
2}{b}^{12}+174262222428333\,{K}^{5}a{b}^{13}+\\
&10821328328184\,{K}^{5}{b
}^{14}-11261563958948\,{K}^{4}{a}^{15}+89249212446876\,{K}^{4}{a}^{14}
b+1210837259888850\,{K}^{4}{a}^{13}{b}^{2}+\\
&4189391894888422\,{K}^{4}{a
}^{12}{b}^{3}+6895753063016786\,{K}^{4}{a}^{11}{b}^{4}+
2814580064017048\,{K}^{4}{a}^{10}{b}^{5}-12639458890823169\,{K}^{4}{a}
^{9}{b}^{6}-\\
&22966133310322834\,{K}^{4}{a}^{8}{b}^{7}-16018386748175986
\,{K}^{4}{a}^{7}{b}^{8}-3897474512992154\,{K}^{4}{a}^{6}{b}^{9}+
1270125106352083\,{K}^{4}{a}^{5}{b}^{10}+\\
&1096703512386140\,{K}^{4}{a}^
{4}{b}^{11}+272735729296492\,{K}^{4}{a}^{3}{b}^{12}+19499559507513\,{K
}^{4}{a}^{2}{b}^{13}-2525285299791\,{K}^{4}a{b}^{14}-\\
&355550893503\,{K}
^{4}{b}^{15}+1240525433820\,{K}^{3}{a}^{16}+30996444517335\,{K}^{3}{a}
^{15}b+161254059835984\,{K}^{3}{a}^{14}{b}^{2}+323480035254894\,{K}^{3
}{a}^{13}{b}^{3}\\
&+139340497312990\,{K}^{3}{a}^{12}{b}^{4}-
850006077901653\,{K}^{3}{a}^{11}{b}^{5}-2322419584135520\,{K}^{3}{a}^{
10}{b}^{6}-2561802646529055\,{K}^{3}{a}^{9}{b}^{7}-\\
&1248015994554538\,{
K}^{3}{a}^{8}{b}^{8}-7430614220767\,{K}^{3}{a}^{7}{b}^{9}+
295552588407043\,{K}^{3}{a}^{6}{b}^{10}+143471865551486\,{K}^{3}{a}^{5
}{b}^{11}+\\
&26870041106910\,{K}^{3}{a}^{4}{b}^{12}+357355312777\,{K}^{3}
{a}^{3}{b}^{13}-520758079751\,{K}^{3}{a}^{2}{b}^{14}-43581243653\,{K}^
{3}a{b}^{15}+932292774\,{K}^{3}{b}^{16}+\\
&281029101328\,{K}^{2}{a}^{17}+
2790875057008\,{K}^{2}{a}^{16}b+9004214636313\,{K}^{2}{a}^{15}{b}^{2}+
9241415030100\,{K}^{2}{a}^{14}{b}^{3}-17561364893706\,{K}^{2}{a}^{13}{
b}^{4}-\\
&84206866299774\,{K}^{2}{a}^{12}{b}^{5}-154505371146174\,{K}^{2}
{a}^{11}{b}^{6}-147155794533792\,{K}^{2}{a}^{10}{b}^{7}-67126980619367
\,{K}^{2}{a}^{9}{b}^{8}-\\
&683626744636\,{K}^{2}{a}^{8}{b}^{9}+
15158608019533\,{K}^{2}{a}^{7}{b}^{10}+7631234344314\,{K}^{2}{a}^{6}{b
}^{11}+1617373647145\,{K}^{2}{a}^{5}{b}^{12}+102616130082\,{K}^{2}{a}^
{4}{b}^{13}-\\
&13253544130\,{K}^{2}{a}^{3}{b}^{14}-1551176922\,{K}^{2}{a}
^{2}{b}^{15}+18396690\,{K}^{2}a{b}^{16}+43062\,{K}^{2}{b}^{17}+
13054956400\,K{a}^{18}+89834852400\,K{a}^{17}b+\\
&222805807756\,K{a}^{16}
{b}^{2}+162505467409\,K{a}^{15}{b}^{3}-551972439658\,K{a}^{14}{b}^{4}-
2285235255513\,K{a}^{13}{b}^{5}-4269981779179\,K{a}^{12}{b}^{6}-\\
&
4588444491892\,K{a}^{11}{b}^{7}-2957997174947\,K{a}^{10}{b}^{8}-
1097913077173\,K{a}^{9}{b}^{9}-177739124652\,K{a}^{8}{b}^{10}+
20064773978\,K{a}^{7}{b}^{11}+\\
&13229913738\,K{a}^{6}{b}^{12}+1728840840
\,K{a}^{5}{b}^{13}+5700582\,K{a}^{4}{b}^{14}-8034939\,K{a}^{3}{b}^{15}
-43062\,K{a}^{2}{b}^{16}+681783828\,{a}^{17}{b}^{2}+\\
&5012412148\,{a}^{
16}{b}^{3}+14948674864\,{a}^{15}{b}^{4}+23306997007\,{a}^{14}{b}^{5}+
20448874409\,{a}^{13}{b}^{6}+9886329587\,{a}^{12}{b}^{7}+1912559211\,{
a}^{11}{b}^{8}-\\
&438777291\,{a}^{10}{b}^{9}-303787639\,{a}^{9}{b}^{10}-
46063842\,{a}^{8}{b}^{11}+2626637\,{a}^{7}{b}^{12}+970084\,{a}^{6}{b}^
{13}+7177\,{a}^{5}{b}^{14}+1614824082688\,{K}^{14}{a}^{4}+\\
&
5519743940372\,{K}^{14}{a}^{3}b+230786784036\,{K}^{14}{a}^{2}{b}^{2}-
791555241926\,{K}^{14}a{b}^{3}-124160225003\,{K}^{14}{b}^{4}+
97929043467680\,{K}^{13}{a}^{5}+\\
&108764092290156\,{K}^{13}{a}^{4}b-
131462140629948\,{K}^{13}{a}^{3}{b}^{2}-88542785125920\,{K}^{13}{a}^{2
}{b}^{3}-6716374959480\,{K}^{13}a{b}^{4}+643233348319\,{K}^{13}{b}^{5}
\\
&+1027078408205536\,{K}^{12}{a}^{6}-859583263625430\,{K}^{12}{a}^{5}b-
3529611326874619\,{K}^{12}{a}^{4}{b}^{2}-1086186158789088\,{K}^{12}{a}
^{3}{b}^{3}+\\
&442429093252751\,{K}^{12}{a}^{2}{b}^{4}+164844088446624\,{
K}^{12}a{b}^{5}+11951095879494\,{K}^{12}{b}^{6}+2150886559698592\,{K}^
{11}{a}^{7}-\\
&17890872284315199\,{K}^{11}{a}^{6}b-20246399144055212\,{K}
^{11}{a}^{5}{b}^{2}+8635005176568966\,{K}^{11}{a}^{4}{b}^{3}+
12892051050045372\,{K}^{11}{a}^{3}{b}^{4}+\\
&3012342425542590\,{K}^{11}{a
}^{2}{b}^{5}-23116584010883\,{K}^{11}a{b}^{6}-38623682161761\,{K}^{11}
{b}^{7}-3962218577885440\,{K}^{10}{a}^{8}-\\
&54214894714616568\,{K}^{10}{
a}^{7}b+34602193697250040\,{K}^{10}{a}^{6}{b}^{2}+149322857531949587\,
{K}^{10}{a}^{5}{b}^{3}+83460921967092918\,{K}^{10}{a}^{4}{b}^{4}+\\
&
396855351645148\,{K}^{10}{a}^{3}{b}^{5}-9107363446484598\,{K}^{10}{a}^
{2}{b}^{6}-1887045030143888\,{K}^{10}a{b}^{7}-93520966182042\,{K}^{10}
{b}^{8}+\\
&1627284760795384\,{K}^{9}{a}^{9}+37033158237054487\,{K}^{9}{a}
^{8}b+361998131827416651\,{K}^{9}{a}^{7}{b}^{2}+362985502015361910\,{K
}^{9}{a}^{6}{b}^{3}-\\
&71892060366817546\,{K}^{9}{a}^{5}{b}^{4}-
232679552490502382\,{K}^{9}{a}^{4}{b}^{5}-105410113570868367\,{K}^{9}{
a}^{3}{b}^{6}-16016766305861726\,{K}^{9}{a}^{2}{b}^{7}+\\
&41929311676730
\,{K}^{9}a{b}^{8}+134786061450036\,{K}^{9}{b}^{9}+1850077795088840\,{K
}^{8}{a}^{10}+8443323766110210\,{K}^{8}{a}^{9}b-\\
&37947272762636197\,{K}
^{8}{a}^{8}{b}^{2}-810178900187248811\,{K}^{8}{a}^{7}{b}^{3}-
1396178310207007158\,{K}^{8}{a}^{6}{b}^{4}-880295099272673493\,{K}^{8}
{a}^{5}{b}^{5}-\\
&167162832065043215\,{K}^{8}{a}^{4}{b}^{6}+
53988668294091505\,{K}^{8}{a}^{3}{b}^{7}+29836603806311020\,{K}^{8}{a}
^{2}{b}^{8}+4677855084070136\,{K}^{8}a{b}^{9}+\\
&235307313839560\,{K}^{8}
{b}^{10}-2100596008043760\,{K}^{7}{a}^{11}-23631537952487111\,{K}^{7}{
a}^{10}b-113786586006435010\,{K}^{7}{a}^{9}{b}^{2}-\\
&304040055521706975
\,{K}^{7}{a}^{8}{b}^{3}+306140767666768896\,{K}^{7}{a}^{7}{b}^{4}+
1289942679419409690\,{K}^{7}{a}^{6}{b}^{5}+1324303533949334025\,{K}^{7
}{a}^{5}{b}^{6}+\\
&652415011215183357\,{K}^{7}{a}^{4}{b}^{7}+
165717791932758393\,{K}^{7}{a}^{3}{b}^{8}+18938588948132107\,{K}^{7}{a
}^{2}{b}^{9}+227884632639923\,{K}^{7}a{b}^{10}-\\
&84502656946761\,{K}^{7}
{b}^{11}-1257588617880416\,{K}^{6}{a}^{12}-3847594481754996\,{K}^{6}{a
}^{11}b+14426948700417274\,{K}^{6}{a}^{10}{b}^{2}+\\
&134545397298508222\,
{K}^{6}{a}^{9}{b}^{3}+515153469066747609\,{K}^{6}{a}^{8}{b}^{4}+
557509213844545830\,{K}^{6}{a}^{7}{b}^{5}+30943238638887114\,{K}^{6}{a
}^{6}{b}^{6}-\\
&360588895263229974\,{K}^{6}{a}^{5}{b}^{7}-
304333571032361316\,{K}^{6}{a}^{4}{b}^{8}-118995812023316076\,{K}^{6}{
a}^{3}{b}^{9}-24929574953717489\,{K}^{6}{a}^{2}{b}^{10}-\\
&
2666389536683722\,{K}^{6}a{b}^{11}-112222341492684\,{K}^{6}{b}^{12}+
2538297716520\,{K}^{5}{a}^{13}+3347416391714541\,{K}^{5}{a}^{12}b+\\
&
19920071276972533\,{K}^{5}{a}^{11}{b}^{2}+55313252285113484\,{K}^{5}{a
}^{10}{b}^{3}+75555867405247854\,{K}^{5}{a}^{9}{b}^{4}-
67159472042812458\,{K}^{5}{a}^{8}{b}^{5}-\\
&243122267317670804\,{K}^{5}{a
}^{7}{b}^{6}-218033484859290454\,{K}^{5}{a}^{6}{b}^{7}-
75280843775858384\,{K}^{5}{a}^{5}{b}^{8}+6200818824054073\,{K}^{5}{a}^
{4}{b}^{9}+\\
&13782248515250288\,{K}^{5}{a}^{3}{b}^{10}+4627826510035521
\,{K}^{5}{a}^{2}{b}^{11}+683488148077625\,{K}^{5}a{b}^{12}+
38908464172393\,{K}^{5}{b}^{13}+\\
&92622876136344\,{K}^{4}{a}^{14}+
1022220688420964\,{K}^{4}{a}^{13}b+2797108808357214\,{K}^{4}{a}^{12}{b
}^{2}+1185678317068771\,{K}^{4}{a}^{11}{b}^{3}-\\
&13686874087704332\,{K}^
{4}{a}^{10}{b}^{4}-47950209446496673\,{K}^{4}{a}^{9}{b}^{5}-
57324944188297218\,{K}^{4}{a}^{8}{b}^{6}-23266367866334726\,{K}^{4}{a}
^{7}{b}^{7}+\\
&8985510483364366\,{K}^{4}{a}^{6}{b}^{8}+13077781424268410
\,{K}^{4}{a}^{5}{b}^{9}+5280145125939507\,{K}^{4}{a}^{4}{b}^{10}+
784959157521536\,{K}^{4}{a}^{3}{b}^{11}-
\end{aligned}
\end{scriptsize}
$$
$$
\begin{scriptsize}
\begin{aligned}
\qquad &58417423386961\,{K}^{4}{a}^{2}
{b}^{12}-31259161729816\,{K}^{4}a{b}^{13}-2655708698282\,{K}^{4}{b}^{
14}+15249366629856\,{K}^{3}{a}^{15}+67554929337906\,{K}^{3}{a}^{14}b-\\
&
137833709047054\,{K}^{3}{a}^{13}{b}^{2}-1292770433853357\,{K}^{3}{a}^{
12}{b}^{3}-3851692745460516\,{K}^{3}{a}^{11}{b}^{4}-6410618710506657\,
{K}^{3}{a}^{10}{b}^{5}-\\
&4299372258492045\,{K}^{3}{a}^{9}{b}^{6}+
1063928497353022\,{K}^{3}{a}^{8}{b}^{7}+3494921421889280\,{K}^{3}{a}^{
7}{b}^{8}+2214951782966338\,{K}^{3}{a}^{6}{b}^{9}+\\
&567633034371841\,{K}
^{3}{a}^{5}{b}^{10}-4584814104768\,{K}^{3}{a}^{4}{b}^{11}-
34354481224039\,{K}^{3}{a}^{3}{b}^{12}-6321327386477\,{K}^{3}{a}^{2}{b
}^{13}-\\
&188568579196\,{K}^{3}a{b}^{14}+27858856338\,{K}^{3}{b}^{15}+
410954761712\,{K}^{2}{a}^{16}-5216455629732\,{K}^{2}{a}^{15}b-
49197073597600\,{K}^{2}{a}^{14}{b}^{2}-\\
&161864662576359\,{K}^{2}{a}^{13
}{b}^{3}-317036422019972\,{K}^{2}{a}^{12}{b}^{4}-364953496174642\,{K}^
{2}{a}^{11}{b}^{5}-99411705279193\,{K}^{2}{a}^{10}{b}^{6}+\\
&
250051260866918\,{K}^{2}{a}^{9}{b}^{7}+318439285188202\,{K}^{2}{a}^{8}
{b}^{8}+164537796710118\,{K}^{2}{a}^{7}{b}^{9}+35050186179610\,{K}^{2}
{a}^{6}{b}^{10}-\\
&2703566292537\,{K}^{2}{a}^{5}{b}^{11}-2831573479843\,{
K}^{2}{a}^{4}{b}^{12}-466591953622\,{K}^{2}{a}^{3}{b}^{13}-12160965805
\,{K}^{2}{a}^{2}{b}^{14}+1907145344\,{K}^{2}a{b}^{15}-\\
&5587767\,{K}^{2}
{b}^{16}-62254358224\,K{a}^{17}-770472809812\,K{a}^{16}b-3379392782476
\,K{a}^{15}{b}^{2}-8147876651744\,K{a}^{14}{b}^{3}-\\
&13236693819432\,K{a
}^{13}{b}^{4}-13361840370126\,K{a}^{12}{b}^{5}-3519734728127\,K{a}^{11
}{b}^{6}+8183717267411\,K{a}^{10}{b}^{7}+10822632016072\,K{a}^{9}{b}^{
8}+\\
&6196545977028\,K{a}^{8}{b}^{9}+1829412618454\,K{a}^{7}{b}^{10}+
211896799424\,K{a}^{6}{b}^{11}-21837978980\,K{a}^{5}{b}^{12}-
7776252183\,K{a}^{4}{b}^{13}-\\
&414798211\,K{a}^{3}{b}^{14}+16297737\,K{a
}^{2}{b}^{15}+64593\,Ka{b}^{16}-2727135312\,{a}^{18}-22062296364\,{a}^
{17}b-79392213736\,{a}^{16}{b}^{2}-\\
&184374576340\,{a}^{15}{b}^{3}-
324668442393\,{a}^{14}{b}^{4}-423133648410\,{a}^{13}{b}^{5}-
369907482509\,{a}^{12}{b}^{6}-200291701128\,{a}^{11}{b}^{7}-\\
&
57029556216\,{a}^{10}{b}^{8}-1499096959\,{a}^{9}{b}^{9}+4189548983\,{a
}^{8}{b}^{10}+1117527129\,{a}^{7}{b}^{11}+59541131\,{a}^{6}{b}^{12}-
12160003\,{a}^{5}{b}^{13}-\\
&1053473\,{a}^{4}{b}^{14}-7177\,{a}^{3}{b}^{
15}+523166929152\,{K}^{14}{a}^{3}-582263613540\,{K}^{14}{a}^{2}b-
393516160888\,{K}^{14}a{b}^{2}-54555931314\,{K}^{14}{b}^{3}-\\
&
6703824234896\,{K}^{13}{a}^{4}-61702716576336\,{K}^{13}{a}^{3}b-
25457094964996\,{K}^{13}{a}^{2}{b}^{2}+2968926865726\,{K}^{13}a{b}^{3}
+1005450782159\,{K}^{13}{b}^{4}-\\
&476267897041792\,{K}^{12}{a}^{5}-
972569576047662\,{K}^{12}{a}^{4}b+165686662778546\,{K}^{12}{a}^{3}{b}^
{2}+466120347138922\,{K}^{12}{a}^{2}{b}^{3}+\\
&104809921907461\,{K}^{12}a
{b}^{4}+5177081618228\,{K}^{12}{b}^{5}-4138622857713832\,{K}^{11}{a}^{
6}-1125781979538292\,{K}^{11}{a}^{5}b+\\
&10576906994896955\,{K}^{11}{a}^{
4}{b}^{2}+7452560812121124\,{K}^{11}{a}^{3}{b}^{3}+626312915027250\,{K
}^{11}{a}^{2}{b}^{4}-354644102279253\,{K}^{11}a{b}^{5}-\\
&51942755632192
\,{K}^{11}{b}^{6}-8663364237320236\,{K}^{10}{a}^{7}+37769329487049634
\,{K}^{10}{a}^{6}b+74558312086146548\,{K}^{10}{a}^{5}{b}^{2}+\\
&
17676360783868193\,{K}^{10}{a}^{4}{b}^{3}-20600049739887781\,{K}^{10}{
a}^{3}{b}^{4}-10729424364136168\,{K}^{10}{a}^{2}{b}^{5}-
1482908722571911\,{K}^{10}a{b}^{6}-\\
&33679930807438\,{K}^{10}{b}^{7}+
6738283536473268\,{K}^{9}{a}^{8}+133758323608241233\,{K}^{9}{a}^{7}b+
53888339905194366\,{K}^{9}{a}^{6}{b}^{2}-\\
&214295598624968657\,{K}^{9}{a
}^{5}{b}^{3}-219327009013586683\,{K}^{9}{a}^{4}{b}^{4}-
67355940946698385\,{K}^{9}{a}^{3}{b}^{5}-1514213957247314\,{K}^{9}{a}^
{2}{b}^{6}+\\
&2276615089100975\,{K}^{9}a{b}^{7}+254170715002320\,{K}^{9}{
b}^{8}+203767459672492\,{K}^{8}{a}^{9}-23940454794729430\,{K}^{8}{a}^{
8}b-\\
&545375697461730317\,{K}^{8}{a}^{7}{b}^{2}-857365742427372677\,{K}^
{8}{a}^{6}{b}^{3}-360180976924310186\,{K}^{8}{a}^{5}{b}^{4}+
115513026118930214\,{K}^{8}{a}^{4}{b}^{5}+\\
&138921354018382773\,{K}^{8}{
a}^{3}{b}^{6}+41743818285078668\,{K}^{8}{a}^{2}{b}^{7}+
4926923941645232\,{K}^{8}a{b}^{8}+169169972401763\,{K}^{8}{b}^{9}-\\
&
3991752731304356\,{K}^{7}{a}^{10}-34689578338228949\,{K}^{7}{a}^{9}b-
132928355759287251\,{K}^{7}{a}^{8}{b}^{2}+561781578633979853\,{K}^{7}{
a}^{7}{b}^{3}+\\
&1511646200314460345\,{K}^{7}{a}^{6}{b}^{4}+
1374432963207443419\,{K}^{7}{a}^{5}{b}^{5}+582998342501689059\,{K}^{7}
{a}^{4}{b}^{6}+104726594089449981\,{K}^{7}{a}^{3}{b}^{7}-\\
&
1501885683355848\,{K}^{7}{a}^{2}{b}^{8}-2917097228636804\,{K}^{7}a{b}^
{9}-267417025463264\,{K}^{7}{b}^{10}+832686997621220\,{K}^{6}{a}^{11}+\\
&
15818232972782436\,{K}^{6}{a}^{10}b+106919671217024530\,{K}^{6}{a}^{9}
{b}^{2}+467279862910883277\,{K}^{6}{a}^{8}{b}^{3}+344995611141785220\,
{K}^{6}{a}^{7}{b}^{4}-\\
&422954178831375360\,{K}^{6}{a}^{6}{b}^{5}-
815865769095384684\,{K}^{6}{a}^{5}{b}^{6}-549163646778006354\,{K}^{6}{
a}^{4}{b}^{7}-191998001382464763\,{K}^{6}{a}^{3}{b}^{8}-\\
&
36317745745231719\,{K}^{6}{a}^{2}{b}^{9}-3397378173585692\,{K}^{6}a{b}
^{10}-113763805686559\,{K}^{6}{b}^{11}+1269254437776660\,{K}^{5}{a}^{
12}+\\
&8437628849153881\,{K}^{5}{a}^{11}b+23512183945971964\,{K}^{5}{a}^{
10}{b}^{2}+14110888588733283\,{K}^{5}{a}^{9}{b}^{3}-226935596423354970
\,{K}^{5}{a}^{8}{b}^{4}-\\
&445461735393919422\,{K}^{5}{a}^{7}{b}^{5}-
297350807047405130\,{K}^{5}{a}^{6}{b}^{6}-20098198984754930\,{K}^{5}{a
}^{5}{b}^{7}+80727940695012664\,{K}^{5}{a}^{4}{b}^{8}+\\
&
49193640052921219\,{K}^{5}{a}^{3}{b}^{9}+13312711315555667\,{K}^{5}{a}
^{2}{b}^{10}+1775566478252544\,{K}^{5}a{b}^{11}+94495364210863\,{K}^{5
}{b}^{12}+\\
&194950664145404\,{K}^{4}{a}^{13}-301960685627314\,{K}^{4}{a}
^{12}b-6340477509499017\,{K}^{4}{a}^{11}{b}^{2}-28598413457245203\,{K}
^{4}{a}^{10}{b}^{3}-\\
&75794221239015260\,{K}^{4}{a}^{9}{b}^{4}-
63261918228686268\,{K}^{4}{a}^{8}{b}^{5}+19729125651706110\,{K}^{4}{a}
^{7}{b}^{6}+65363113894261356\,{K}^{4}{a}^{6}{b}^{7}+\\
&43249734775404058
\,{K}^{4}{a}^{5}{b}^{8}+11400699018541118\,{K}^{4}{a}^{4}{b}^{9}-
199122852450284\,{K}^{4}{a}^{3}{b}^{10}-817175792486743\,{K}^{4}{a}^{2
}{b}^{11}-\\
&175510220603759\,{K}^{4}a{b}^{12}-12253125734325\,{K}^{4}{b}
^{13}-24298305089244\,{K}^{3}{a}^{14}-443458042479681\,{K}^{3}{a}^{13}
b-\\
&1876636236634742\,{K}^{3}{a}^{12}{b}^{2}-4415698696303203\,{K}^{3}{a
}^{11}{b}^{3}-5556422479047489\,{K}^{3}{a}^{10}{b}^{4}+
3450792674167897\,{K}^{3}{a}^{9}{b}^{5}+\\
&14638478231500723\,{K}^{3}{a}^
{8}{b}^{6}+13624256031115082\,{K}^{3}{a}^{7}{b}^{7}+4973598433723032\,
{K}^{3}{a}^{6}{b}^{8}-224614488424755\,{K}^{3}{a}^{5}{b}^{9}-\\
&
791926314693483\,{K}^{3}{a}^{4}{b}^{10}-246815727308231\,{K}^{3}{a}^{3
}{b}^{11}-23712223279200\,{K}^{3}{a}^{2}{b}^{12}+1664695605143\,{K}^{3
}a{b}^{13}+\\
&324351166694\,{K}^{3}{b}^{14}-6573603371344\,{K}^{2}{a}^{15
}-49410807620672\,{K}^{2}{a}^{14}b-116401469323350\,{K}^{2}{a}^{13}{b}
^{2}-106493810435854\,{K}^{2}{a}^{12}{b}^{3}\\
&+227197756805605\,{K}^{2}{
a}^{11}{b}^{4}+1157143458766160\,{K}^{2}{a}^{10}{b}^{5}+
1682429928684150\,{K}^{2}{a}^{9}{b}^{6}+1068539216327870\,{K}^{2}{a}^{
8}{b}^{7}+\\
&178605190911666\,{K}^{2}{a}^{7}{b}^{8}-158020103057932\,{K}^
{2}{a}^{6}{b}^{9}-103202984052664\,{K}^{2}{a}^{5}{b}^{10}-
23132357064922\,{K}^{2}{a}^{4}{b}^{11}-\\
&963219341342\,{K}^{2}{a}^{3}{b}
^{12}+368416083517\,{K}^{2}{a}^{2}{b}^{13}+36348776217\,{K}^{2}a{b}^{
14}-713149515\,{K}^{2}{b}^{15}-310677290800\,K{a}^{16}-
\end{aligned}
\end{scriptsize}
$$
$$
\begin{scriptsize}
\begin{aligned}
\qquad &963191484776\,K
{a}^{15}b+1869125870004\,K{a}^{14}{b}^{2}+13157968164222\,K{a}^{13}{b}
^{3}+38441627982918\,K{a}^{12}{b}^{4}+75872939747004\,K{a}^{11}{b}^{5}
+\\
&83555042905797\,K{a}^{10}{b}^{6}+46045384280516\,K{a}^{9}{b}^{7}+
5688596768885\,K{a}^{8}{b}^{8}-7709262373734\,K{a}^{7}{b}^{9}-
4737374864728\,K{a}^{6}{b}^{10}-\\
&1122794999625\,K{a}^{5}{b}^{11}-
85849180846\,K{a}^{4}{b}^{12}+7862263012\,K{a}^{3}{b}^{13}+1082607681
\,K{a}^{2}{b}^{14}-12056327\,Ka{b}^{15}-28708\,K{b}^{16}+\\
&3290945392\,{
a}^{17}+61463881652\,{a}^{16}b+297760143344\,{a}^{15}{b}^{2}+
745892371932\,{a}^{14}{b}^{3}+1372640215129\,{a}^{13}{b}^{4}+
2032496622804\,{a}^{12}{b}^{5}+\\
&2109078302595\,{a}^{11}{b}^{6}+
1401344859150\,{a}^{10}{b}^{7}+551690342179\,{a}^{9}{b}^{8}+
100856078874\,{a}^{8}{b}^{9}-6687048368\,{a}^{7}{b}^{10}-6457273203\,{
a}^{6}{b}^{11}-\\
&915591190\,{a}^{5}{b}^{12}-7170228\,{a}^{4}{b}^{13}+
3985173\,{a}^{3}{b}^{14}+21531\,{a}^{2}{b}^{15}-137300806080\,{K}^{14}
{a}^{2}-70999227860\,{K}^{14}ab-\\
&10330598764\,{K}^{14}{b}^{2}-
7716460351184\,{K}^{13}{a}^{3}+305935831784\,{K}^{13}{a}^{2}b+
3437887955868\,{K}^{13}a{b}^{2}+601273750899\,{K}^{13}{b}^{3}-\\
&
32016723369648\,{K}^{12}{a}^{4}+242771178616752\,{K}^{12}{a}^{3}b+
199570628377350\,{K}^{12}{a}^{2}{b}^{2}+26264194126371\,{K}^{12}a{b}^{
3}-770060392865\,{K}^{12}{b}^{4}\\
&+1192188381120040\,{K}^{11}{a}^{5}+
3890792056536405\,{K}^{11}{a}^{4}b+1494821278795288\,{K}^{11}{a}^{3}{b
}^{2}-758673405941100\,{K}^{11}{a}^{2}{b}^{3}-\\
&376869093016859\,{K}^{11
}a{b}^{4}-38203576082594\,{K}^{11}{b}^{5}+10531470692119248\,{K}^{10}{
a}^{6}+13339794799160492\,{K}^{10}{a}^{5}b-\\
&13008290240273101\,{K}^{10}
{a}^{4}{b}^{2}-19187393206564211\,{K}^{10}{a}^{3}{b}^{3}-
6092039292983170\,{K}^{10}{a}^{2}{b}^{4}-369353926241752\,{K}^{10}a{b}
^{5}+\\
&46071954837733\,{K}^{10}{b}^{6}+22734550849908424\,{K}^{9}{a}^{7}
-40855997523289843\,{K}^{9}{a}^{6}b-147373271077439705\,{K}^{9}{a}^{5}
{b}^{2}-\\
&99243967184484638\,{K}^{9}{a}^{4}{b}^{3}-7017041302352190\,{K}
^{9}{a}^{3}{b}^{4}+11638020047399391\,{K}^{9}{a}^{2}{b}^{5}+
3329852098575535\,{K}^{9}a{b}^{6}+\\
&252721176794942\,{K}^{9}{b}^{7}-
5086191849139816\,{K}^{8}{a}^{8}-218840180616499058\,{K}^{8}{a}^{7}b-
268209127820607841\,{K}^{8}{a}^{6}{b}^{2}+\\
&74692704519889508\,{K}^{8}{a
}^{5}{b}^{3}+255648875149961271\,{K}^{8}{a}^{4}{b}^{4}+
141409266020791897\,{K}^{8}{a}^{3}{b}^{5}+30638417717393858\,{K}^{8}{a
}^{2}{b}^{6}+\\
&2101378504207296\,{K}^{8}a{b}^{7}-51741624407282\,{K}^{8}
{b}^{8}-4233155861784192\,{K}^{7}{a}^{9}-33993396202810107\,{K}^{7}{a}
^{8}b+\\
&490749736219117708\,{K}^{7}{a}^{7}{b}^{2}+1105673262827306517\,{
K}^{7}{a}^{6}{b}^{3}+860602010226395943\,{K}^{7}{a}^{5}{b}^{4}+
243473222608791141\,{K}^{7}{a}^{4}{b}^{5}-\\
&25146479261576095\,{K}^{7}{a
}^{3}{b}^{6}-29606343001992830\,{K}^{7}{a}^{2}{b}^{7}-5988530544782415
\,{K}^{7}a{b}^{8}-390821426720579\,{K}^{7}{b}^{9}+\\
&3549684638564912\,{K
}^{6}{a}^{10}+39774183051797904\,{K}^{6}{a}^{9}b+264936658589440290\,{
K}^{6}{a}^{8}{b}^{2}-23253124498207591\,{K}^{6}{a}^{7}{b}^{3}-\\
&
836053898052980210\,{K}^{6}{a}^{6}{b}^{4}-1084726865057884595\,{K}^{6}
{a}^{5}{b}^{5}-636423483060620483\,{K}^{6}{a}^{4}{b}^{6}-
196479957133281562\,{K}^{6}{a}^{3}{b}^{7}-\\
&30869425922772329\,{K}^{6}{a
}^{2}{b}^{8}-1918482700649038\,{K}^{6}a{b}^{9}+5002639556868\,{K}^{6}{
b}^{10}+595927387647512\,{K}^{5}{a}^{11}+\\
&741619486053055\,{K}^{5}{a}^{
10}b-23149404108027253\,{K}^{5}{a}^{9}{b}^{2}-294678512276563627\,{K}^
{5}{a}^{8}{b}^{3}-459961595068021312\,{K}^{5}{a}^{7}{b}^{4}-\\
&
158552034362702252\,{K}^{5}{a}^{6}{b}^{5}+191629479153897981\,{K}^{5}{
a}^{5}{b}^{6}+225899297451072438\,{K}^{5}{a}^{4}{b}^{7}+
104039864883373946\,{K}^{5}{a}^{3}{b}^{8}+\\
&24978090414519088\,{K}^{5}{a
}^{2}{b}^{9}+3069239124835531\,{K}^{5}a{b}^{10}+151737846709549\,{K}^{
5}{b}^{11}-579061113273400\,{K}^{4}{a}^{12}-\\
&4700502461022914\,{K}^{4}{
a}^{11}b-20568042022593357\,{K}^{4}{a}^{10}{b}^{2}-60274756704193639\,
{K}^{4}{a}^{9}{b}^{3}-3101750325638679\,{K}^{4}{a}^{8}{b}^{4}+\\
&
127711457076060492\,{K}^{4}{a}^{7}{b}^{5}+151814783913016234\,{K}^{4}{
a}^{6}{b}^{6}+68971573407733366\,{K}^{4}{a}^{5}{b}^{7}+
4225213632019738\,{K}^{4}{a}^{4}{b}^{8}-\\
&8454619259063026\,{K}^{4}{a}^{
3}{b}^{9}-3543360611949354\,{K}^{4}{a}^{2}{b}^{10}-591871354551764\,{K
}^{4}a{b}^{11}-37188099932363\,{K}^{4}{b}^{12}-\\
&140161902453640\,{K}^{3
}{a}^{13}-536808732218418\,{K}^{3}{a}^{12}b-835239284585174\,{K}^{3}{a
}^{11}{b}^{2}+1749854428584227\,{K}^{3}{a}^{10}{b}^{3}+\\
&
20009652118969842\,{K}^{3}{a}^{9}{b}^{4}+33260818632285462\,{K}^{3}{a}
^{8}{b}^{5}+19240469188117760\,{K}^{3}{a}^{7}{b}^{6}-1551389218247314
\,{K}^{3}{a}^{6}{b}^{7}-\\
&7333550810124668\,{K}^{3}{a}^{5}{b}^{8}-
3615603839244534\,{K}^{3}{a}^{4}{b}^{9}-667099642921947\,{K}^{3}{a}^{3
}{b}^{10}+14428865703118\,{K}^{3}{a}^{2}{b}^{11}+\\
&21146813446302\,{K}^{
3}a{b}^{12}+2035451816443\,{K}^{3}{b}^{13}-1882772930176\,{K}^{2}{a}^{
14}+66733263713456\,{K}^{2}{a}^{13}b+344482894661647\,{K}^{2}{a}^{12}{
b}^{2}+\\
&1087169055246268\,{K}^{2}{a}^{11}{b}^{3}+2600689850371348\,{K}^
{2}{a}^{10}{b}^{4}+2442191086809981\,{K}^{2}{a}^{9}{b}^{5}+
40589252330701\,{K}^{2}{a}^{8}{b}^{6}-\\
&1580272837369054\,{K}^{2}{a}^{7}
{b}^{7}-1218934532449430\,{K}^{2}{a}^{6}{b}^{8}-370704879164434\,{K}^{
2}{a}^{5}{b}^{9}-15022616065481\,{K}^{2}{a}^{4}{b}^{10}+\\
&18590923075742
\,{K}^{2}{a}^{3}{b}^{11}+3999074018231\,{K}^{2}{a}^{2}{b}^{12}+
148009937590\,{K}^{2}a{b}^{13}-17510672155\,{K}^{2}{b}^{14}+
1180651440784\,K{a}^{15}+\\
&10331219189092\,K{a}^{14}b+30751924953232\,K{
a}^{13}{b}^{2}+64074353547984\,K{a}^{12}{b}^{3}+100404901353600\,K{a}^
{11}{b}^{4}+\\
&38443036612706\,K{a}^{10}{b}^{5}-89368546584226\,K{a}^{9}{
b}^{6}-131090693485656\,K{a}^{8}{b}^{7}-75124065103040\,K{a}^{7}{b}^{8
}-18245551274216\,K{a}^{6}{b}^{9}+\\
&671246901080\,K{a}^{5}{b}^{10}+
1334697811842\,K{a}^{4}{b}^{11}+241519135774\,K{a}^{3}{b}^{12}+
7444599539\,K{a}^{2}{b}^{13}-974813822\,Ka{b}^{14}+2740056\,K{b}^{15}\\
&+
48195689744\,{a}^{16}+249404555204\,{a}^{15}b+436134213556\,{a}^{14}{b
}^{2}+531558923976\,{a}^{13}{b}^{3}+441559680121\,{a}^{12}{b}^{4}-\\
&
1180898729833\,{a}^{11}{b}^{5}-3355686082906\,{a}^{10}{b}^{6}-
3572326521080\,{a}^{9}{b}^{7}-2019293841156\,{a}^{8}{b}^{8}-
616288403582\,{a}^{7}{b}^{9}-\\
&77948704502\,{a}^{6}{b}^{10}+6083547961\,
{a}^{5}{b}^{11}+2598909529\,{a}^{4}{b}^{12}+146370877\,{a}^{3}{b}^{13}
-5382340\,{a}^{2}{b}^{14}-21531\,a{b}^{15}+31311936\,{K}^{14}a-\\
&
233545852\,{K}^{14}b+966225874128\,{K}^{13}{a}^{2}+941530601520\,{K}^{
13}ab+155420643480\,{K}^{13}{b}^{2}+41081000247952\,{K}^{12}{a}^{3}+\\
&
25877600391112\,{K}^{12}{a}^{2}b-5465913679753\,{K}^{12}a{b}^{2}-
1908653950478\,{K}^{12}{b}^{3}+296661955152208\,{K}^{11}{a}^{4}-
334213561312056\,{K}^{11}{a}^{3}b-\\
&588931540051323\,{K}^{11}{a}^{2}{b}^
{2}-170704745386947\,{K}^{11}a{b}^{3}-13190965731240\,{K}^{11}{b}^{4}-
1393094881364988\,{K}^{10}{a}^{5}-\\
&8260171869295074\,{K}^{10}{a}^{4}b-
6757615359344263\,{K}^{10}{a}^{3}{b}^{2}-879458173306536\,{K}^{10}{a}^
{2}{b}^{3}+376128503121118\,{K}^{10}a{b}^{4}+\\
&73085367538555\,{K}^{10}{
b}^{5}-17385491256982476\,{K}^{9}{a}^{6}-38111834216188087\,{K}^{9}{a}
^{5}b-7145652655765302\,{K}^{9}{a}^{4}{b}^{2}+\\
&20132690630012112\,{K}^{
9}{a}^{3}{b}^{3}+11954103505398883\,{K}^{9}{a}^{2}{b}^{4}+
2205287792979065\,{K}^{9}a{b}^{5}+115913612400164\,{K}^{9}{b}^{6}-
\end{aligned}
\end{scriptsize}
$$
$$
\begin{scriptsize}
\begin{aligned}
\qquad &40214067009703812\,{K}^{8}{a}^{7}+131639843169544\,{K}^{8}{a}^{6}b+
156558622175754131\,{K}^{8}{a}^{5}{b}^{2}+174156412173136580\,{K}^{8}{
a}^{4}{b}^{3}+\\
&66450821963204562\,{K}^{8}{a}^{3}{b}^{4}+
6095034005529188\,{K}^{8}{a}^{2}{b}^{5}-1488474504978533\,{K}^{8}a{b}^
{6}-244373592177005\,{K}^{8}{b}^{7}-\\
&4262609640794292\,{K}^{7}{a}^{8}+
232044622223120291\,{K}^{7}{a}^{7}b+450566932400389617\,{K}^{7}{a}^{6}
{b}^{2}+223907430855339989\,{K}^{7}{a}^{5}{b}^{3}-\\
&76252311956570662\,{
K}^{7}{a}^{4}{b}^{4}-109626476115719167\,{K}^{7}{a}^{3}{b}^{5}-
38945629332253468\,{K}^{7}{a}^{2}{b}^{6}-5727951900575235\,{K}^{7}a{b}
^{7}-\\
&291248294307704\,{K}^{7}{b}^{8}+5817530118780564\,{K}^{6}{a}^{9}+
90456025289724496\,{K}^{6}{a}^{8}b-211722802621223320\,{K}^{6}{a}^{7}{
b}^{2}-\\
&804260518427049252\,{K}^{6}{a}^{6}{b}^{3}-869733715790870347\,{
K}^{6}{a}^{5}{b}^{4}-438230331370874864\,{K}^{6}{a}^{4}{b}^{5}-
103383746035998504\,{K}^{6}{a}^{3}{b}^{6}-\\
&6367471759271661\,{K}^{6}{a}
^{2}{b}^{7}+1483742496739426\,{K}^{6}a{b}^{8}+198232240580240\,{K}^{6}
{b}^{9}-921628104356508\,{K}^{5}{a}^{10}-\\
&14675339280309847\,{K}^{5}{a}
^{9}b-212169540522227624\,{K}^{5}{a}^{8}{b}^{2}-242390052256145101\,{K
}^{5}{a}^{7}{b}^{3}+142387137490652383\,{K}^{5}{a}^{6}{b}^{4}+\\
&
437705095527040895\,{K}^{5}{a}^{5}{b}^{5}+346770220231186390\,{K}^{5}{
a}^{4}{b}^{6}+138494478959211193\,{K}^{5}{a}^{3}{b}^{7}+
30073038024652238\,{K}^{5}{a}^{2}{b}^{8}+\\
&3331095984744577\,{K}^{5}a{b}
^{9}+143933703893023\,{K}^{5}{b}^{10}-618571777557828\,{K}^{4}{a}^{11}
-5356142866137578\,{K}^{4}{a}^{10}b-\\
&24708287467462518\,{K}^{4}{a}^{9}{
b}^{2}+70712016882151816\,{K}^{4}{a}^{8}{b}^{3}+209044763547371381\,{K
}^{4}{a}^{7}{b}^{4}+175324553438318195\,{K}^{4}{a}^{6}{b}^{5}+\\
&
33620373342382331\,{K}^{4}{a}^{5}{b}^{6}-37350348148744791\,{K}^{4}{a}
^{4}{b}^{7}-28616057049689490\,{K}^{4}{a}^{3}{b}^{8}-8747338664124168
\,{K}^{4}{a}^{2}{b}^{9}-\\
&1290791335967142\,{K}^{4}a{b}^{10}-
75627600062918\,{K}^{4}{b}^{11}+119150612609500\,{K}^{3}{a}^{12}+
833840582698667\,{K}^{3}{a}^{11}b+\\
&4579583581115468\,{K}^{3}{a}^{10}{b}
^{2}+26469740772937682\,{K}^{3}{a}^{9}{b}^{3}+31640137736783012\,{K}^{
3}{a}^{8}{b}^{4}-1392010795845594\,{K}^{3}{a}^{7}{b}^{5}-\\
&
27774586393551336\,{K}^{3}{a}^{6}{b}^{6}-21979728328610634\,{K}^{3}{a}
^{5}{b}^{7}-6852220022473196\,{K}^{3}{a}^{4}{b}^{8}-221287702147389\,{
K}^{3}{a}^{3}{b}^{9}+\\
&418502522351477\,{K}^{3}{a}^{2}{b}^{10}+
103476980956710\,{K}^{3}a{b}^{11}+7853645364594\,{K}^{3}{b}^{12}+
42080378186256\,{K}^{2}{a}^{13}+\\
&190032726967756\,{K}^{2}{a}^{12}b+
646961067553627\,{K}^{2}{a}^{11}{b}^{2}+1832518528326428\,{K}^{2}{a}^{
10}{b}^{3}-355556998553669\,{K}^{2}{a}^{9}{b}^{4}-\\
&4834936965943107\,{K
}^{2}{a}^{8}{b}^{5}-5403946926848664\,{K}^{2}{a}^{7}{b}^{6}-
2311436963468304\,{K}^{2}{a}^{6}{b}^{7}-28837070317252\,{K}^{2}{a}^{5}
{b}^{8}+\\
&338780578433516\,{K}^{2}{a}^{4}{b}^{9}+119453042900413\,{K}^{2
}{a}^{3}{b}^{10}+13007062996650\,{K}^{2}{a}^{2}{b}^{11}-716455710136\,
{K}^{2}a{b}^{12}-\\
&169649709771\,{K}^{2}{b}^{13}+1660767110752\,K{a}^{14
}-2207968158376\,K{a}^{13}b-13696747630000\,K{a}^{12}{b}^{2}-
39597663442185\,K{a}^{11}{b}^{3}-\\
&250679658638620\,K{a}^{10}{b}^{4}-
443148312946451\,K{a}^{9}{b}^{5}-318849573722046\,K{a}^{8}{b}^{6}-
60459908744136\,K{a}^{7}{b}^{7}+\\
&52843969331026\,K{a}^{6}{b}^{8}+
37295884685596\,K{a}^{5}{b}^{9}+8932801190136\,K{a}^{4}{b}^{10}+
437851974124\,K{a}^{3}{b}^{11}-143471711864\,K{a}^{2}{b}^{12}-\\
&
15117869866\,Ka{b}^{13}+288406903\,K{b}^{14}-82866047536\,{a}^{15}-
831590288964\,{a}^{14}b-2265215755284\,{a}^{13}{b}^{2}-4220031694728\,
{a}^{12}{b}^{3}-\\
&9069845696385\,{a}^{11}{b}^{4}-11321910217349\,{a}^{10
}{b}^{5}-6269545516932\,{a}^{9}{b}^{6}+135893147933\,{a}^{8}{b}^{7}+
2118624548324\,{a}^{7}{b}^{8}+\\
&1178848132690\,{a}^{6}{b}^{9}+
281051305918\,{a}^{5}{b}^{10}+22523506033\,{a}^{4}{b}^{11}-1894139614
\,{a}^{3}{b}^{12}-276297630\,{a}^{2}{b}^{13}+2987168\,a{b}^{14}+7177\,
{b}^{15}\\
&+50891409424\,{K}^{13}a+11615301204\,{K}^{13}b-1723217091152\,
{K}^{12}{a}^{2}-3641406369942\,{K}^{12}ab-717406251381\,{K}^{12}{b}^{2
}-\\
&106611546435664\,{K}^{11}{a}^{3}-125631676152171\,{K}^{11}{a}^{2}b-
22276913546612\,{K}^{11}a{b}^{2}+149932087428\,{K}^{11}{b}^{3}-
921477476876608\,{K}^{10}{a}^{4}-\\
&335321749015788\,{K}^{10}{a}^{3}b+
675513168416462\,{K}^{10}{a}^{2}{b}^{2}+343951420930887\,{K}^{10}a{b}^
{3}+42312060567697\,{K}^{10}{b}^{4}-\\
&224822721439848\,{K}^{9}{a}^{5}+
9113703426775205\,{K}^{9}{a}^{4}b+12011579714385819\,{K}^{9}{a}^{3}{b}
^{2}+4545321284411010\,{K}^{9}{a}^{2}{b}^{3}+\\
&459286827939476\,{K}^{9}a
{b}^{4}-16563525204644\,{K}^{9}{b}^{5}+18269988327740104\,{K}^{8}{a}^{
6}+57345016724065526\,{K}^{8}{a}^{5}b+\\
&43380938324270807\,{K}^{8}{a}^{4
}{b}^{2}+2016283661002179\,{K}^{8}{a}^{3}{b}^{3}-7866364168161881\,{K}
^{8}{a}^{2}{b}^{4}-2536815806374494\,{K}^{8}a{b}^{5}-\\
&229611275189739\,
{K}^{8}{b}^{6}+48204132263729968\,{K}^{7}{a}^{7}+62130773991526781\,{K
}^{7}{a}^{6}b-64408460228007020\,{K}^{7}{a}^{5}{b}^{2}-\\
&
143588285939373137\,{K}^{7}{a}^{4}{b}^{3}-86865297183814146\,{K}^{7}{a
}^{3}{b}^{4}-22343649081768145\,{K}^{7}{a}^{2}{b}^{5}-2237789223581738
\,{K}^{7}a{b}^{6}-\\
&34917319832253\,{K}^{7}{b}^{7}+14906794428731648\,{K
}^{6}{a}^{8}-151148659919697864\,{K}^{6}{a}^{7}b-411536847144982190\,{
K}^{6}{a}^{6}{b}^{2}-\\
&370097713378620282\,{K}^{6}{a}^{5}{b}^{3}-
122511186080227620\,{K}^{6}{a}^{4}{b}^{4}+9555876039698707\,{K}^{6}{a}
^{3}{b}^{5}+16423510144303433\,{K}^{6}{a}^{2}{b}^{6}+\\
&3820620328826261
\,{K}^{6}a{b}^{7}+286660719055213\,{K}^{6}{b}^{8}-2610796839447768\,{K
}^{5}{a}^{9}-87401641373975649\,{K}^{5}{a}^{8}b-\\
&13874478240722049\,{K}
^{5}{a}^{7}{b}^{2}+311353562263468600\,{K}^{5}{a}^{6}{b}^{3}+
472127870330715962\,{K}^{5}{a}^{5}{b}^{4}+314216093424549000\,{K}^{5}{
a}^{4}{b}^{5}+\\
&110637300590115018\,{K}^{5}{a}^{3}{b}^{6}+
20561305908347595\,{K}^{5}{a}^{2}{b}^{7}+1746436083406920\,{K}^{5}a{b}
^{8}+38493099860522\,{K}^{5}{b}^{9}-\\
&455230966520168\,{K}^{4}{a}^{10}-
5801625691725300\,{K}^{4}{a}^{9}b+83093264190707412\,{K}^{4}{a}^{8}{b}
^{2}+169318093020829647\,{K}^{4}{a}^{7}{b}^{3}+\\
&78545237600244913\,{K}^
{4}{a}^{6}{b}^{4}-68045554649666063\,{K}^{4}{a}^{5}{b}^{5}-
98595796539665656\,{K}^{4}{a}^{4}{b}^{6}-50809862654400173\,{K}^{4}{a}
^{3}{b}^{7}-\\
&13473876719245864\,{K}^{4}{a}^{2}{b}^{8}-1829146000444877
\,{K}^{4}a{b}^{9}-100608385231027\,{K}^{4}{b}^{10}+146010094862064\,{K
}^{3}{a}^{11}+\\
&1782900675391474\,{K}^{3}{a}^{10}b+17576459949313190\,{K
}^{3}{a}^{9}{b}^{2}+6573177586457609\,{K}^{3}{a}^{8}{b}^{3}-
38834161644773778\,{K}^{3}{a}^{7}{b}^{4}-\\
&55196631271773857\,{K}^{3}{a}
^{6}{b}^{5}-28653877190479683\,{K}^{3}{a}^{5}{b}^{6}-2944935310456052
\,{K}^{3}{a}^{4}{b}^{7}+3371369752069724\,{K}^{3}{a}^{3}{b}^{8}+\\
&
1608968476689814\,{K}^{3}{a}^{2}{b}^{9}+292072474434394\,{K}^{3}a{b}^{
10}+19716626463697\,{K}^{3}{b}^{11}-16522558463216\,{K}^{2}{a}^{12}+\\
&
20485835386584\,{K}^{2}{a}^{11}b+432107080967014\,{K}^{2}{a}^{10}{b}^{
2}-3928401281891521\,{K}^{2}{a}^{9}{b}^{3}-8840897211758225\,{K}^{2}{a
}^{8}{b}^{4}-\\
&5973620079452524\,{K}^{2}{a}^{7}{b}^{5}+281481293270066\,
{K}^{2}{a}^{6}{b}^{6}+2450757075317472\,{K}^{2}{a}^{5}{b}^{7}+
1321492213423618\,{K}^{2}{a}^{4}{b}^{8}+
\end{aligned}
\end{scriptsize}
$$
$$
\begin{scriptsize}
\begin{aligned}
\qquad &266435769885126\,{K}^{2}{a}^{3
}{b}^{9}-2811134174860\,{K}^{2}{a}^{2}{b}^{10}-8571743336077\,{K}^{2}a
{b}^{11}-882010936046\,{K}^{2}{b}^{12}-6834964299808\,K{a}^{13}-\\
&
21883995529772\,K{a}^{12}b-60265969337672\,K{a}^{11}{b}^{2}-
365304924430456\,K{a}^{10}{b}^{3}-449369348564134\,K{a}^{9}{b}^{4}+
9425121998155\,K{a}^{8}{b}^{5}+\\
&390569580462290\,K{a}^{7}{b}^{6}+
317251342237782\,K{a}^{6}{b}^{7}+101265060459192\,K{a}^{5}{b}^{8}+
3954386963972\,K{a}^{4}{b}^{9}-5736801991082\,K{a}^{3}{b}^{10}\\
&-1278647277536\,K{a}^{2}{b}^{11}-49634209706\,Ka{b}^{12}+5883396002\,K{
b}^{13}-218271124832\,{a}^{14}-78278981072\,{a}^{13}b+131693094764\,{a
}^{12}{b}^{2}\\
&-4352390062248\,{a}^{11}{b}^{3}-313920381978\,{a}^{10}{b}
^{4}+13870135960045\,{a}^{9}{b}^{5}+18976137756795\,{a}^{8}{b}^{6}+
10925040635808\,{a}^{7}{b}^{7}+\\
&2500684816842\,{a}^{6}{b}^{8}-
283130954660\,{a}^{5}{b}^{9}-271343648696\,{a}^{4}{b}^{10}-47975867793
\,{a}^{3}{b}^{11}-1497673484\,{a}^{2}{b}^{12}+196492312\,a{b}^{13}-\\
&
543705\,{b}^{14}-331815807344\,{K}^{12}a-76467775798\,{K}^{12}b-
2790297882264\,{K}^{11}{a}^{2}+4760473687220\,{K}^{11}ab+1257339729339
\,{K}^{11}{b}^{2}+\\
&140090606088644\,{K}^{10}{a}^{3}+252813156465054\,{K
}^{10}{a}^{2}b+90347442696682\,{K}^{10}a{b}^{2}+9282403825083\,{K}^{10
}{b}^{3}+1523235978285556\,{K}^{9}{a}^{4}\\
&+1747323682264597\,{K}^{9}{a}
^{3}b+130106919840946\,{K}^{9}{a}^{2}{b}^{2}-230089975479143\,{K}^{9}a
{b}^{3}-44560927895597\,{K}^{9}{b}^{4}+\\
&3032160274139372\,{K}^{8}{a}^{5
}-3260374132158388\,{K}^{8}{a}^{4}b-10030885667822625\,{K}^{8}{a}^{3}{
b}^{2}-5894264280303247\,{K}^{8}{a}^{2}{b}^{3}-\\
&1239663820065337\,{K}^{
8}a{b}^{4}-82867030162130\,{K}^{8}{b}^{5}-11207606290403492\,{K}^{7}{a
}^{6}-49732171270576889\,{K}^{7}{a}^{5}b-\\
&57628615946520527\,{K}^{7}{a}
^{4}{b}^{2}-24432046441639155\,{K}^{7}{a}^{3}{b}^{3}-2496629431034273
\,{K}^{7}{a}^{2}{b}^{4}+680228031819000\,{K}^{7}a{b}^{5}+\\
&
124037907644439\,{K}^{7}{b}^{6}-38899037836711300\,{K}^{6}{a}^{7}-
86900634074912402\,{K}^{6}{a}^{6}b-36002532611395668\,{K}^{6}{a}^{5}{b
}^{2}+\\
&43898139679411123\,{K}^{6}{a}^{4}{b}^{3}+48327630886662571\,{K}^
{6}{a}^{3}{b}^{4}+18190064053584356\,{K}^{6}{a}^{2}{b}^{5}+
2994901429609894\,{K}^{6}a{b}^{6}+\\
&179545495072828\,{K}^{6}{b}^{7}-
16682998481963588\,{K}^{5}{a}^{8}+53366413954257167\,{K}^{5}{a}^{7}b+
220683169365287680\,{K}^{5}{a}^{6}{b}^{2}+\\
&267107353451448127\,{K}^{5}{
a}^{5}{b}^{3}+150891749838209906\,{K}^{5}{a}^{4}{b}^{4}+
40634275294646461\,{K}^{5}{a}^{3}{b}^{5}+3461483753509527\,{K}^{5}{a}^
{2}{b}^{6}-\\
&491748467655435\,{K}^{5}a{b}^{7}-83428128113164\,{K}^{5}{b}
^{8}-805195408290764\,{K}^{4}{a}^{9}+44872401191679682\,{K}^{4}{a}^{8}
b+\\
&60111338170463787\,{K}^{4}{a}^{7}{b}^{2}-42458774294396965\,{K}^{4}{
a}^{6}{b}^{3}-139226880011355943\,{K}^{4}{a}^{5}{b}^{4}-
120744538465589397\,{K}^{4}{a}^{4}{b}^{5}-\\
&52963024951441167\,{K}^{4}{a
}^{3}{b}^{6}-12712340491161802\,{K}^{4}{a}^{2}{b}^{7}-1577697169146556
\,{K}^{4}a{b}^{8}-78518959226643\,{K}^{4}{b}^{9}+\\
&254364149522892\,{K}^
{3}{a}^{10}+6909699638376441\,{K}^{3}{a}^{9}b-13278665232723742\,{K}^{
3}{a}^{8}{b}^{2}-51346447417517381\,{K}^{3}{a}^{7}{b}^{3}-\\
&
48110141215772905\,{K}^{3}{a}^{6}{b}^{4}-10132849652978844\,{K}^{3}{a}
^{5}{b}^{5}+11989320109516061\,{K}^{3}{a}^{4}{b}^{6}+9919526555896783
\,{K}^{3}{a}^{3}{b}^{7}+\\
&3262772207324032\,{K}^{3}{a}^{2}{b}^{8}+
517657351374835\,{K}^{3}a{b}^{9}+32603871364075\,{K}^{3}{b}^{10}+
7177700519536\,{K}^{2}{a}^{11}+\\
&105158309268388\,{K}^{2}{a}^{10}b-
4265957633573160\,{K}^{2}{a}^{9}{b}^{2}-6394519529905170\,{K}^{2}{a}^{
8}{b}^{3}+442297630926712\,{K}^{2}{a}^{7}{b}^{4}+\\
&6983258017463205\,{K}
^{2}{a}^{6}{b}^{5}+5879351836026305\,{K}^{2}{a}^{5}{b}^{6}+
1944801114291594\,{K}^{2}{a}^{4}{b}^{7}+64266774421718\,{K}^{2}{a}^{3}
{b}^{8}-\\
&134235549150456\,{K}^{2}{a}^{2}{b}^{9}-34947156999062\,{K}^{2}
a{b}^{10}-2788274935124\,{K}^{2}{b}^{11}+4197138751392\,K{a}^{12}-
588050753616\,K{a}^{11}b-\\
&243382331903712\,K{a}^{10}{b}^{2}+
10374396265872\,K{a}^{9}{b}^{3}+784594573752836\,K{a}^{8}{b}^{4}+
959487257647224\,K{a}^{7}{b}^{5}+\\
&421999182743562\,K{a}^{6}{b}^{6}-
13853625159704\,K{a}^{5}{b}^{7}-81631066805374\,K{a}^{4}{b}^{8}-
28920031307866\,K{a}^{3}{b}^{9}-3199106510532\,K{a}^{2}{b}^{10}+\\
&
200901041755\,Ka{b}^{11}+46733544880\,K{b}^{12}+599874377376\,{a}^{13}
+1346228911576\,{a}^{12}b-1891789990720\,{a}^{11}{b}^{2}+\\
&
13865058512804\,{a}^{10}{b}^{3}+35027217927938\,{a}^{9}{b}^{4}+
26681722029374\,{a}^{8}{b}^{5}+2205534004606\,{a}^{7}{b}^{6}-
8310149282404\,{a}^{6}{b}^{7}-\\
&5159144609900\,{a}^{5}{b}^{8}-
1196466842268\,{a}^{4}{b}^{9}-42486346632\,{a}^{3}{b}^{10}+24433073863
\,{a}^{2}{b}^{11}+2461495059\,a{b}^{12}-48205412\,{b}^{13}+\\
&
833734201368\,{K}^{11}a+200376725161\,{K}^{11}b+15113330502352\,{K}^{
10}{a}^{2}+2102278604160\,{K}^{10}ab-330127702027\,{K}^{10}{b}^{2}-\\
&
68175874034856\,{K}^{9}{a}^{3}-244833720758751\,{K}^{9}{a}^{2}b-
124525035959009\,{K}^{9}a{b}^{2}-17255651698575\,{K}^{9}{b}^{3}-
1453299752200920\,{K}^{8}{a}^{4}-\\
&2483515219406910\,{K}^{8}{a}^{3}b-
1086919792328737\,{K}^{8}{a}^{2}{b}^{2}-107267484729990\,{K}^{8}a{b}^{
3}+8308840127846\,{K}^{8}{b}^{4}-\\
&4372368855637312\,{K}^{7}{a}^{5}-
3706209619817403\,{K}^{7}{a}^{4}b+2531012245507132\,{K}^{7}{a}^{3}{b}^
{2}+3199732525381513\,{K}^{7}{a}^{2}{b}^{3}+\\
&978283327372673\,{K}^{7}a{
b}^{4}+94805458896700\,{K}^{7}{b}^{5}+2654551452016288\,{K}^{6}{a}^{6}
+23990690978815296\,{K}^{6}{a}^{5}b+\\
&38235760396857966\,{K}^{6}{a}^{4}{
b}^{2}+24039633280008185\,{K}^{6}{a}^{3}{b}^{3}+6769684712645941\,{K}^
{6}{a}^{2}{b}^{4}+770995003754564\,{K}^{6}a{b}^{5}+\\
&18559400953944\,{K}
^{6}{b}^{6}+20952236279031096\,{K}^{5}{a}^{7}+62046854037649907\,{K}^{
5}{a}^{6}b+60589903701206363\,{K}^{5}{a}^{5}{b}^{2}+\\
&17602296566710185
\,{K}^{5}{a}^{4}{b}^{3}-6934977455660400\,{K}^{5}{a}^{3}{b}^{4}-
5798515195449386\,{K}^{5}{a}^{2}{b}^{5}-1355986782924447\,{K}^{5}a{b}^
{6}-\\
&108653678121387\,{K}^{5}{b}^{7}+10316915681438312\,{K}^{4}{a}^{8}-
4539691137654438\,{K}^{4}{a}^{7}b-68410018711790509\,{K}^{4}{a}^{6}{b}
^{2}-\\
&107769760403040527\,{K}^{4}{a}^{5}{b}^{3}-77982145794302954\,{K}^
{4}{a}^{4}{b}^{4}-30219901352570438\,{K}^{4}{a}^{3}{b}^{5}-
6273825258236269\,{K}^{4}{a}^{2}{b}^{6}-\\
&620542621947919\,{K}^{4}a{b}^{
7}-19600463029102\,{K}^{4}{b}^{8}+1329211200554424\,{K}^{3}{a}^{9}-
12701454995480078\,{K}^{3}{a}^{8}b-\\
&28409991120082018\,{K}^{3}{a}^{7}{b
}^{2}-11888175676997805\,{K}^{3}{a}^{6}{b}^{3}+18325892632230832\,{K}^
{3}{a}^{5}{b}^{4}+25365969878241814\,{K}^{3}{a}^{4}{b}^{5}+\\
&
13796141272443772\,{K}^{3}{a}^{3}{b}^{6}+3923989294278279\,{K}^{3}{a}^
{2}{b}^{7}+575343977003974\,{K}^{3}a{b}^{8}+34360078320678\,{K}^{3}{b}
^{9}+\\
&15257958402528\,{K}^{2}{a}^{10}-2350991541982052\,{K}^{2}{a}^{9}b
-760043059987533\,{K}^{2}{a}^{8}{b}^{2}+7012311225803830\,{K}^{2}{a}^{
7}{b}^{3}+\\
&10290547876831079\,{K}^{2}{a}^{6}{b}^{4}+5488669915643763\,{
K}^{2}{a}^{5}{b}^{5}+402438954013971\,{K}^{2}{a}^{4}{b}^{6}-
875362166196421\,{K}^{2}{a}^{3}{b}^{7}-\\
&418988052136348\,{K}^{2}{a}^{2}
{b}^{8}-79428084646802\,{K}^{2}a{b}^{9}-5636891090451\,{K}^{2}{b}^{10}
-3719303338592\,K{a}^{11}-125198241989132\,K{a}^{10}b+
\end{aligned}
\end{scriptsize}
$$
$$
\begin{scriptsize}
\begin{aligned}
\qquad &367427759494564
\,K{a}^{9}{b}^{2}+1041462628782392\,K{a}^{8}{b}^{3}+705292351331846\,K
{a}^{7}{b}^{4}-152632671328755\,K{a}^{6}{b}^{5}-\\
&443449407162599\,K{a}^
{5}{b}^{6}-234987789436475\,K{a}^{4}{b}^{7}-46961274404656\,K{a}^{3}{b
}^{8}+1493197477268\,K{a}^{2}{b}^{9}+1896974325539\,Ka{b}^{10}+\\
&
196670547417\,K{b}^{11}-527743347424\,{a}^{12}-3131084819288\,{a}^{11}
b+17502326531768\,{a}^{10}{b}^{2}+26511779799944\,{a}^{9}{b}^{3}-\\
&
7367463850226\,{a}^{8}{b}^{4}-36980225208650\,{a}^{7}{b}^{5}-
28880944820860\,{a}^{6}{b}^{6}-8674690794024\,{a}^{5}{b}^{7}+
202581566274\,{a}^{4}{b}^{8}+\\
&760322422442\,{a}^{3}{b}^{9}+157371465658
\,{a}^{2}{b}^{10}+5476055307\,a{b}^{11}-815235128\,{b}^{12}-
1011318037260\,{K}^{10}a-247134263374\,{K}^{10}b-\\
&24581993242412\,{K}^{
9}{a}^{2}-12208307141651\,{K}^{9}ab-1629370229805\,{K}^{9}{b}^{2}-
49437930610788\,{K}^{8}{a}^{3}+88732252082716\,{K}^{8}{a}^{2}b+\\
&
76908223180679\,{K}^{8}a{b}^{2}+13280414273971\,{K}^{8}{b}^{3}+
778947624133500\,{K}^{7}{a}^{4}+1794584809059095\,{K}^{7}{a}^{3}b+\\
&
1141962854108765\,{K}^{7}{a}^{2}{b}^{2}+263546280381661\,{K}^{7}a{b}^{
3}+19528749443108\,{K}^{7}{b}^{4}+3282764773940396\,{K}^{6}{a}^{5}+\\
&
5423043455524678\,{K}^{6}{a}^{4}b+2294549768662472\,{K}^{6}{a}^{3}{b}^
{2}-175900158955480\,{K}^{6}{a}^{2}{b}^{3}-267103231890026\,{K}^{6}a{b
}^{4}-\\
&38261110532684\,{K}^{6}{b}^{5}+1329832353294732\,{K}^{5}{a}^{6}-
4488343686840185\,{K}^{5}{a}^{5}b-12885211421709632\,{K}^{5}{a}^{4}{b}
^{2}-\\
&10935364726466451\,{K}^{5}{a}^{3}{b}^{3}-4165460139552913\,{K}^{5
}{a}^{2}{b}^{4}-733602668249600\,{K}^{5}a{b}^{5}-48402315171979\,{K}^{
5}{b}^{6}-\\
&7408920965841420\,{K}^{4}{a}^{7}-26520063925924244\,{K}^{4}{
a}^{6}b-34848430722441086\,{K}^{4}{a}^{5}{b}^{2}-21328377937291598\,{K
}^{4}{a}^{4}{b}^{3}-\\
&5899395505307096\,{K}^{4}{a}^{3}{b}^{4}-
328109767989118\,{K}^{4}{a}^{2}{b}^{5}+153465217347093\,{K}^{4}a{b}^{6
}+22065161899644\,{K}^{4}{b}^{7}-\\
&3836264164060940\,{K}^{3}{a}^{8}-
3871539447147435\,{K}^{3}{a}^{7}b+10544521639455256\,{K}^{3}{a}^{6}{b}
^{2}+24847350458635066\,{K}^{3}{a}^{5}{b}^{3}+\\
&22078801176961278\,{K}^{
3}{a}^{4}{b}^{4}+10359055572822511\,{K}^{3}{a}^{3}{b}^{5}+
2698635194066088\,{K}^{3}{a}^{2}{b}^{6}+367323597732681\,{K}^{3}a{b}^{
7}+\\
&20323821906886\,{K}^{3}{b}^{8}-565264896005456\,{K}^{2}{a}^{9}+
1785234114339432\,{K}^{2}{a}^{8}b+6207407398520831\,{K}^{2}{a}^{7}{b}^
{2}+\\
&5659079033855276\,{K}^{2}{a}^{6}{b}^{3}+497846569155692\,{K}^{2}{a
}^{5}{b}^{4}-2523695454566092\,{K}^{2}{a}^{4}{b}^{5}-1934498236354658
\,{K}^{2}{a}^{3}{b}^{6}-\\
&655919117148928\,{K}^{2}{a}^{2}{b}^{7}-
109569770690549\,{K}^{2}a{b}^{8}-7313343005164\,{K}^{2}{b}^{9}-
26786456881616\,K{a}^{10}+342579386922976\,K{a}^{9}b+\\
&493807507190380\,
K{a}^{8}{b}^{2}-232472134388689\,K{a}^{7}{b}^{3}-908781845082514\,K{a}
^{6}{b}^{4}-732347171996855\,K{a}^{5}{b}^{5}-\\
&231991134385239\,K{a}^{4}
{b}^{6}+5666667500028\,K{a}^{3}{b}^{7}+24043836485425\,K{a}^{2}{b}^{8}
+6072917895801\,Ka{b}^{9}+495112255696\,K{b}^{10}-\\
&217870735968\,{a}^{
11}+14121088737528\,{a}^{10}b-796720941748\,{a}^{9}{b}^{2}-
47608859543596\,{a}^{8}{b}^{3}-57389080311634\,{a}^{7}{b}^{4}-\\
&
21927246674779\,{a}^{6}{b}^{5}+5763955072681\,{a}^{5}{b}^{6}+
7939066033165\,{a}^{4}{b}^{7}+2596653864697\,{a}^{3}{b}^{8}+
260414813589\,{a}^{2}{b}^{9}-\\
&26874560871\,a{b}^{10}-5156249015\,{b}^{
11}+567904300808\,{K}^{9}a+133164681967\,{K}^{9}b+21014582924024\,{K}^
{8}{a}^{2}+14238813875284\,{K}^{8}ab+\\
&2376108556354\,{K}^{8}{b}^{2}+
96224337465776\,{K}^{7}{a}^{3}+40758471016020\,{K}^{7}{a}^{2}b-
10593073855823\,{K}^{7}a{b}^{2}-3866406950985\,{K}^{7}{b}^{3}-\\
&
172609427644256\,{K}^{6}{a}^{4}-667968429605572\,{K}^{6}{a}^{3}b-
561286043934010\,{K}^{6}{a}^{2}{b}^{2}-170012407763142\,{K}^{6}a{b}^{3
}-17358827265311\,{K}^{6}{b}^{4}\\
&-1453031376276008\,{K}^{5}{a}^{5}-
3157267253199025\,{K}^{5}{a}^{4}b-2316736275650001\,{K}^{5}{a}^{3}{b}^
{2}-697700047477932\,{K}^{5}{a}^{2}{b}^{3}-\\
&73416869055526\,{K}^{5}a{b}
^{4}+45629552406\,{K}^{5}{b}^{5}-1362580068550424\,{K}^{4}{a}^{6}-
1434084336512452\,{K}^{4}{a}^{5}b+1085678954165590\,{K}^{4}{a}^{4}{b}^
{2}+\\
&2145499853923277\,{K}^{4}{a}^{3}{b}^{3}+1129856076985994\,{K}^{4}{
a}^{2}{b}^{4}+253070645685405\,{K}^{4}a{b}^{5}+20942913448202\,{K}^{4}
{b}^{6}+\\
&1657055511607808\,{K}^{3}{a}^{7}+6948020300041926\,{K}^{3}{a}^
{6}b+10935629262529494\,{K}^{3}{a}^{5}{b}^{2}+8571255526529717\,{K}^{3
}{a}^{4}{b}^{3}+\\
&3637799081961792\,{K}^{3}{a}^{3}{b}^{4}+
827633648654025\,{K}^{3}{a}^{2}{b}^{5}+90118434383785\,{K}^{3}a{b}^{6}
+3249523242246\,{K}^{3}{b}^{7}+853240918448912\,{K}^{2}{a}^{8}\\
&+1607156497750060\,{K}^{2}{a}^{7}b-152109263958412\,{K}^{2}{a}^{6}{b}^{
2}-2974544539369497\,{K}^{2}{a}^{5}{b}^{3}-3453815917971886\,{K}^{2}{a
}^{4}{b}^{4}-\\
&1910451825690418\,{K}^{2}{a}^{3}{b}^{5}-574351246968321\,
{K}^{2}{a}^{2}{b}^{6}-90213663849566\,{K}^{2}a{b}^{7}-5804648946916\,{
K}^{2}{b}^{8}+108840267252464\,K{a}^{9}-\\
&69360337481996\,K{a}^{8}b-
632740775706748\,K{a}^{7}{b}^{2}-827434438945328\,K{a}^{6}{b}^{3}-
394158755385236\,K{a}^{5}{b}^{4}+39729765841096\,K{a}^{4}{b}^{5}+\\
&
124555029069677\,K{a}^{3}{b}^{6}+54993285202071\,K{a}^{2}{b}^{7}+
10581054760448\,Ka{b}^{8}+778324633288\,K{b}^{9}+3914330667504\,{a}^{
10}-\\
&17083371284620\,{a}^{9}b-41645257962800\,{a}^{8}{b}^{2}-
20374996143828\,{a}^{7}{b}^{3}+20739585118259\,{a}^{6}{b}^{4}+
30315183992272\,{a}^{5}{b}^{5}+\\
&14695307029065\,{a}^{4}{b}^{6}+
2577468275032\,{a}^{3}{b}^{7}-298295774506\,{a}^{2}{b}^{8}-
169949362649\,a{b}^{9}-16988878331\,{b}^{10}-9359770176\,{K}^{8}a+\\
&
10915809561\,{K}^{8}b-10187011727752\,{K}^{7}{a}^{2}-8311473289500\,{K
}^{7}ab-1525557453524\,{K}^{7}{b}^{2}-65428247692852\,{K}^{6}{a}^{3}-\\
&
58620349631551\,{K}^{6}{a}^{2}b-14398138060749\,{K}^{6}a{b}^{2}-
947647295443\,{K}^{6}{b}^{3}-45218791476740\,{K}^{5}{a}^{4}+
69278485870675\,{K}^{5}{a}^{3}b+\\
&119136834812380\,{K}^{5}{a}^{2}{b}^{2}
+48268108509533\,{K}^{5}a{b}^{3}+5969134650570\,{K}^{5}{b}^{4}+
368670160687220\,{K}^{4}{a}^{5}+981553240283514\,{K}^{4}{a}^{4}b+\\
&
908129310897681\,{K}^{4}{a}^{3}{b}^{2}+379638499758119\,{K}^{4}{a}^{2}
{b}^{3}+72942585399018\,{K}^{4}a{b}^{4}+5242284395526\,{K}^{4}{b}^{5}+
523557687695596\,{K}^{3}{a}^{6}+\\
&1096696458554245\,{K}^{3}{a}^{5}b+
756415710477622\,{K}^{3}{a}^{4}{b}^{2}+131423920464871\,{K}^{3}{a}^{3}
{b}^{3}-64914841434319\,{K}^{3}{a}^{2}{b}^{4}-\\
&29497328559319\,{K}^{3}a
{b}^{5}-3290523531095\,{K}^{3}{b}^{6}-216189731323888\,{K}^{2}{a}^{7}-
1072115048720760\,{K}^{2}{a}^{6}b-1955857227655694\,{K}^{2}{a}^{5}{b}^
{2}\\
&-1787774603372450\,{K}^{2}{a}^{4}{b}^{3}-907382022282339\,{K}^{2}{a
}^{3}{b}^{4}-259332524337026\,{K}^{2}{a}^{2}{b}^{5}-38984144639072\,{K
}^{2}a{b}^{6}-\\
&2398228288446\,{K}^{2}{b}^{7}-104728017254512\,K{a}^{8}-
256686176681480\,K{a}^{7}b-171608609959724\,K{a}^{6}{b}^{2}+
120617735022502\,K{a}^{5}{b}^{3}+\\
&263942014958862\,K{a}^{4}{b}^{4}+
179209993414076\,K{a}^{3}{b}^{5}+61486917035185\,K{a}^{2}{b}^{6}+
10741146972564\,Ka{b}^{7}+757326482729\,K{b}^{8}-\\
&8161904144976\,{a}^{9
}-5562030977756\,{a}^{8}b+21552564643888\,{a}^{7}{b}^{2}+
41202471901524\,{a}^{6}{b}^{3}+29225205676557\,{a}^{5}{b}^{4}+
\end{aligned}
\end{scriptsize}
$$
$$
\begin{scriptsize}
\begin{aligned}
\qquad &7351641798664\,{a}^{4}{b}^{5}-1939375689153\,{a}^{3}{b}^{6}-
1702959978554\,{a}^{2}{b}^{7}-400742340479\,a{b}^{8}-32815542186\,{b}^
{9}-189026579872\,{K}^{7}a-\\
&58212264218\,{K}^{7}b+2523237000760\,{K}^{6
}{a}^{2}+2531681866046\,{K}^{6}ab+503734666884\,{K}^{6}{b}^{2}+
23979158968552\,{K}^{5}{a}^{3}+\\
&27296320495436\,{K}^{5}{a}^{2}b+
9599362694476\,{K}^{5}a{b}^{2}+1090294195676\,{K}^{5}{b}^{3}+
42538471201464\,{K}^{4}{a}^{4}+42074693837002\,{K}^{4}{a}^{3}b+\\
&
8097446989341\,{K}^{4}{a}^{2}{b}^{2}-2303911716257\,{K}^{4}a{b}^{3}-
605113652267\,{K}^{4}{b}^{4}-41087915794552\,{K}^{3}{a}^{5}-
149506115754198\,{K}^{3}{a}^{4}b-\\
&170582790113698\,{K}^{3}{a}^{3}{b}^{2
}-86134927287435\,{K}^{3}{a}^{2}{b}^{3}-20074062671374\,{K}^{3}a{b}^{4
}-1767653846566\,{K}^{3}{b}^{5}-108016462309312\,{K}^{2}{a}^{6}-\\
&
285840552549800\,{K}^{2}{a}^{5}b-293480390259879\,{K}^{2}{a}^{4}{b}^{2
}-148765538254888\,{K}^{2}{a}^{3}{b}^{3}-38459447679506\,{K}^{2}{a}^{2
}{b}^{4}-\\
&4583930172025\,{K}^{2}a{b}^{5}-173335998053\,{K}^{2}{b}^{6}+
13595184706384\,K{a}^{7}+86095537533340\,K{a}^{6}b+183958941853768\,K{
a}^{5}{b}^{2}+\\
&191662558936560\,K{a}^{4}{b}^{3}+110034868933364\,K{a}^{
3}{b}^{4}+35599357420344\,K{a}^{2}{b}^{5}+6082999069932\,Ka{b}^{6}+
426508476718\,K{b}^{7}+\\
&5468044164048\,{a}^{8}+15447633995764\,{a}^{7}b
+16182033399812\,{a}^{6}{b}^{2}+3862708193224\,{a}^{5}{b}^{3}-
6425359506547\,{a}^{4}{b}^{4}-\\
&6441019453773\,{a}^{3}{b}^{5}-
2600347442170\,{a}^{2}{b}^{6}-504207547600\,a{b}^{7}-38349728870\,{b}^
{8}+122724961740\,{K}^{6}a+35807472855\,{K}^{6}b-\\
&75489197492\,{K}^{5}{
a}^{2}-264558877753\,{K}^{5}ab-65041831722\,{K}^{5}{b}^{2}-
4710070971404\,{K}^{4}{a}^{3}-6241553418239\,{K}^{4}{a}^{2}b-\\
&
2565957790031\,{K}^{4}a{b}^{2}-333289673685\,{K}^{4}{b}^{3}-
12306446883500\,{K}^{3}{a}^{4}-18327404147223\,{K}^{3}{a}^{3}b-
9687819156084\,{K}^{3}{a}^{2}{b}^{2}-\\
&2143749637690\,{K}^{3}a{b}^{3}-
170087002944\,{K}^{3}{b}^{4}-2634017494480\,{K}^{2}{a}^{5}+
2561896630684\,{K}^{2}{a}^{4}b+9176598422829\,{K}^{2}{a}^{3}{b}^{2}+\\
&
6751036270532\,{K}^{2}{a}^{2}{b}^{3}+1971853517875\,{K}^{2}a{b}^{4}+
203243158989\,{K}^{2}{b}^{5}+11711672619776\,K{a}^{6}+35681672761880\,
K{a}^{5}b+\\
&43305441108728\,K{a}^{4}{b}^{2}+27011547227561\,K{a}^{3}{b}^
{3}+9167796339800\,K{a}^{2}{b}^{4}+1609451438003\,Ka{b}^{5}+
114619208064\,K{b}^{6}-\\
&210679824496\,{a}^{7}-2533427160948\,{a}^{6}b-
6829946224908\,{a}^{5}{b}^{2}-8226600450528\,{a}^{4}{b}^{3}-
5284628337933\,{a}^{3}{b}^{4}-1882878507327\,{a}^{2}{b}^{5}\\
&-350582250254\,a{b}^{6}-26481178125\,{b}^{7}-37965515184\,{K}^{5}a-
10872581466\,{K}^{5}b-126416626224\,{K}^{4}{a}^{2}-63613213824\,{K}^{4
}ab-\\
&8534936874\,{K}^{4}{b}^{2}+346727531568\,{K}^{3}{a}^{3}+
569883646590\,{K}^{3}{a}^{2}b+273872487600\,{K}^{3}a{b}^{2}+
39274316289\,{K}^{3}{b}^{3}+1619164542960\,{K}^{2}{a}^{4}\\
&+2844394737888\,{K}^{2}{a}^{3}b+1824341930406\,{K}^{2}{a}^{2}{b}^{2}+
505411593921\,{K}^{2}a{b}^{3}+51067461459\,{K}^{2}{b}^{4}+
1188251075328\,K{a}^{5}+\\
&2250005110284\,K{a}^{4}b+1617213822096\,K{a}^{
3}{b}^{2}+539394318024\,K{a}^{2}{b}^{3}+80618618082\,Ka{b}^{4}+
4110544155\,K{b}^{5}-518566234944\,{a}^{6}-\\
&1768196637240\,{a}^{5}b-
2400596482644\,{a}^{4}{b}^{2}-1677971587536\,{a}^{3}{b}^{3}-
639613818768\,{a}^{2}{b}^{4}-126121969059\,a{b}^{5}-10020245607\,{b}^{
6}+\\
&6053534640\,{K}^{4}a+1714260060\,{K}^{4}b+30872125872\,{K}^{3}{a}^{
2}+21702774420\,{K}^{3}ab+3710033460\,{K}^{3}{b}^{2}+28232631504\,{K}^
{2}{a}^{3}+\\
&24646868820\,{K}^{2}{a}^{2}b+6153974280\,{K}^{2}a{b}^{2}+
460361034\,{K}^{2}{b}^{3}-63821875968\,K{a}^{4}-130420080576\,K{a}^{3}
b-96091123176\,K{a}^{2}{b}^{2}-\\
&30122952084\,Ka{b}^{3}-3352338072\,K{b}
^{4}-90453825216\,{a}^{5}-213233602032\,{a}^{4}b-198833168160\,{a}^{3}
{b}^{2}-91402582068\,{a}^{2}{b}^{3}-\\
&20627686584\,a{b}^{4}-1816925706\,
{b}^{5}-400347360\,{K}^{3}a-112332420\,{K}^{3}b-2402084160\,{K}^{2}{a}
^{2}-1875036600\,{K}^{2}ab-336997260\,{K}^{2}{b}^{2}-\\
&4804168320\,K{a}^
{3}-6152157360\,K{a}^{2}b-2549031120\,Ka{b}^{2}-336997260\,K{b}^{3}-
3202778880\,{a}^{4}-5702827680\,{a}^{3}b-3750073200\,{a}^{2}{b}^{2}-\\
&
1074341880\,a{b}^{3}-112332420\,{b}^{4}.
\end{aligned}
\end{scriptsize}
$$

\end{document}